\documentclass[11pt,a4paper]{article}

\usepackage{soul}
\usepackage{tikz}
\usepackage{amsthm}
\usepackage{amsmath}
\usepackage{amssymb}
\usepackage{mathrsfs}

\usepackage{geometry}
\usepackage{graphicx}
\usepackage{amsfonts}
\usepackage{epstopdf}
\usepackage{placeins}
\usepackage{enumerate}
\usepackage{enumitem}
\usepackage{subcaption}
\usepackage{microtype}
\usepackage[hidelinks]{hyperref} 
\usepackage[normalem]{ulem} 

\usetikzlibrary{calc,arrows.meta, positioning}

\numberwithin{equation}{section}

\allowdisplaybreaks

\newcommand{\esup}{\operatorname*{ess\,sup}}
\newcommand{\einf}{\operatorname*{ess\,inf}}

\newcommand{\Rmnum}[1]{\uppercase\expandafter{\romannumeral#1}} 

\def\Xint#1{\mathchoice
{\XXint\displaystyle\textstyle{#1}}%
{\XXint\textstyle\scriptstyle{#1}}%
{\XXint\scriptstyle\scriptscriptstyle{#1}}%
{\XXint\scriptscriptstyle\scriptscriptstyle{#1}}%
\!\int}
\def\XXint#1#2#3{{\setbox0=\hbox{$#1{#2#3}{\int}$ }
\vcenter{\hbox{$#2#3$ }}\kern-.6\wd0}}

\def\dashint{\Xint-}

\theoremstyle{plain}
\newtheorem{theorem}{Theorem}[section]
\newtheorem{proposition}[theorem]{Proposition}
\newtheorem{lemma}[theorem]{Lemma}
\newtheorem{corollary}[theorem]{Corollary}

\theoremstyle{definition}
\newtheorem{remark}[theorem]{Remark}

\renewcommand{\thefootnote}{}

\AtEndDocument{
\bigskip{\footnotesize%
   \textsc{Department of Mathematics, Aarhus University, 8000 Aarhus C, Denmark} \par
  \textit{E-mail address}: \texttt{yangmengqh@gmail.com}
}}

\begin{document}

\title{\large{Energy inequalities for cutoff functions of $p$-energies on metric measure spaces}}
\author{Meng Yang}
\date{}

\maketitle

\abstract{For $p>1$, and for a $p$-energy on a volume doubling metric measure space, we provide several geometric and functional conditions for the validity of the cutoff Sobolev inequality. In particular, we prove that the elliptic Harnack inequality, two-sided capacity bounds, and some additional geometric and analytic assumptions imply the cutoff Sobolev inequality, \emph{without} assuming the Poincar\'e inequality. Conversely, we show that the Poincar\'e inequality together with the cutoff Sobolev inequality recovers the analytic inputs used in the first implication. Moreover, in a lower-dimensional regime, the Poincar\'e inequality and the capacity upper bound imply the cutoff Sobolev inequality. As an application, we prove that the $p$-energy measure is singular with respect to the Hausdorff measure on the Sierpi\'nski carpet for any $p>1$, thereby resolving a problem posed by Murugan and Shimizu [\textit{Comm. Pure Appl. Math.} \textbf{78} (2025), no. 9, 1523--1608].}

\footnote{\textsl{Date}: \today}
\footnote{The author is very grateful to Ryosuke Shimizu for stimulating discussions and for communicating their recent results relevant to this work during the conference \textit{FRACTALS 8} at Cornell, to Alexander Grigor'yan and Jiaxin Hu for valuable discussions during a visit to Bielefeld, and to Sylvester Eriksson-Bique for sharing his paper \cite{Eri26}, which resolves the resistance conjecture.}

\renewcommand{\thefootnote}{\arabic{footnote}}
\setcounter{footnote}{0}

\section{Introduction}

In the classical PDE theory on $\mathbb{R}^{d}$, the study of regularity for weak solutions often relies on suitable test functions. In many arguments, these test functions are built from the weak solution itself and from appropriate cutoff functions. Consider two concentric balls $B(x_{0},R)\subseteq B(x_{0},R+r)$. By smoothing the distance function, one can readily construct a cutoff function $\phi\in C_{c}^{\infty}(\mathbb{R}^{d})$ such that $0\le \phi\le 1$ on $\mathbb{R}^d$, $\phi=1$ on $B(x_{0},R)$, $\phi=0$ on $\mathbb{R}^{d}\backslash B(x_{0},R+r)$, and $|\nabla \phi|\le 2/r$ on $\mathbb{R}^{d}$. The estimate on $\nabla \phi$ gives the intrinsic scale control needed in the analysis. Analogous cutoff functions with controlled gradients can also be constructed on Riemannian manifolds, where they play a fundamental role.

On general metric measure spaces, particularly on fractals where a gradient operator is typically unavailable, it is natural to ask what kind of estimates for cutoff functions should replace such pointwise gradient bounds. In 2004, Barlow and Bass \cite{BB04} introduced the cutoff Sobolev inequality, which provides energy estimates for cutoff functions associated with Dirichlet forms. This condition has since become a key tool in the study of sub-Gaussian heat kernel estimates. More precisely, the two-sided estimates
\begin{align*}\label{eq_HKbeta}\tag*{HK($\beta$)}
\frac{C_1}{V(x,t^{1/\beta})}\exp\left(-C_2\left(\frac{d(x,y)}{t^{1/\beta}}\right)^{\frac{\beta}{\beta-1}}\right)\le p_t(x,y)\le\frac{C_3}{V(x,t^{1/\beta})}\exp\left(-C_4\left(\frac{d(x,y)}{t^{1/\beta}}\right)^{\frac{\beta}{\beta-1}}\right),\nonumber
\end{align*}
where $\beta\ge2$ is the walk dimension, are closely tied to the conjunction of the volume doubling condition, the Poincar\'e inequality, and the cutoff Sobolev inequality. In fractal examples, one typically has $\beta>2$. This should be compared with the classical case of complete non-compact Riemannian manifolds: it was proved independently by Grigor'yan \cite{Gri92} and by Saloff-Coste \cite{Sal92,Sal95} that the two-sided Gaussian heat kernel estimates
$$\frac{C_1}{V(x,\sqrt{t})}\exp\left(-C_2\frac{d(x,y)^2}{t}\right)\le p_t(x,y)\le\frac{C_3}{V(x,\sqrt{t})}\exp\left(-C_4\frac{d(x,y)^2}{t}\right),$$
are equivalent to the conjunction of the volume doubling condition and the Poincar\'e inequality. This corresponds to the case $\beta=2$, where classical cutoff functions constructed from distance functions yield the required cutoff estimates, and hence the cutoff Sobolev inequality follows trivially. In the genuinely sub-Gaussian case $\beta>2$, the cutoff Sobolev inequality becomes an essential additional ingredient in the study of elliptic and parabolic Harnack inequalities, heat kernel estimates, and stability of analytic estimates. The following theorem summarizes the standard picture, where $\Psi$ is a suitable general scaling function.

\begin{theorem}\label{thm_DF}
For a strongly local regular Dirichlet form on a volume doubling metric measure space, the followings are equivalent.
\begin{enumerate}[label=(\alph*),ref=(\alph*)]
\item\label{item_PHI} The parabolic Harnack inequality $\text{PHI}(\Psi)$.
\item\label{item_HK} The two-sided heat kernel estimates $\text{HK}(\Psi)$.
\item\label{item_EHIcap} The elliptic Harnack inequality $\text{EHI}$ and the two-sided capacity bounds $\text{cap}(\Psi)$.
\item\label{item_PICS} The Poincar\'e inequality $\text{PI}(\Psi)$ and the cutoff Sobolev inequality $\text{CS}(\Psi)$.
\item\label{item_PIcap} The Poincar\'e inequality $\text{PI}(\Psi)$ and the capacity upper bound $\text{cap}(\Psi)_{\leq}$.
\end{enumerate}
\end{theorem}

The equivalences among \ref{item_PHI}, \ref{item_HK}, \ref{item_EHIcap}, and \ref{item_PICS} have long been established; see, for example, \cite{BBK06,GT12,GH14,AB15,GHL15} and the references therein. The remaining implication needed to include \ref{item_PIcap} is the so-called resistance conjecture; see \cite{Bar13, GHL14, Gri22}. This conjecture was very recently resolved by Eriksson-Bique \cite{Eri26}, who proved that, under the volume doubling condition,
\begin{equation}\label{eq_PIcapCS}
\text{PI}(\Psi)+\text{cap}(\Psi)_{\le}\Rightarrow \text{CS}(\Psi).
\end{equation}

The preceding results are formulated in the Dirichlet form framework, which generalizes the classical Dirichlet integral $\int_{\mathbb{R}^{d}}|\nabla f(x)|^{2}\mathrm{d}x$ on $\mathbb{R}^{d}$. For general $p>1$, extending the classical $p$-energy $\int_{\mathbb{R}^{d}}|\nabla f(x)|^{p}\mathrm{d}x$, the study of $p$-energies on fractals and general metric measure spaces, initiated by \cite{HPS04}, has advanced rapidly in recent years; see \cite{CGQ22,Shi24,BC23,MS25,Kig23,CGYZ26,AB25,AES25a}. In this setting, a new parameter $\beta_{p}$, called the $p$-walk dimension, naturally arises in connection with a $p$-energy. In particular, $\beta_{2}$ coincides with the walk dimension $\beta$ in \ref{eq_HKbeta}. It is therefore natural to expect that an appropriate $p$-version of the cutoff Sobolev inequality should provide essential functional tools in the analysis of $p$-energies.

In our previous works \cite{Yan25a,Yan25c,Yan25e}, we introduced a $p$-version of the cutoff Sobolev inequality. As applications, we investigated the singularity of the associated $p$-energy measure with respect to the underlying measure and revealed a dichotomy in the dependence of the $p$-walk dimension $\beta_{p}$ on $p$.

A natural and important question is the following.

\begin{quote}
Can one establish, for \emph{nonlinear} $p$-energies, analogues of the equivalences in Theorem \ref{thm_DF}?
\end{quote}

In the nonlinear setting, a heat kernel is generally not available, and therefore one should formulate the expected equivalences among capacity bounds, Harnack, Poincar\'e, and cutoff Sobolev inequalities. The analogue of the equivalence between \ref{item_PICS} and \ref{item_PIcap} for $p$-energies is also contained in \cite{Eri26}, through the implication (\ref{eq_PIcapCS}). The main purpose of this paper is to make progress toward the analogue of the equivalence between \ref{item_EHIcap} and \ref{item_PICS} for $p$-energies.

The main result of this paper proves the cutoff Sobolev inequality from the elliptic Harnack inequality and the two-sided capacity bounds. More precisely, for a $p$-energy on a volume doubling metric measure space, and under additional geometric and analytic assumptions, but \emph{without} assuming $\text{PI}(\Psi)$, we prove that
\begin{equation}\label{eq_EHIcapCS}
\text{EHI}+\text{cap}(\Psi)\Rightarrow \text{CS}(\Psi).
\end{equation}
See Theorem \ref{thm_main} for the precise statement and the additional assumptions. The proof is based on a Wolff potential estimate for superharmonic functions, established by adapting a technique of Trudinger and Wang \cite{TW02AJM}; in this argument, the elliptic Harnack inequality plays a central role.

We also prove the converse implication from \ref{item_PICS} to \ref{item_EHIcap} for $p$-energies. More precisely, assuming $\text{PI}(\Psi)$ and $\text{CS}(\Psi)$, we obtain $\text{EHI}$, $\text{cap}(\Psi)$, and the analytic inputs used in (\ref{eq_EHIcapCS}); see Theorem \ref{thm_main_conv}. The proof uses the standard Nash--Moser--De Giorgi iteration technique, where the cutoff Sobolev inequality supplies the required cutoff functions.

Moreover, in a lower-dimensional regime, we prove that the Poincar\'e inequality and the capacity upper bound already imply the cutoff Sobolev inequality. In particular, on an Ahlfors $d_h$-regular space, this corresponds to the range $d_h-1<\beta_p\le d_h$. The proof proceeds by first deriving the elliptic Harnack inequality from the Poincar\'e inequality and the capacity upper bound on suitable approximating spaces, and then applying (\ref{eq_EHIcapCS}). This technique for proving the elliptic Harnack inequality, which is restricted to such lower-dimensional regime, goes back to \cite{Hol03}, and was also used recently in \cite{Mur23a,MS25}. Thus, in this regime, we obtain the implication (\ref{eq_PIcapCS}); see Theorem \ref{thm_main_RSVR}. In the case $p=2$, the same method also yields two-sided heat kernel estimates; see Theorem \ref{thm_HK}.

As an application, we prove that the $p$-energy measure is singular with respect to the Hausdorff measure on the Sierpi\'nski carpet for any $p>1$; see Corollary \ref{cor_SC_singular}. This gives a complete answer, in the setting of the Sierpi\'nski carpet, to a problem posed by Murugan and Shimizu \cite[Problem 10.5]{MS25}.

Throughout this paper, $p\in(1,+\infty)$ is fixed. The letters $C$, $C_1$, $C_2$, $C_A$, $C_B$ will always refer to some positive constants and may change at each occurrence. The sign $\asymp$ means that the ratio of the two sides is bounded from above and below by positive constants. The sign $\lesssim$ ($\gtrsim$) means that the LHS is bounded by positive constant times the RHS from above (below). We use $x_+$ to denote the positive part of $x\in\mathbb{R}$, that is, $x_+=\max\{x,0\}$. For $x$, $y\in \mathbb{R}$, let $x\vee y=\max\{x,y\}$ and $x\wedge y=\min\{x,y\}$. We use $\#A$ to denote the cardinality of a set $A$. We use $l(V)$ to denote the family of all real-valued functions on a set $V$.

\section{Statement of main results}

We say that a function $\Phi:[0,+\infty)\to[0,+\infty)$ is doubling if $\Phi$ is a homeomorphism, which implies that $\Phi$ is strictly increasing continuous and $\Phi(0)=0$, and there exists $C_\Phi>1$, called a doubling constant of $\Phi$, such that $\Phi(2r)\le C_\Phi\Phi(r)$ for any $r>0$. Throughout this paper, we always assume that $\Phi, \Psi$ are two doubling functions with doubling constants $C_\Phi, C_\Psi$, respectively, and there exist $\beta_*, \beta^*>0$ with $\beta_*\le\beta^*$ such that
$$\frac{1}{C_\Psi}\left(\frac{R}{r}\right)^{\beta_*}\le \frac{\Psi(R)}{\Psi(r)}\le C_{\Psi}\left(\frac{R}{r}\right)^{\beta^*}\text{ for any }r\le R.$$
Indeed, we can take $\beta^*=\log_2C_\Psi$.

Let $(X,d,m)$ be a complete unbounded metric measure space, that is, $(X,d)$ is a complete unbounded locally compact separable metric space and $m$ is a positive Radon measure on $X$ with full support. Throughout this paper, we always assume that all metric balls are relatively compact. For any $x\in X$ and any $r>0$, denote $B(x,r)=\{y\in X:d(x,y)<r\}$ and $V(x,r)=m(B(x,r))$. If $B=B(x,r)$, then denote $\delta B=B(x,\delta r)$ for any $\delta>0$. Let $\mathcal{B}(X)$ be the family of all Borel measurable subsets of $X$. Let $C(X)$ be the family of all continuous functions on $X$. Let $C_c(X)$ be the family of all continuous functions on $X$ with compact support. Denote $\dashint_A=\frac{1}{m(A)}\int_A$ and $u_A=\dashint_Au\mathrm{d} m$ for any measurable set $A$ with $m(A)\in(0,+\infty)$ and any function $u$ such that the integral $\int_Au\mathrm{d} m$ is well-defined.

Let $\varepsilon>0$. We say that $V\subseteq X$ is an $\varepsilon$-net if for any distinct $x,y\in V$, we have $d(x,y)\ge\varepsilon$, and for any $z\in X$, there exists $x\in V$ such that $d(x,z)<\varepsilon$. Since $(X,d)$ is separable, all $\varepsilon$-nets are countable.

We say that the chain condition \ref{eq_CC} holds if there exists $C_{cc}>0$ such that for any $x,y\in X$ and any integer $n\ge1$, there exists a sequence $\{x_k:0\le k\le n\}$ of points in $X$ with $x_0=x$ and $x_n=y$ such that
\begin{equation*}\label{eq_CC}\tag*{CC}
d(x_k,x_{k-1})\le C_{cc} \frac{d(x,y)}{n}\text{ for any }k=1,\ldots,n.
\end{equation*}
Throughout this paper, we always assume \ref{eq_CC}.

We say that the linearly locally connected condition \hypertarget{eq_LLC}{\text{LLC}} holds if there exists $A_{LLC}>1$ such that for any $x_0\in X$ and any $r>0$, for any $x,y\in B(x_0,r)\backslash B(x_0,r/2)$, there exists a connected compact set $K\subseteq B(x_0,A_{LLC}r)\backslash B(x_0,r/(2A_{LLC}))$ containing $x,y$. See also \cite[3.12]{HK98} or \cite[Page 234]{HKST15}.

We say that the volume doubling condition \ref{eq_VD} holds if there exists $C_{VD}>0$ such that
\begin{equation*}\label{eq_VD}\tag*{VD}
V(x,2r)\le C_{VD}V(x,r)\text{ for any }x\in X,r>0.
\end{equation*}

We say that $(\mathcal{E},\mathcal{F})$ is a $p$-energy on $(X,d,m)$ if $\mathcal{F}$ is a dense subspace of $L^p(X;m)$ and $\mathcal{E}:\mathcal{F}\to[0,+\infty)$ satisfies the following conditions.

\begin{enumerate}[label=(\arabic*)]
\item $\mathcal{E}^{1/p}$ is a semi-norm on $\mathcal{F}$, that is, for any $f,g\in\mathcal{F}$, $c\in\mathbb{R}$, we have $\mathcal{E}(f)\ge0$, $\mathcal{E}(cf)^{1/p}=|c|\mathcal{E}(f)^{1/p}$ and $\mathcal{E}(f+g)^{1/p}\le\mathcal{E}(f)^{1/p}+\mathcal{E}(g)^{1/p}$.
\item (Closed property) $(\mathcal{F},\mathcal{E}(\cdot)^{1/p}+\lVert {\cdot}\rVert_{L^p(X;m)})$ is a Banach space.
\item (Markovian property) For any $\varphi\in C(\mathbb{R})$ with $\varphi(0)=0$ and $|\varphi(t)-\varphi(s)|\le|t-s|$ for any $t,s\in\mathbb{R}$, for any $f\in\mathcal{F}$, we have $\varphi(f)\in\mathcal{F}$ and $\mathcal{E}(\varphi(f))\le\mathcal{E}(f)$.
\item (Regular property) $\mathcal{F}\cap C_c(X)$ is uniformly dense in $C_c(X)$ and $(\mathcal{E}(\cdot)^{1/p}+\lVert {\cdot}\rVert_{L^p(X;m)})$-dense in $\mathcal{F}$.
\item (Strongly local property) For any $f,g\in\mathcal{F}$ with compact support and $g$ constant in an open neighborhood of $\mathrm{supp}(f)$, we have $\mathcal{E}(f+g)=\mathcal{E}(f)+\mathcal{E}(g)$.
\item ($p$-Clarkson's inequality) For any $f,g\in\mathcal{F}$, we have
\begin{equation*}\label{eq_Cla}\tag*{Cla}
\begin{cases}
\mathcal{E}(f+g)+\mathcal{E}(f-g)\ge2 \left(\mathcal{E}(f)^{\frac{1}{p-1}}+\mathcal{E}(g)^{\frac{1}{p-1}}\right)^{p-1}&\text{if }p\in(1,2],\\
\mathcal{E}(f+g)+\mathcal{E}(f-g)\le2 \left(\mathcal{E}(f)^{\frac{1}{p-1}}+\mathcal{E}(g)^{\frac{1}{p-1}}\right)^{p-1}&\text{if }p\in[2,+\infty).\\
\end{cases}
\end{equation*}
\end{enumerate}
Moreover, we also always assume the following condition.
\begin{itemize}
\item ($\mathcal{F}\cap L^\infty(X;m)$ is an algebra) For any $f,g\in\mathcal{F}\cap L^\infty(X;m)$, we have $fg\in\mathcal{F}$ and
\begin{equation*}\label{eq_Alg}\tag*{Alg}
\mathcal{E}(fg)^{{1}/{p}}\le \lVert f\rVert_{L^\infty(X;m)}\mathcal{E}(g)^{1/p}+\lVert g\rVert_{L^\infty(X;m)}\mathcal{E}(f)^{1/p}.
\end{equation*}
\end{itemize}
Denote $\mathcal{E}_\lambda(\cdot)=\mathcal{E}(\cdot)+\lambda \lVert {\cdot}\rVert^p_{L^p(X;m)}$ for any $\lambda>0$. Indeed, a general condition called the generalized $p$-contraction property was introduced in \cite{KS24a}, which implies \ref{eq_Cla}, \ref{eq_Alg}, and holds on a large family of metric measure spaces.

We list some basic properties of $p$-energies as follows. For any $f,g\in\mathcal{F}$, the derivative
$$\mathcal{E}(f;g)=\frac{1}{p}\frac{\mathrm{d}}{\mathrm{d} t}\mathcal{E}(f+tg)|_{t=0}\in\mathbb{R}$$
exists, the map $\mathcal{E}(f;\cdot):\mathcal{F}\to\mathbb{R}$ is linear, $\mathcal{E}(f;f)=\mathcal{E}(f)$. Moreover, for any $f,g\in\mathcal{F}$, for any $a\in\mathbb{R}$, we have
\begin{equation}\label{eq_quasi_strict}
\mathbb{R}\ni t\mapsto\mathcal{E}(f+tg;g)\in\mathbb{R}\text{ is strictly increasing if and only if }\mathcal{E}(g)>0,
\end{equation}
$$\mathcal{E}(af;g)=\mathrm{sgn}(a)|a|^{p-1}\mathcal{E}(f;g),$$
$$|\mathcal{E}(f;g)|\le\mathcal{E}(f)^{(p-1)/p}\mathcal{E}(g)^{1/p}.$$
Moreover, for any $\lambda>0$, all of the above results remain valid with $\mathcal{E}$ replaced by $\mathcal{E}_\lambda$, and for any $f,g\in\mathcal{F}$, we have
$$\mathcal{E}_\lambda(f;g)=\mathcal{E}(f;g)+\lambda\int_X\mathrm{sgn}(f)|f|^{p-1}g\mathrm{d} m.$$
See \cite[Theorem 3.7, Corollary 3.25]{KS24a} for the proofs of these results.

By \cite[Theorem 2.4]{Sas26}, a $p$-energy $(\mathcal{E},\mathcal{F})$ corresponds to a (canonical) $p$-energy measure $\Gamma:\mathcal{F}\times\mathcal{B}(X)\to[0,+\infty)$, $(f,A)\mapsto\Gamma(f)(A)$ satisfying the following conditions.
\begin{enumerate}[label=(\alph*),ref=(\alph*)]
\item\label{item_meas1} For any $f\in\mathcal{F}$, $\Gamma(f)(\cdot)$ is a positive Radon measure on $X$ with $\Gamma(f)(X)=\mathcal{E}(f)$.
\item\label{item_meas2} For any $A\in\mathcal{B}(X)$, $\Gamma(\cdot)(A)^{1/p}$ is a semi-norm on $\mathcal{F}$.
\item\label{item_meas3} For any $f,g\in\mathcal{F}\cap C_c(X)$, $A\in\mathcal{B}(X)$, if $f-g$ is constant on $A$, then $\Gamma(f)(A)=\Gamma(g)(A)$.
\item\label{item_meas4} ($p$-Clarkson's inequality) For any $f,g\in\mathcal{F}$, for any $A\in\mathcal{B}(X)$, we have
\begin{equation*}
\begin{cases}
\Gamma(f+g)(A)+\Gamma(f-g)(A)\ge2 \left(\Gamma(f)(A)^{\frac{1}{p-1}}+\Gamma(g)(A)^{\frac{1}{p-1}}\right)^{p-1}&\text{if }p\in(1,2],\\
\Gamma(f+g)(A)+\Gamma(f-g)(A)\le2 \left(\Gamma(f)(A)^{\frac{1}{p-1}}+\Gamma(g)(A)^{\frac{1}{p-1}}\right)^{p-1}&\text{if }p\in[2,+\infty).\\
\end{cases}
\end{equation*}
\item\label{item_meas5} (Chain rule) For any $f,g\in\mathcal{F}\cap C_c(X)$, for any piecewise $C^1$ functions $\varphi,\psi:\mathbb{R}\to\mathbb{R}$ with $\varphi(0)=\psi(0)=0$, we have
$$\mathrm{d}\Gamma\left(\varphi(f);\psi(g)\right)=\mathrm{sgn}(\varphi'(f))|\varphi'(f)|^{p-1}\psi'(g)\mathrm{d}\Gamma(f;g),$$
where $\Gamma(f;g)$ is a signed measure given by $\Gamma(f;g)=\frac{1}{p}\frac{\mathrm{d}}{\mathrm{d}t}\Gamma(f+tg)|_{t=0}$.
\item\label{item_meas6} (Leibniz rule) For any $f,g,h\in \mathcal{F}\cap C_c(X)$, we have $\mathrm{d}\Gamma(f;gh)=g \mathrm{d}\Gamma(f;h)+h \mathrm{d}\Gamma(f;g)$.
\end{enumerate}
Using the chain rule, we have the following condition.
\begin{itemize}
\item (Strong sub-additivity) For any $f,g\in\mathcal{F}$, we have $f\vee g, f\wedge g\in\mathcal{F}$ and
\begin{equation*}\label{eq_SubAdd}\tag*{SubAdd}
\mathcal{E}(f\vee g)+\mathcal{E}(f\wedge g)\le\mathcal{E}(f)+\mathcal{E}(g).
\end{equation*}
\end{itemize}

Let
\begin{align*}
&\mathcal{F}_{\mathrm{loc}}=\left\{u:
\begin{array}{l}
\text{for any relatively compact open set }U,\\
\text{there exists }u^\#\in \mathcal{F}\text{ such that }u=u^\#\text{ }m\text{-a.e. in }U
\end{array}
\right\}.
\end{align*}
For any $u\in \mathcal{F}_{\mathrm{loc}}$, let $\Gamma(u)|_U=\Gamma(u^\#)|_U$, where $u^\#, U$ are given as above, then $\Gamma(u)$ is a well-defined positive Radon measure on $X$, as follows from \ref{item_meas3} and \ref{item_meas6} together with an argument similar to that in \cite[Corollary 3.2.1]{FOT11}.

Let $U,V$ be two open subsets of $X$ satisfying $U\subseteq\overline{U}\subseteq V$. We say that $\phi\in\mathcal{F}$ is a cutoff function for $U\subseteq V$ if $0\le\phi\le1$ in $X$, $\phi=1$ in an open neighborhood of $\overline{U}$ and $\mathrm{supp}(\phi)\subseteq V$, where $\mathrm{supp}(f)$ refers to the support of the measure of $|f|\mathrm{d} m$ for any given function $f$.

We say that the cutoff Sobolev inequality \ref{eq_CS} holds if there exist $C_{1},C_{2}>0$, $A_{S}>1$ such that for any ball $B(x,r)$, there exists a cutoff function $\phi\in\mathcal{F}$ for $B(x,r)\subseteq B(x,A_Sr)$ such that for any $f\in\mathcal{F}$, we have
\begin{equation*}\label{eq_CS}\tag*{CS($\Psi$)}
\int_{B(x,A_{S}r)}|\widetilde{f}|^p\mathrm{d}\Gamma(\phi)\le C_{1}\int_{B(x,A_{S}r)}\mathrm{d}\Gamma(f)+\frac{C_{2}}{\Psi(r)}\int_{B(x,A_{S}r)}|f|^p\mathrm{d} m,
\end{equation*}
where $\widetilde{f}$ is a quasi-continuous modification of $f$, such that $\widetilde{f}$ is uniquely determined $\Gamma(\phi)$-a.e. in $X$, see \cite[Section 8]{Yan25a} for more details. For $\beta_p>0$, we say that the cutoff Sobolev inequality \hypertarget{eq_CSbeta}{CS($\beta_p$)} holds if \ref{eq_CS} holds with $\Psi:r\mapsto r^{\beta_p}$.

Let $A_1,A_2\in\mathcal{B}(X)$. We define the capacity between $A_1,A_2$ as
\begin{align*}
&\mathrm{cap}(A_1,A_2)=\inf\left\{\mathcal{E}(\varphi):\varphi\in\mathcal{F},
\begin{array}{l}
\varphi=1\text{ in an open neighborhood of }A_1,\\
\varphi=0\text{ in an open neighborhood of }A_2
\end{array}
\right\},
\end{align*}
here we use the convention that $\inf\emptyset=+\infty$.

We say that the two-sided capacity bounds \hypertarget{eq_cap}{$\text{cap}(\Psi)$} hold if both the capacity upper bound \ref{eq_ucap} and the capacity lower bound \ref{eq_lcap} hold as follows. There exist $C_{cap}>0$, $A_{cap}>1$ such that for any ball $B(x,r)$, we have
\begin{align*}
\mathrm{cap}\left(B(x,r),X\backslash B(x,A_{cap}r)\right)&\le C_{cap} \frac{V(x,r)}{\Psi(r)},\label{eq_ucap}\tag*{$\text{cap}(\Psi)_{\le}$}\\
\mathrm{cap}\left(B(x,r),X\backslash B(x,A_{cap}r)\right)&\ge \frac{1}{C_{cap}} \frac{V(x,r)}{\Psi(r)}.\label{eq_lcap}\tag*{$\text{cap}(\Psi)_{\ge}$}
\end{align*}
For $\beta_p>0$, we say that \hypertarget{eq_capbeta}{$\text{cap}(\beta_p)$} (resp. \hypertarget{eq_ucapbeta}{$\text{cap}(\beta_p)_{\le}$}, \hypertarget{eq_lcapbeta}{$\text{cap}(\beta_p)_{\ge}$}) holds if \hyperlink{eq_cap}{$\text{cap}(\Psi)$} (resp. \ref{eq_ucap}, \ref{eq_lcap}) holds with $\Psi:r\mapsto r^{\beta_p}$. Under \ref{eq_VD}, by taking $f\equiv1$ in $B(x,A_Sr)$, it is easy to see that \ref{eq_CS} (resp. \hyperlink{eq_CSbeta}{$\text{CS}(\beta_p)$}) implies \ref{eq_ucap} (resp. \hyperlink{eq_ucapbeta}{$\text{cap}(\beta_p)_{\le}$}).

Let $U$ be an open subset of $X$. Let
$$\mathcal{F}(U)=\text{the }\mathcal{E}_1\text{-closure of }\mathcal{F}\cap C_c(U).$$
We say that $u\in\mathcal{F}$ is harmonic in $U$ if $\mathcal{E}(u;v)=0$ for any $v\in\mathcal{F}\cap C_c(U)$, denoted by $-\Delta_pu=0$ in $U$. We say that $u\in\mathcal{F}$ is superharmonic in $U$ (resp. subharmonic in $U$) if $\mathcal{E}(u;v)\ge0$ (resp. $\mathcal{E}(u;v)\le0$) for any non-negative $v\in\mathcal{F}\cap C_c(U)$, denoted by $-\Delta_pu\ge0$ in $U$ (resp. $-\Delta_pu\le0$ in $U$). Moreover, if $U$ is bounded, then by \cite[Lemma 6.2]{Yan25c}, we have
$$\mathcal{F}(U)=\left\{u\in\mathcal{F}:\widetilde{u}=0\text{ q.e. on }X\backslash U\right\}.$$
Assuming that $U$ is bounded, then the above equality for harmonic functions, as well as the corresponding inequalities for superharmonic and subharmonic functions, also hold for all $v\in\mathcal{F}(U)$.

We say that the bottom spectrum positivity condition \ref{eq_BSP} holds if for any bounded open subset $U$ of $X$, we have
\begin{equation*}\label{eq_BSP}\tag*{\text{BSP}}
\lambda_1(U)=\inf \left\{\frac{\mathcal{E}(u)}{\lVert u\rVert_{L^p(X;m)}^p}:u\in \mathcal{F}(U)\backslash \{0\}\right\}>0.
\end{equation*}

We say that the lower semi-continuous condition \hypertarget{eq_LSC}{\text{LSC}} holds if for any bounded open subset $U$ of $X$, for any $u\in \mathcal{F}$ which is bounded from below superharmonic in $U$, we have $u$ has a modification which is lower semi-continuous in $U$.

We say that the elliptic Harnack inequality \ref{eq_EHI} holds if there exist $C_H>0$, $A_H>1$ such that for any ball $B(x,r)$, for any $u\in\mathcal{F}$ which is non-negative harmonic in $B(x,A_Hr)$, we have
\begin{equation*}\label{eq_EHI}\tag*{\text{EHI}}
\esup_{B(x,r)}u\le C_H\einf_{B(x,r)}u.
\end{equation*}

We say that the annulus elliptic Harnack inequality \ref{eq_EHIann} holds if there exist $A_1$, $A_2$, $A_3>1$ with $A_1<A_2<A_3$ and $C_H>0$ such that the following holds: for any ball $B(x,r)$ and any $u\in \mathcal{F}$ which is non-negative superharmonic in $B(x,A_3r)$ and harmonic in $B(x,A_3r)\backslash \overline{B(x,r)}$, we have
\begin{equation*}\label{eq_EHIann}\tag*{$\text{EHI}_{\text{ann}}$}
\esup_{B(x,A_2r)\backslash B(x,A_1r)}u\le C_H\einf_{B(x,A_2r)\backslash B(x,A_1r)}u.
\end{equation*}
See also \cite[Theorem 3.1]{Bla01} and \cite[LEMMA 6.3]{GS05} for similar results. However, in addition to harmonicity in the annulus, we also require superharmonicity in the entire ball. By a standard chaining argument, as in these references, it is easy to see that
\begin{center}
\ref{eq_VD} + \hyperlink{eq_LLC}{\text{LLC}} + \ref{eq_EHI} $\Rightarrow$ \ref{eq_EHIann}.
\end{center}

The main result of this paper is as follows.

\begin{theorem}\label{thm_main}
Assume \ref{eq_VD}, \hyperlink{eq_LSC}{\text{LSC}}, \ref{eq_BSP}, \ref{eq_EHI}, \hyperlink{eq_cap}{$\text{cap}(\Psi)$}. If either
\begin{enumerate}[label=(\alph*),ref=(\alph*)]
\item\label{cond_main_LLC} \hyperlink{eq_LLC}{\text{LLC}}, or
\item\label{cond_main_EHIann} \ref{eq_EHIann},
\end{enumerate}
holds, then \ref{eq_CS} holds.
\end{theorem}

\begin{remark}
\begin{enumerate}[label=(\arabic*)]
\item Since \ref{eq_VD}, \hyperlink{eq_LLC}{\text{LLC}}, \ref{eq_EHI} imply \ref{eq_EHIann}, it suffices here to treat case \ref{cond_main_EHIann}.
\item For the case $p=2$, it is well-known in the theory of heat kernel estimates that, under \ref{eq_VD}, the conjunction of \ref{eq_EHI} and \hyperlink{eq_cap}{$\text{cap}(\Psi)$} already implies \ref{eq_CS}; see \cite[THEOREM 1.2]{GHL15}. However, for general $p>1$, due to the \emph{nonlinearity} and the limitation of our techniques, \hyperlink{eq_LSC}{\text{LSC}}, \ref{eq_BSP}, and either \hyperlink{eq_LLC}{\text{LLC}} or \ref{eq_EHIann} are additionally required.
\end{enumerate}
\end{remark}

In a complementary direction, the condition \ref{eq_CS}, together with only \ref{eq_VD} and the Poincar\'e inequality \ref{eq_PI} introduced below, implies the remaining functional inequalities. We now recall the latter condition.

We say that the Poincar\'e inequality \ref{eq_PI} holds if there exist $C_{PI}>0$, $A_{PI}\ge1$ such that for any ball $B(x,r)$, for any $f\in\mathcal{F}$, we have
\begin{equation*}\label{eq_PI}\tag*{PI($\Psi$)}
\int_{B(x,r)}\lvert f-f_{B(x,r)}\rvert^p\mathrm{d} m\le C_{PI}\Psi(r)\int_{{B(x,A_{PI}r)}}\mathrm{d}\Gamma(f).
\end{equation*}
For $\beta_p>0$, we say that the Poincar\'e inequality \hypertarget{eq_PIbeta}{PI($\beta_p$)} holds if \ref{eq_PI} holds with $\Psi:r\mapsto r^{\beta_p}$. Under \ref{eq_VD} and \ref{eq_PI}, both \ref{eq_BSP} and \ref{eq_lcap} hold; these follow from \cite[Lemma 4.1]{Yan25c} and \cite[Lemma 5.1 (a)]{BB04}, respectively.

With these preliminaries in place, we can now state the complementary result.

\begin{theorem}\label{thm_main_conv}
Assume \ref{eq_VD}, \ref{eq_PI}, \ref{eq_CS}. Then \hyperlink{eq_LSC}{\text{LSC}}, \ref{eq_BSP}, \ref{eq_EHI}, \ref{eq_EHIann}, \hyperlink{eq_cap}{$\text{cap}(\Psi)$} hold.
\end{theorem}

Moreover, in a lower-dimensional regime, \ref{eq_PI} and \ref{eq_ucap} already imply \ref{eq_CS}. To formulate this precisely, we need the following conditions.

We say that the volume regular condition \ref{eq_VPhi} holds if there exists $C_{VR}>0$ such that
\begin{equation*}\label{eq_VPhi}\tag*{V($\Phi$)}
\frac{1}{C_{VR}}\Phi(r)\le V(x,r)\le C_{VR}\Phi(r)\text{ for any }x\in X,r>0.
\end{equation*}
For $d_h>0$, we say that the Ahlfors regular condition \hypertarget{eq_Vdh}{V($d_h$)} holds if \ref{eq_VPhi} holds with $\Phi:r\mapsto r^{d_h}$.

We say that the fast volume regular condition \ref{eq_FVR} holds if \ref{eq_VPhi} holds and there exists $C_{FVR}>0$ such that
\begin{equation*}\label{eq_FVR}\tag*{FVR($\Phi,\Psi$)}
\frac{\Psi(R)}{\Psi(r)}\le C_{FVR}\frac{\Phi(R)}{\Phi(r)}\text{ for any }r\le R.
\end{equation*}
For $d_h,\beta_p>0$, we say that the fast volume regular condition \hypertarget{eq_FVRdb}{FVR($d_h,\beta_p$)} holds if \hyperlink{eq_Vdh}{V($d_h$)} holds and $d_h\ge\beta_p$.

We say that the relatively slow volume regular condition \ref{eq_RSVR} holds if there exist $\tau<1$, $C_{RSVR}>0$ such that
\begin{equation*}\label{eq_RSVR}\tag*{RSVR($\Phi,\Psi$)}
\frac{1}{C_{RSVR}}\left(\frac{r}{R}\right)^{\tau}\frac{\Phi(R)}{\Phi(r)}\le \frac{\Psi(R)}{\Psi(r)}\text{ for any }r\le R.
\end{equation*}
For $d_h,\beta_p>0$, we say that the relatively slow volume regular condition \hypertarget{eq_RSVRdb}{RSVR($d_h,\beta_p$)} holds if \hyperlink{eq_Vdh}{V($d_h$)} holds and $\beta_p>d_h-1$. If \ref{eq_RSVR} holds with $\tau<0$, then the slow volume regular condition $\mathrm{SVR}(\Phi,\Psi)$ introduced in \cite{Yan25c} is satisfied. This justifies the terminology ``relatively" slow.

If both \ref{eq_FVR} and \ref{eq_RSVR} hold, then we have the following result.

\begin{theorem}\label{thm_main_RSVR}
Assume \ref{eq_FVR}, \ref{eq_RSVR}, \ref{eq_PI}, \ref{eq_ucap}. Then \ref{eq_CS} holds.
\end{theorem}

In the case $p=2$, assuming \ref{eq_RSVR} holds, we obtain the following result.

\begin{theorem}\label{thm_HK}
Let $p=2$ and let $(\mathcal{E},\mathcal{F})$ be a strongly local regular Dirichlet form on $L^2(X;m)$. Assume \ref{eq_RSVR}, \ref{eq_PI}, \ref{eq_ucap}. Then \ref{eq_EHI}, \ref{eq_CS} hold, hence the associated heat kernel $p_t(x,y)$ satisfies the following two-sided estimates. There exist $C_1,C_2,C_3,C_4>0$ such that for any $x,y\in X$, for any $t>0$, we have
\begin{equation*}\label{eq_HK}\tag*{\text{HK}($\Psi$)}
\frac{C_1}{V\left(x,\Psi^{-1}(t)\right)}\exp\left(-\Upsilon\left(C_2d(x,y),t\right)\right)\le p_t(x,y)\le\frac{C_3}{V\left(x,\Psi^{-1}(t)\right)}\exp\left(-\Upsilon\left(C_4d(x,y),t\right)\right),
\end{equation*}
where
$$\Upsilon(R,t)=\sup_{s\in(0,+\infty)}\left(\frac{R}{s}-\frac{t}{\Psi(s)}\right).$$
\end{theorem}

\begin{remark}
\begin{enumerate}[label=(\arabic*)]
\item To prove Theorems \ref{thm_main_RSVR} and \ref{thm_HK}, we use the technique from the proofs of \cite[Theorem 1.1]{Mur23a} and \cite[Theorem 5.4]{MS25} to establish the \ref{eq_EHI} condition on certain spaces that approximate $X$, under \ref{eq_RSVR}. See also \cite[Theorem 4.3]{Hol03}, where this technique was applied to prove \ref{eq_EHI} in the classical setting. Notably, both the \hyperlink{eq_LSC}{\text{LSC}} and \ref{eq_BSP} conditions hold trivially on these approximating spaces. In the proof of Theorem \ref{thm_main_RSVR}, the condition \ref{eq_FVR} is used to establish \hyperlink{eq_LLC}{\text{LLC}} on these approximating spaces, so that Theorem \ref{thm_main} can be applied. By contrast, in the proof of Theorem \ref{thm_HK}, the \ref{eq_EHI} condition together with two-sided capacity bounds is sufficient to obtain the CS condition on the approximating spaces, and hence on $X$.
\item The above two results are related to the so-called resistance conjecture; see \cite[Remark 3.17(1)]{Bar13}, \cite[Conjecture 4.16]{GHL14}, and \cite[Conjecture on Page 155]{Gri22}. This conjecture asserts that, under \ref{eq_VD}, the conjunction of \ref{eq_PI} and \ref{eq_ucap} implies \ref{eq_HK}. Under \ref{eq_VD}, the condition \ref{eq_HK} is previously known to be equivalent to the conjunction of \ref{eq_PI} and \ref{eq_CS}; see \cite[Theorem 1.2]{GHL15}. This conjecture was very recently solved by Eriksson-Bique \cite{Eri26}, where it was shown that, under \ref{eq_VD}, \ref{eq_PI} and \ref{eq_ucap} imply \ref{eq_CS}. Here, we prove related results by first establishing the \ref{eq_EHI} condition in Theorem \ref{thm_main}, which may be of independent interest.
\end{enumerate}
\end{remark}

Building on our previous results in \cite{Yan25a,Yan25c}, we obtain the following result for homogeneous metric measure spaces.

\begin{corollary}\label{cor_main_homo}
Assume \hyperlink{eq_Vdh}{V($d_h$)}, $\beta_p>d_h-1$, \hyperlink{eq_PIbeta}{PI($\beta_p$)}, \hyperlink{eq_capbeta}{cap($\beta_p$)$_\le$}. Then \hyperlink{eq_CSbeta}{CS($\beta_p$)} holds.
\end{corollary}

\begin{proof}
By \cite[Proposition 2.1]{Yan25a}, the conjunction of \hyperlink{eq_Vdh}{V($d_h$)}, \hyperlink{eq_PIbeta}{PI($\beta_p$)} and \hyperlink{eq_capbeta}{cap($\beta_p$)$_\le$} implies $p\le\beta_p\le d_h+(p-1)$. If $\beta_p>d_h$, then \hyperlink{eq_CSbeta}{CS($\beta_p$)} follows from \cite[Theorem 2.3]{Yan25c}. If $\beta_p\le d_h$, then \hyperlink{eq_CSbeta}{CS($\beta_p$)} follows from Theorem \ref{thm_main_RSVR}.
\end{proof}

As a direct consequence of the above result, we establish that, on the Sierpi\'nski carpet, for \emph{all} $p>1$, the cutoff Sobolev inequality holds and the associated $p$-energy measure is singular with respect to the Hausdorff measure, as stated below. This sharpens our earlier result \cite[Corollary 2.5]{Yan25c} and provides a complete answer to \cite[Problem 10.5]{MS25} in the special case of the Sierpi\'nski carpet.

\begin{corollary}\label{cor_SC_singular}
Let $X$ be the Sierpi\'nski carpet, $d$ the Euclidean metric, $d_h=\frac{\log8}{\log3}$ the Hausdorff dimension, and $m$ the $d_h$-dimensional Hausdorff measure. For any $p>1$, let $(\mathcal{E},\mathcal{F})$ and $\Gamma$ be the $p$-energy and the associated $p$-energy measure as given in \cite{MS25}. Then \hyperlink{eq_PIbeta}{PI($\beta_p$)}, \hyperlink{eq_capbeta}{cap($\beta_p$)$_\le$}, \hyperlink{eq_CSbeta}{CS($\beta_p$)}, \ref{eq_EHI} all hold. In particular, for any $f\in\mathcal{F}$, we have $\Gamma(f)\perp m$.
\end{corollary}

\begin{proof}
For any $p>1$, by \cite[Proposition 8.8 (i)]{MS25}, we have $\beta_p\ge p>1>d_h-1$, by \cite[Proposition 6.21, Theorem 8.21]{MS25}, \hyperlink{eq_capbeta}{cap($\beta_p$)$_\le$}, \hyperlink{eq_PIbeta}{PI($\beta_p$)} hold. By Corollary \ref{cor_main_homo} and Theorem \ref{thm_main_conv}, \hyperlink{eq_CSbeta}{CS($\beta_p$)}, \ref{eq_EHI} hold. By \cite[Theorem 2.1]{Yan25c} and either \cite[Theorem 9.8]{KS24a} or \cite[Corollary 2.2]{Yan25e}, we have $\Gamma(f)\perp m$ for any $f\in\mathcal{F}$.
\end{proof}

This paper is organized as follows. In Section \ref{sec_CS}, we prove Theorem \ref{thm_main}. In Section \ref{sec_EHI}, we prove Theorem \ref{thm_main_conv}. In Section \ref{sec_RSVR}, we prove Theorems \ref{thm_main_RSVR} and \ref{thm_HK}. In Section \ref{sec_RBA}, we provide the proof of a result used in Section \ref{sec_CS}.

\section{Proof of Theorem \ref{thm_main}}\label{sec_CS}

We only need to prove the following result.

\begin{proposition}\label{prop_main_CS}
Assume \ref{eq_VD}, \hyperlink{eq_LSC}{\text{LSC}}, \ref{eq_BSP}, \ref{eq_EHI}, \ref{eq_EHIann}, \hyperlink{eq_cap}{$\text{cap}(\Psi)$}. Then \ref{eq_CS} holds.
\end{proposition}

We establish several preparatory results concerning harmonic and superharmonic functions. To begin with, we have the solvability of the boundary value problem for harmonic functions as follows.

\begin{proposition}[{\cite[Proposition 4.2]{Yan25c}}]\label{prop_harm_exist}
Assume \ref{eq_VD}, \ref{eq_BSP}. Let $U$ be a bounded open set and $u\in\mathcal{F}$. Then there exists a unique function $h\in\mathcal{F}$ such that $h$ is harmonic in $U$ and $\widetilde{h}=\widetilde{u}$ q.e. on $X\backslash U$. We denote this function by $H^Uu$. Indeed, $H^Uu$ is the unique function $v\in \mathcal{F}$ satisfying $\mathcal{E}(v)=\inf \{\mathcal{E}(w):w\in \mathcal{F},\widetilde{w}=\widetilde{u}\text{ q.e. on }X\backslash U\}$. Moreover, if $0\le u\le M$ $m$-a.e. in $X$, where $M>0$ is some constant, then $0\le H^Uu\le M$ $m$-a.e. in $X$.
\end{proposition}

Every superharmonic function induces a positive Radon measure, known as the Riesz measure, as follows.

\begin{proposition}\label{prop_super_meas}
Let $U$ be an open set and $u\in\mathcal{F}$ superharmonic in $U$. Then there exists a positive Radon measure $\mu[u]$ on $U$ such that $\mathcal{E}(u;v)=\int_Uv\mathrm{d}\mu[u]$ for any $v\in\mathcal{F}\cap C_c(U)$, denoted by $-\Delta_pu=\mu[u]$ in $U$. Moreover, let $u_1$, $u_2\in\mathcal{F}$ be superharmonic in $U$, and $W$, $V$ two open sets satisfying that $W\subseteq \overline{W}\subseteq V\subseteq U$ and $u_1|_{V\backslash\overline{W}}=u_2|_{V\backslash\overline{W}}$, then $\mu[u_1](V)=\mu[u_2](V)$, that is, if $u_1=u_2$ near $V$, then $\mu[u_1](V)=\mu[u_2](V)$.
\end{proposition}

\begin{proof}
The proof of the existence of $\mu[u]$ is essentially the same as that of \cite[Proposition 8.9, ``(3)$\Rightarrow$(1)"]{Yan25a}, with $\mathcal{E}_1$ replaced with $\mathcal{E}$, and is therefore omitted. Moreover, let $u_1$, $u_2$, $W$, $V$ be given as above. Let $\{W_n\}_{n\ge0}$ be a sequence of open sets satisfying that $W_0=W$ and $W_n\subseteq\overline{W}_n\subseteq W_{n+1}\uparrow V$, then for any $n\ge0$, there exists $\varphi_n\in\mathcal{F}\cap C_c(W_{n+1})$ satisfying that $0\le\varphi_n\le1$ in $W_{n+1}$ and $\varphi_n=1$ on $\overline{W}_n$. By replacing $\varphi_n$ with $\varphi_0\vee\ldots\vee\varphi_n$, we may assume that $\{\varphi_n\}_{n\ge0}$ is increasing, then it follows obviously that $\varphi_n\uparrow1_V$. By the monotone convergence theorem and the strongly local property, we have
\begin{align*}
&\mu[u_1](V)=\lim_{n\to+\infty}\int_U\varphi_n\mathrm{d}\mu[u_1]=\lim_{n\to+\infty}\mathcal{E}(u_1;\varphi_n)\\
&=\lim_{n\to+\infty}\mathcal{E}(u_2;\varphi_n)=\lim_{n\to+\infty}\int_U\varphi_n\mathrm{d}\mu[u_2]=\mu[u_2](V).
\end{align*}
\end{proof}

The following result serves as a useful tool for constructing new superharmonic functions. The proof follows the same argument as in the proof of \cite[Lemma 7.13]{BB11}, with the integration of the upper gradient replaced by that of the energy measure and by using the strong local property, and is therefore omitted.

\begin{lemma}[Pasting lemma]\label{lem_pasting}
Let $U_1\subseteq U_2$ be two bounded open sets, $u_1\in\mathcal{F}$ superharmonic in $U_1$, $u_2\in\mathcal{F}$ superharmonic in $U_2$, and $\widetilde{u}_1=\widetilde{u}_2$ q.e. on $U_2\backslash U_1$. Then $u_1\wedge u_2\in\mathcal{F}$ is superharmonic in $U_2$.
\end{lemma}

The following comparison principle will play an important role in the subsequent analysis.

\begin{proposition}[Comparison principle]\label{prop_comparison}
Assume \ref{eq_VD}, \ref{eq_BSP}. Let $\lambda\ge0$. Let $U$ be a bounded open set and $u,v\in\mathcal{F}$ satisfying that
$$-\Delta_pu+\lambda|u|^{p-2}u\ge-\Delta_pv+\lambda|v|^{p-2}v\text{ in }U,$$
that is,
$$\mathcal{E}(u;\varphi)+\lambda\int_X|u|^{p-2}u\varphi\mathrm{d} m\ge\mathcal{E}(v;\varphi)+\lambda\int_X|v|^{p-2}v\varphi\mathrm{d} m$$
for any non-negative $\varphi\in\mathcal{F}(U)$, and $(v-u)_+\in \mathcal{F}(U)$, or
$$\widetilde{u}\ge \widetilde{v}\text{ q.e. on }W\backslash U\text{ for some open set }W\supseteq \overline{U}\supseteq U.$$
Then $u\ge v$ in $U$.
\end{proposition}

\begin{proof}
By the strongly local property, we may assume that $(v-u)_+\in\mathcal{F}(U)$, then
$$\mathcal{E}(u;(v-u)_+)+\lambda\int_X|u|^{p-2}u(v-u)_+\mathrm{d} m\ge\mathcal{E}(v;(v-u)_+)+\lambda\int_X|v|^{p-2}v(v-u)_+\mathrm{d} m,$$
hence
\begin{align*}
&\mathcal{E}(u;(v-u)_+)-\mathcal{E}(v;(v-u)_+)\ge\lambda\int_X\left(|v|^{p-2}v-|u|^{p-2}u\right)(v-u)_+\mathrm{d} m\\
&=\lambda\int_{\{v>u\}}\left(|v|^{p-2}v-|u|^{p-2}u\right)(v-u)_+\mathrm{d} m=\lambda\int_{\{v>u\}}\left(|v|^{p-2}v-|u|^{p-2}u\right)(v-u)\mathrm{d} m\ge0,
\end{align*}
where in the last inequality, we use the facts that $\lambda\ge0$ and that
$$\left(|x|^{p-2}x-|y|^{p-2}y\right)\cdot(x-y)\ge0\text{ for any }x,y\in\mathbb{R}^d,d\ge1.$$
By the strongly local property, we have
\begin{align*}
&\Gamma(u;v-u)(\{v>u\})=\Gamma(u;(v-u)_+)(\{v>u\})=\mathcal{E}(u;(v-u)_+)\\
&\ge\mathcal{E}(v;(v-u)_+)=\Gamma(v;(v-u)_+)(\{v>u\})=\Gamma(v;v-u)(\{v>u\}).
\end{align*}
Let $\psi(t)=\Gamma(u+t(v-u);v-u)(\{v>u\})$ for any $t\in\mathbb{R}$, then by (\ref{eq_quasi_strict}), we have $\Gamma(u;v-u)(\{v>u\})=\psi(0)\le\psi(1)=\Gamma(v;v-u)(\{v>u\})$, hence $\Gamma(u;v-u)(\{v>u\})=\Gamma(v;v-u)(\{v>u\})$ and
$$0=\Gamma(v-u)(\{v>u\})=\Gamma((v-u)_+)(\{v>u\})=\mathcal{E}((v-u)_+).$$
Since $(v-u)_+\in\mathcal{F}(U)$, by \cite[Lemma 4.1]{Yan25c}, we have $(v-u)_+=0$, that is, $u\ge v$ in $U$.
\end{proof}

\begin{proposition}\label{prop_Pois_modi}
Assume \ref{eq_VD}, \ref{eq_BSP}. Let $U$ be a bounded open set, $u\in\mathcal{F}$ superharmonic in $U$, and $V$ an open subset of $U$. Then $H^{V}u\in\mathcal{F}$ is superharmonic in $U$, harmonic in $V$, and $H^{V}u\le u$ in $U$.
\end{proposition}

\begin{remark}
In the classical potential theory, the function $H^{V}u$ is called the Poisson modification of the superharmonic function $u$ in $V$, and it plays an important role in the theory of Perron solutions; see \cite[7.13]{HKM06} and \cite[Theorem 9.44, Section 10.9]{BB11}.
\end{remark}

\begin{proof}
By Proposition \ref{prop_harm_exist}, we have $H^{V}u\in\mathcal{F}$ is harmonic in $V$ and $\widetilde{H^{V}u}=\widetilde{u}$ q.e. on $X\backslash V$. Since $u$ is superharmonic in $V$, by Proposition \ref{prop_comparison}, we have $H^{V}u\le u$ in $V$, which gives $H^{V}u\le u$ in $U$. Applying Lemma \ref{lem_pasting} with $U_1=V$, $U_2=U$, $u_1=H^{V}u$, $u_2=u$, we conclude that $H^{V}u=u_1\wedge u_2$ is superharmonic in $U=U_2$.
\end{proof}

Let $U$ be a bounded open subset of $X$. We say that $U$ is regular if for any $u\in \mathcal{F}(U)$ which is bounded and harmonic in $U\backslash K$ for some compact set $K\subseteq U$, we have $\lim_{U\ni y\to x}u(y)=0$ uniformly for all $x\in \partial U$.

We say that the regular ball approximation condition \hypertarget{eq_RBA}{\text{RBA}} holds if for any ball $B$, for any $\delta\in(0,1)$, there exists a regular bounded open set $U$ such that $\delta B\subseteq U\subseteq B$.

We now present the relationship between the measure $\mu[\cdot]$ and the capacity $\mathrm{cap}(\cdot,\cdot)$, in particular, the variational characterizations given in (\ref{eq_cap_supchar}) and (\ref{eq_cap_infchar}) will be crucial for subsequent arguments.

\begin{lemma}\label{lem_cap_dual}
Assume \ref{eq_VD}, \ref{eq_BSP}. Let $U$ be a bounded open set.
\begin{enumerate}[label=(\arabic*)]
\item For any compact set $K\subseteq U$, there exists a unique $e_K\in\mathcal{F}(U)$ with $0\le e_K\le1$ in $X$, $\widetilde{e}_K=1$ q.e. on ${K}$, $e_K$ is superharmonic in $U$ and harmonic in $U\backslash {K}$, such that $\mathcal{E}(e_K)=\mathrm{cap}(K,X\backslash U)$. In particular, $\mathrm{supp}(\mu[e_K])\subseteq {K}$.
\item For any subset $K\subseteq U$, $\mathrm{cap}(K,X\backslash U)=0$ if and only if $\mathrm{cap}_1(K)=0$. Hence $\mathrm{cap}(\cdot,X\backslash U)$-q.e. in $U$ coincides with ($\mathrm{cap}_1$-)q.e. in $U$. In particular, for any $u\in\mathcal{F}$, its $\mathrm{cap}(\cdot,X\backslash U)$-quasi-continuous modification in $U$ coincides with its ($\mathrm{cap}_1$-)quasi-continuous modification in $U$.
\item For any $u\in\mathcal{F}$ superharmonic in $U$, for any compact set $K\subseteq U$, we have
\begin{equation}\label{eq_mu_cap}
\mu[u]({K})\le\mathcal{E}(u)^{(p-1)/p}\mathrm{cap}(K,X\backslash U)^{1/p}.
\end{equation}
In particular, $\mu[u]$ charges no set of zero capacity. Hence $\{\widetilde{v}:v\in\mathcal{F}(U)\}\subseteq L^1(U;\mu[u])$ and
$$\mathcal{E}(u;v)=\int_U\widetilde{v}\mathrm{d}\mu[u]\text{ for any }v\in\mathcal{F}(U).$$
\item For any compact set $K\subseteq U$, we have
\begin{align}
&\mathrm{cap}(K,X\backslash U)\nonumber\\
&=\sup\left\{\mu[u](K)\middle|
\begin{array}{l}
u\in\mathcal{F}(U)\text{ is non-negative superharmonic in }U,\\
\mathrm{supp}(\mu[u])\subseteq K,u\le 1\text{ in }U
\end{array}
\right\}.\label{eq_cap_supchar}
\end{align}
Moreover, assume \hyperlink{eq_LSC}{\text{LSC}} and $U$ is regular, then
\begin{align}
&\mathrm{cap}(K,X\backslash U)\nonumber\\
&=\inf\left\{\mu[v](U)\middle|
\begin{array}{l}
v\in\mathcal{F}(U)\text{ is non-negative superharmonic in }U,\\
\lim_{U\ni y\to x}v(y)=0\text{ uniformly for all }x\in\partial U,\\
\widetilde{v}\ge1\text{ q.e. on }K
\end{array}
\right\}.\label{eq_cap_infchar}
\end{align}
\end{enumerate}
\end{lemma}

\begin{remark}
In (1), the function $e_K$ and the measure $\mu[e_K]$ are also referred to as the capacitory potential and the capacitory measure (or equilibrium measure) of $(K,X\backslash U)$, see \cite[Section 6]{GH14}. In (2), this result corresponds to the classical fact that the sets of zero Sobolev capacity coincide with the sets of zero variational capacity, see \cite[Lemma 6.15]{BB11}. In (4), (\ref{eq_cap_infchar}) originates from \cite[(4.19)]{TW02AJM}. Similar to the classical case in \cite{TW02AJM}, boundary regularity needs to be taken into account in our setting.
\end{remark}

\begin{proof}
(1) Take $w\in\mathcal{F}\cap C_c(U)$ with $0\le w\le 1$ in $U$ and $w=1$ on $K$. Since $U\backslash K$ is a bounded open set, by Proposition \ref{prop_harm_exist}, there exists a unique $e_K\in\mathcal{F}$ such that $e_K$ is harmonic in $U\backslash K$ and $\widetilde{e}_K={w}$ q.e. on $X\backslash(U\backslash K)$, that is, $\widetilde{e}_K=1$ q.e. on $K$ and $\widetilde{e}_K=0$ q.e. on $X\backslash U$, moreover, $0\le e_K\le 1$ in $X$. By \cite[Lemma 3.2]{Yan25c}, we have
\begin{align*}
\mathcal{E}(e_K)&=\inf\{\mathcal{E}(u):u\in\mathcal{F},\widetilde{u}=1\text{ q.e. on }K,\widetilde{u}=0\text{ q.e. on }X\backslash U\}\\
&=\inf\{\mathcal{E}(u):u\in\mathcal{F}(U),\widetilde{u}=1\text{ q.e. on }K\}.
\end{align*}
By the definition of $\mathrm{cap}(\cdot,\cdot)$, we have $\mathcal{E}(e_K)=\mathrm{cap}(K,X\backslash U)$. For any non-negative $v\in\mathcal{F}\cap C_c(U)$, for any $t>0$, we have $(\widetilde{e}_K+tv)\wedge 1=1$ q.e. on $K$ and $(\widetilde{e}_K+tv)\wedge 1=0$ q.e. on $X\backslash U$, hence $\mathcal{E}(e_K)\le\mathcal{E}((e_K+tv)\wedge1)\le\mathcal{E}(e_K+tv)$, which gives $\mathcal{E}(e_K;v)=\frac{1}{p}\lim_{t\downarrow0}\left(\mathcal{E}(e_K+tv)-\mathcal{E}(e_K)\right)\ge0$, that is, $e_K\in\mathcal{F}$ is superharmonic in $U$. Moreover, for any $v\in\mathcal{F}\cap C_c(U)$ with $\mathrm{supp}(v)\subseteq U\backslash K$, we have $\int_Uv\mathrm{d}\mu[e_K]=\mathcal{E}(e_K;v)=0$, hence $\mathrm{supp}(\mu[e_K])\subseteq K$.

(2) By the same argument as in the proof of \cite[Lemma 8.3]{Yan25a}, we only need to prove this result for any compact set $K\subseteq U$. ``$\Rightarrow$": Assume $\mathrm{cap}(K,X\backslash U)=0$, let $e_K\in\mathcal{F}$ be given by (1), then $\mathcal{E}(e_K)=0$, by \cite[Lemma 4.1]{Yan25c}, we have $\lVert {e_K}\rVert_{L^p(X;m)}=0$, hence $\mathcal{E}_1(e_K)=0$, by \cite[Lemma 8.8]{Yan25a}, we have $\mathrm{cap}_1(K)\le\mathcal{E}_1(e_K)=0$. ``$\Leftarrow$": Assume $\mathrm{cap}_1(K)=0$, by \cite[Lemma 8.8]{Yan25a}, there exists $e\in\mathcal{F}$ with $0\le e\le1$ in $X$ and $\widetilde{e}=1$ q.e. on $K$ such that $\mathcal{E}_1(e)=\mathrm{cap}_1(K)=0$, then $e=0$ in $X$, in particular, $e\in\mathcal{F}(U)$, which gives $\mathrm{cap}(K,X\backslash U)\le\mathcal{E}(e)=0$.

(3) For any $\varphi\in\mathcal{F}\cap C_c(U)$ with $0\le\varphi\le1$ in $U$ and $\varphi=1$ on ${K}$, we have
$$\mu[u]({K})\le\int_U\varphi\mathrm{d}\mu[u]=\mathcal{E}(u;\varphi)\le\mathcal{E}(u)^{(p-1)/p}\mathcal{E}(\varphi)^{1/p}.$$
Taking the infimum with respect to $\varphi$, we have (\ref{eq_mu_cap}). Hence $\mu[u]$ charges no set of zero $\mathrm{cap}(\cdot,X\backslash U)$-capacity, by (2), $\mu[u]$ charges no set of zero ($\mathrm{cap}_1$-)capacity. For any $v\in\mathcal{F}(U)$, by \cite[Proposition 8.5]{Yan25a}, there exists $\{v_n\}\subseteq\mathcal{F}\cap C_c(U)$ such that $\{v_n\}$ is $\mathcal{E}_1$-convergent to $v$ and $\{v_n\}$ converges to $\widetilde{v}$ q.e. on $X$, then $\{v_n\}$ converges to $\widetilde{v}$ $\mu[u]$-a.e. in $U$. By Fatou's lemma, for any $n$, we have
\begin{align*}
&\int_U|v_n-\widetilde{v}|\mathrm{d}\mu[u]\le\varliminf_{m\to+\infty}\int_U|v_n-v_m|\mathrm{d}\mu[u]\\
&=\varliminf_{m\to+\infty}\mathcal{E}(u;|v_n-v_m|)\le\varliminf_{m\to+\infty}\mathcal{E}(u)^{(p-1)/p}\mathcal{E}(v_n-v_m)^{1/p},
\end{align*}
which implies $\widetilde{v}\in L^1(U;\mu[u])$ and $\{v_n\}$ is $L^1(U;\mu[u])$-convergent to $\widetilde{v}$, hence
$$\int_U\widetilde{v}\mathrm{d}\mu[u]=\lim_{n\to+\infty}\int_Uv_n\mathrm{d}\mu[u]=\lim_{n\to+\infty}\mathcal{E}(u;v_n)=\mathcal{E}(u;v).$$

(4) Let $e_K\in\mathcal{F}(U)$ be given by (1), then
\begin{align*}
&\mu[e_K](U)=\mu[e_K](K)=\int_K\widetilde{e}_K\mathrm{d}\mu[e_K]=\int_U\widetilde{e}_K\mathrm{d}\mu[e_K]=\mathcal{E}(e_K)=\mathrm{cap}(K,X\backslash U).
\end{align*}

Firstly, we prove (\ref{eq_cap_supchar}). ``$\le$": It is obvious since $e_K$ is one such function $u$ and $\mathrm{cap}(K,X\backslash U)=\mu[e_K](K)$. ``$\ge$": For any such function $u$, we have
\begin{align*}
&\mu[u](K)=\int_K\mathrm{d}\mu[u]=\int_K\widetilde{e}_K\mathrm{d}\mu[u]=\int_U\widetilde{e}_K\mathrm{d}\mu[u]=\mathcal{E}(u;e_K)\le\mathcal{E}(u)^{(p-1)/p}\mathcal{E}(e_K)^{1/p},
\end{align*}
where
$$\mathcal{E}(u)=\int_U\widetilde{u}\mathrm{d}\mu[u]=\int_K\widetilde{u}\mathrm{d}\mu[u]\le\int_K\mathrm{d}\mu[u]=\mu[u](K),$$
hence $\mu[u](K)\le\mathcal{E}(e_K)=\mathrm{cap}(K,X\backslash U)$.

Secondly, we prove (\ref{eq_cap_infchar}). ``$\ge$": Since $U$ is regular and $e_K\in \mathcal{F}(U)$ is bounded and harmonic in $U\backslash K$, we have $e_K$ is one such function $v$, then the result follows directly by the fact $\mathrm{cap}(K,X\backslash U)=\mu[e_K](U)$. ``$\le$": For any such function $v$, we only need to show that $\mu[e_K](K)\le \mu[v](U)$. Since $e_K, v\in \mathcal{F}(U)$ are non-negative superharmonic in $U$, by \hyperlink{eq_LSC}{\text{LSC}}, we may assume that $e_K, v$ are lower semi-continuous in $U$. Let $\eta\in \mathcal{F}\cap C_c(X)$ satisfy $\eta=1$ on $\overline{U}$. For any $\varepsilon\in(0,\frac{1}{2})$, let $\phi_\varepsilon:\mathbb{R}\to \mathbb{R}$ be given by $t\mapsto (1-2\varepsilon)t+\varepsilon$, then $\phi_\varepsilon(e_K)\eta\in \mathcal{F}$ satisfies $\varepsilon\le \phi_\varepsilon(e_K)\eta\le 1-\varepsilon$ in $U$, is lower semi-continuous, and superharmonic in $U$. Let $w=(\phi_\varepsilon(e_K)\eta)\wedge v$, then by Lemma \ref{lem_pasting}, we have $w\in \mathcal{F}$ is superharmonic in $U$. By the lower semi-continuity, we have $\{v>1-\varepsilon\}\supseteq K$, $\{v>\varepsilon\}\subseteq U$ are open sets. Moreover, since $\lim_{U\ni y\to x}v(y)=0$ uniformly for all $x\in \partial U$, there exist open sets $W,V$ with $K\subseteq W\subseteq \overline{W}\subseteq V\subseteq \overline{V}\subseteq U$ such that
$$w=\phi_\varepsilon(e_K)\eta=(1-2\varepsilon)e_K+\varepsilon\text{ in }W,$$
and $w=v\text{ in }U\backslash \overline{V}$. By Proposition \ref{prop_super_meas}, we have
$$(1-2\varepsilon)^{p-1}\mu[e_K](K)\le (1-2\varepsilon)^{p-1}\mu[e_K](W)=\mu[w](W)\le \mu[w](U)=\mu[v](U).$$
Letting $\varepsilon\downarrow0$, we have $\mu[e_K](K)\le\mu[v](U)$.
\end{proof}

\begin{lemma}\label{lem_FU_modi}
Assume \ref{eq_VD}, \ref{eq_BSP}. Let $U$ be a bounded open set and $u\in\mathcal{F}$ superharmonic in $U$. Then there exists $v\in\mathcal{F}(U)$ which is non-negative superharmonic in $U$ such that $\mu[v]=\mu[u]$ in $U$. Moreover, if $(-u)_+\in \mathcal{F}(U)$ or $\widetilde{u}\ge0$ q.e. on $W\backslash U$ for some open set $W\supseteq\overline{U}\supseteq U$, then $v\le u$ in $U$.
\end{lemma}

\begin{proof}
By \cite[Lemma 4.1]{Yan25c}, we have $(\mathcal{F}(U),\mathcal{E})$ is a Banach space. Since $u\in\mathcal{F}$ is superharmonic in $U$, by Lemma \ref{lem_cap_dual} (3), we have $\varphi\mapsto\int_U\widetilde{\varphi}\mathrm{d}\mu[u]$ is a bounded linear functional on $(\mathcal{F}(U),\mathcal{E})$, by \cite[Theorem 3.24]{KS24a}, there exists $v\in\mathcal{F}(U)$ such that $\mathcal{E}(v;\varphi)=\int_U\widetilde{\varphi}\mathrm{d}\mu[u]$ for any $\varphi\in \mathcal{F}(U)$, which gives $v$ is superharmonic in $U$ and $\mu[v]=\mu[u]$ in $U$. Since $\widetilde{v}=0$ q.e. on $X\backslash U$, by Proposition \ref{prop_comparison}, we have $v$ is non-negative in $U$. Moreover, if $(-u)_+\in \mathcal{F}(U)$, which implies $(v-u)_+\in \mathcal{F}(U)$, or if $\widetilde{u}\ge0=\widetilde{v}$ q.e. on $W\backslash U$, then by Proposition \ref{prop_comparison} again, we have $v\le u$ in $U$.
\end{proof}

We establish the validity of the \hyperlink{eq_RBA}{\text{RBA}} condition as follows; the proof is postponed to Section \ref{sec_RBA} to keep this section concise.

\begin{proposition}\label{prop_RBA}
Assume \ref{eq_VD}, \hyperlink{eq_LSC}{\text{LSC}}, \ref{eq_BSP}, \ref{eq_EHI}, \ref{eq_EHIann}, \hyperlink{eq_cap}{$\text{cap}(\Psi)$}. Then \hyperlink{eq_RBA}{\text{RBA}} holds.
\end{proposition}

We define the Wolff potential as follows. Let $\mu$ be a positive Radon measure on $X$. For any $x\in X$, for any $R>0$, let
$$\mathcal{W}^\mu(x,R)=\sum_{n=0}^{+\infty} \left(\frac{\mu(B(x,\frac{1}{2^n}R))}{\mathrm{cap}(B(x,\frac{1}{2^{n+1}}R),X\backslash B(x,\frac{1}{2^n}R))}\right)^{{1}/{(p-1)}}.$$
Due to the measurability issue of the function $r\mapsto \mathrm{cap}(B(x,\frac{1}{2}r),X\backslash B(x,r))$, we define the potential in the form of a summation rather than an integral, as in the classical case in \cite{KM92Pisa,KM94Acta}. The key result of this section is as follows.

\begin{theorem}\label{thm_Wolff}
Assume \ref{eq_VD}, \hyperlink{eq_LSC}{\text{LSC}}, \ref{eq_BSP}, \ref{eq_EHI}, \ref{eq_EHIann}, \hyperlink{eq_cap}{$\text{cap}(\Psi)$}. Then there exists $C>0$ such that for any bounded open set $U$, for any $u\in\mathcal{F}$ which is non-negative superharmonic in $U$, for any Lebesgue point $x_0\in U$ of $u$, for any $R>0$ with $B(x_0,4R)\subseteq U$, we have
\begin{equation*}\label{eq_Wolff}\tag*{$\text{Wol}$}
\frac{1}{C}\mathcal{W}^{\mu[u]}(x_0,R)\le u(x_0)\le C \left(\einf_{B(x_0,R)}u+\mathcal{W}^{\mu[u]}(x_0,2R)\right).
\end{equation*}
\end{theorem}

\begin{remark}
In classical quasi-linear elliptic potential theory, the above pointwise estimate, known as the Wolff potential estimate, was established by Kilpel\"ainen and Mal\'y \cite{KM92Pisa,KM94Acta}. It was later extended to the sub-elliptic setting by Trudinger and Wang \cite{TW02AJM}, whose argument relies on the elliptic Harnack inequality. We follow the approach of \cite{TW02AJM}.
\end{remark}

For simplicity, by appropriately adjusting the constants in the following results, we always assume that the annuli appearing therein satisfy the assumption as in the \ref{eq_EHIann} condition. The proof relies on the following result, which corresponds to \cite[LEMMA 5.1]{TW02AJM}.

\begin{proposition}\label{prop_Wolff}
Assume \ref{eq_VD}, \hyperlink{eq_LSC}{\text{LSC}}, \ref{eq_BSP}, \ref{eq_EHI}, \ref{eq_EHIann}, \hyperlink{eq_cap}{$\text{cap}(\Psi)$}. Then there exists $C>0$ such that for any ball $B=B(x_0,R)$, the following holds:
\begin{enumerate}[label=(\alph*),ref=(\alph*)]
\item\label{item_Wolff1} If $u\in\mathcal{F}$ is non-negative superharmonic in $B$, then
\begin{equation}\label{eq_Wolff1}
\einf_{\frac{1}{2}B}u\ge \frac{1}{C}\left(\frac{\mu[u]\left(\frac{9}{10}B\right)}{\mathrm{cap}\left(\frac{1}{2}B,X\backslash B\right)}\right)^{1/(p-1)}.
\end{equation}
\item\label{item_Wolff2} If $u\in \mathcal{F}$ is non-negative superharmonic in $2B$ and satisfies
$$\mu[u]=0\text{ in }\omega=\left(\frac{5}{8}B\backslash\overline{\frac{3}{8}B}\right)\cup \left(\frac{5}{4}B\backslash \overline{\frac{3}{4}B}\right),$$
then
\begin{equation}\label{eq_Wolff2}
\esup_{\frac{11}{20}B\backslash {\frac{9}{20}B}}u-\esup_{\frac{11}{10}B\backslash {\frac{9}{10}B}}u\le C\left(\frac{\mu[u](B)}{\mathrm{cap}\left(\frac{1}{2}B,X\backslash B\right)}\right)^{1/(p-1)}.
\end{equation}
\end{enumerate}
\end{proposition}

\begin{proof}
\ref{item_Wolff1} By Lemma \ref{lem_FU_modi}, there exists $v\in\mathcal{F}(\frac{99}{100}B)$ which is non-negative superharmonic in $\frac{99}{100}B$ such that $\mu[v]=\mu[u]$ in $\frac{99}{100}B$ and $v\le u$ in $\frac{99}{100}B$, by replacing $u$ with $v$, we may assume that $u\in\mathcal{F}(\frac{99}{100}B)\subseteq \mathcal{F}(B)$, hence assume $u\in \mathcal{F}(B)$. Since $\widetilde{u}=\widetilde{H^{B\backslash\overline{\frac{9}{10}B}}u}$ q.e. on $\frac{9}{10}B$, by Proposition \ref{prop_super_meas}, we have $\mu[u](\frac{9}{10}B)=\mu[H^{B\backslash\overline{\frac{9}{10}B}}u](\frac{9}{10}B)$, by replacing $u$ with $H^{B\backslash\overline{\frac{9}{10}B}}u$, we may assume that $u$ is harmonic in $B\backslash\overline{\frac{9}{10}B}$, by \ref{eq_EHI}, we have $u\in C(B\backslash \overline{\frac{9}{10}B})$. By \ref{eq_EHIann}, there exists $C_1>0$ such that
$$\esup_{\frac{39}{40}B\backslash {\frac{19}{20}B}}u\le C_1\einf_{\frac{39}{40}B\backslash{\frac{19}{20}B}}u<+\infty.$$
Let $u_1=H^{\frac{19}{20}B}u\in\mathcal{F}(B)$. Since $\widetilde{u}_1=\widetilde{u}$ q.e. on $B\backslash {\frac{19}{20}B}$, by Proposition \ref{prop_super_meas}, we have $\mu[u_1](B)=\mu[u](B)$. Since $u_1$ is harmonic in $\frac{19}{20}B$, by Proposition \ref{prop_comparison}, we have
$$\esup_{\frac{19}{20}B}u_1\le\esup_{\frac{39}{40}B\backslash \frac{19}{20}B}u_1=\esup_{\frac{39}{40}B\backslash \frac{19}{20}B}u\le C_1\einf_{\frac{39}{40}B\backslash \frac{19}{20}B}u=C_1\einf_{\frac{39}{40}B\backslash \frac{19}{20}B}u_1\le C_1\einf_{\frac{19}{20}B}u_1.$$
By Proposition \ref{prop_Pois_modi}, we have $u_1\le u$ in $\frac{19}{20}B$, in particular, we have
$$\einf_{\frac{1}{2}B}u\ge\einf_{\frac{1}{2}B}u_1\ge\einf_{{\frac{19}{20}B}}u_1.$$
Let $t=\esup_{\frac{19}{20}B}u_1<+\infty$, then by \hyperlink{eq_LSC}{\text{LSC}}, we have $\esup_{\overline{\frac{19}{20}B}}u_1=t$. Since $u\in\mathcal{F}(B)$ is harmonic in $B\backslash\overline{\frac{9}{10}B}$ and $\widetilde{u}=\widetilde{u}_1$ q.e. on $B\backslash \frac{19}{20}B$, we have $u_1\in \mathcal{F}(B)$ is harmonic in $B\backslash \overline{\frac{19}{20}B}$, that is, $\mathrm{supp}(\mu[u_1])\subseteq \overline{\frac{19}{20}B}$. By Proposition \ref{prop_comparison}, we have $u_1\le\esup_{\overline{\frac{19}{20}B}}u_1=t$ in $B\backslash \overline{\frac{19}{20}B}$, hence $u_1\le t$ in $B$. Without loss of generality, we may assume that $t>0$; otherwise $u_1=0$ in $B$, which gives $0=\mu[u_1](B)=\mu[u](B)\ge \mu[u](\frac{9}{10}B)$, then (\ref{eq_Wolff1}) is obvious. By (\ref{eq_cap_supchar}),
$$\frac{1}{t^{p-1}}\mu[u_1]\left(\overline{\frac{19}{20}B}\right)=\mu\left[\frac{u_1}{t}\right]\left(\overline{\frac{19}{20}B}\right)\le \mathrm{cap}\left(\overline{\frac{19}{20}B},X\backslash B\right).$$
By \ref{eq_VD}, \hyperlink{eq_cap}{$\text{cap}(\Psi)$}, there exists $C_2>0$ depending only on $C_\Psi, C_{VD}, C_{cap}, A_{cap}$ such that
$$\mathrm{cap}\left(\overline{\frac{19}{20}B},X\backslash B\right)\le \mathrm{cap}\left({\frac{39}{40}B},X\backslash B\right)\le C_2 \mathrm{cap}\left(\frac{1}{2}B,X\backslash B\right).$$
In summary, we have
$$\mu[u]\left(\frac{9}{10}B\right)\le\mu[u](B)=\mu[u_1](B)=\mu[u_1]\left(\overline{\frac{19}{20}B}\right)\le C_2t^{p-1}\mathrm{cap}\left(\frac{1}{2}B,X\backslash B\right),$$
hence
$$\einf_{\frac{1}{2}B}u\ge\einf_{{\frac{19}{20}B}}u_1\ge \frac{1}{C_1}t\ge \frac{1}{C_1C_2^{1/(p-1)}}\left(\frac{\mu[u]\left(\frac{9}{10}B\right)}{\mathrm{cap}\left(\frac{1}{2}B,X\backslash B\right)}\right)^{1/(p-1)}.$$

\ref{item_Wolff2} Let $\omega^*=\left(\frac{4}{9}B\right)\cup \left(\frac{4}{5}B\backslash \overline{\frac{5}{9}B}\right)$, then $B\backslash \overline{\omega}\subseteq \omega^*$. Let $u_2=H^{\omega^*}u$, then $\widetilde{u}_2=\widetilde{u}$ q.e. on $\left(\overline{\frac{5}{9}B}\backslash {\frac{4}{9}B}\right)\cup \left(2B\backslash {\frac{4}{5}B}\right)$, $\mathrm{supp}(\mu[u_2])\cap \left(\frac{5}{4}B\right)\subseteq \partial \omega^*$, and
$$\esup_{\frac{11}{20}B\backslash {\frac{9}{20}B}}u-\esup_{\frac{11}{10}B\backslash {\frac{9}{10}B}}u=\esup_{\frac{11}{20}B\backslash {\frac{9}{20}B}}u_2-\esup_{\frac{11}{10}B\backslash {\frac{9}{10}B}}u_2.$$
Since $u=u_2$ in $B\backslash \overline{\frac{4}{5}B}$, by Proposition \ref{prop_super_meas}, we have $\mu[u](B)=\mu[u_2](B)$. By \hyperlink{eq_RBA}{\text{RBA}}, there exists a regular bounded open set $U$ such that $\frac{9}{10}B\subseteq U\subseteq B$. By Lemma \ref{lem_FU_modi}, there exists $u_3\in \mathcal{F}(U)$ which is non-negative superharmonic in $U$ such that $\mu[u_3]=\mu[u_2]$ in $U$, in particular, $\mu[u_3]=0$ in $\left(\frac{5}{9}B\backslash \overline{\frac{4}{9}B}\right)\cup \left(U\backslash \overline{\frac{4}{5}B}\right)$, hence by \ref{eq_EHI}, we have $u_3$ is continuous in $\left(\frac{5}{9}B\backslash \overline{\frac{4}{9}B}\right)\cup \left(U\backslash \overline{\frac{4}{5}B}\right)$. Since $\widetilde{u}_3=0\ge \widetilde{u}_2-\esup_{\frac{11}{10}B\backslash U}u_2$ q.e. on $\frac{11}{10}B\backslash U$, by Proposition \ref{prop_comparison}, we have
$$u_3\ge u_2-\esup_{\frac{11}{10}B\backslash {U}}u_2\ge u_2-\esup_{\frac{11}{10}B\backslash {\frac{9}{10}B}}u_2\text{ in }U,$$
which gives
$$\esup_{\frac{11}{20}B\backslash \frac{9}{20}B}u_3\ge \esup_{\frac{11}{20}B\backslash \frac{9}{20}B}u_2-\esup_{\frac{11}{10}B\backslash \frac{9}{10}B}u_2.$$
Since $u_3$ is non-negative harmonic in $\frac{5}{9}B\backslash \overline{\frac{4}{9}B}$, by \ref{eq_EHIann}, there exists $C_3>0$ such that
$$\esup_{\frac{11}{20}B\backslash {\frac{9}{20}B}}u_3\le C_3\einf_{\frac{11}{20}B\backslash {\frac{9}{20}B}}u_3<+\infty.$$
Let $s=\einf_{\frac{11}{20}B\backslash \frac{9}{20}B}u_3$, without loss of generality, we may assume that $s>0$; otherwise, the LHS of (\ref{eq_Wolff2}) satisfies $\le C_3s=0$, and the statement is trivial. By Proposition \ref{prop_comparison}, we have $u_3\ge s$ in $\frac{9}{20}B$. Since $u_3\in C\left(\frac{5}{9}B\backslash \overline{\frac{4}{9}B}\right)$, we have $\widetilde{u}_3\ge s$ q.e. on $\overline{\frac{9}{20}B}$. Since $U$ is regular and $u_3\in \mathcal{F}(U)$ is bounded and harmonic in $U\backslash \overline{\frac{17}{20}B}$, by (\ref{eq_cap_infchar}), we have
\begin{align*}
&\mathrm{cap}\left(\overline{\frac{9}{20}B},X\backslash U\right)\le \mu \left[\frac{u_3}{s}\right](U)=\frac{1}{s^{p-1}}\mu[u_3](U)\\
&=\frac{1}{s^{p-1}}\mu[u_2](U)\le \frac{1}{s^{p-1}}\mu[u_2](B)=\frac{1}{s^{p-1}}\mu[u](B).
\end{align*}
By \ref{eq_VD}, \hyperlink{eq_cap}{$\text{cap}(\Psi)$}, there exists $C_4>0$ depending only on $C_\Psi, C_{VD}, C_{cap}, A_{cap}$ such that
$$\mathrm{cap}\left(\overline{\frac{9}{20}B},X\backslash U\right)\ge \mathrm{cap}\left({\frac{9}{20}B},X\backslash B\right)\ge \frac{1}{C_4}\mathrm{cap}\left(\frac{1}{2}B,X\backslash B\right).$$
In summary, we have
$$\esup_{\frac{11}{20}B\backslash {\frac{9}{20}B}}u-\esup_{\frac{11}{10}B\backslash {\frac{9}{10}B}}u\le\esup_{\frac{11}{20}B\backslash {\frac{9}{20}B}}u_3\le C_3s\le C_3C_4^{1/(p-1)}\left(\frac{\mu[u](B)}{\mathrm{cap}\left(\frac{1}{2}B,X\backslash B\right)}\right)^{1/(p-1)}.$$
\end{proof}

We give the proof of Theorem \ref{thm_Wolff} as follows.

\begin{proof}[Proof of Theorem \ref{thm_Wolff}]
By \hyperlink{eq_LSC}{\text{LSC}}, we may assume that $u$ is lower semi-continuous in $U$. Firstly, we prove the lower bound. For any $r\in(0,2R]$, let $B=B(x_0,r)$, then $u$ is non-negative superharmonic in $2B$. By Lemma \ref{lem_FU_modi}, there exists $u_1\in \mathcal{F}(B)$ which is non-negative superharmonic in $B$ such that $\mu[u_1]=\mu[u]$ in $B$, moreover, $u_1\le u$ in $B$. Since $\widetilde{u}_1=0\le \widetilde{u}-\einf_{\frac{3}{2}B}u$ q.e. on $\frac{3}{2}B\backslash B$, by Proposition \ref{prop_comparison}, we have
$$\einf_{\frac{3}{4}B}u_1\le \einf_{\frac{3}{4}B}u-\einf_{\frac{3}{2}B}u.$$
Let $u_2=H^{\frac{7}{8}B}u_1$, then by Proposition \ref{prop_Pois_modi}, we have $u_2$ is non-negative superharmonic in $B$, harmonic in $\frac{7}{8}B$, $\widetilde{u}_2=\widetilde{u}_1$ q.e. on $B\backslash \frac{7}{8}B$, and $u_2\le u_1$ in $B$. By (\ref{eq_Wolff1}), we have
$$\einf_{\frac{1}{2}B}u_2\ge \frac{1}{C_1}\left(\frac{\mu[u_2]\left(\frac{9}{10}B\right)}{\mathrm{cap}\left(\frac{1}{2}B,X\backslash B\right)}\right)^{1/(p-1)},$$
where $C_1$ is the positive constant appearing therein. Since $u_2=u_1$ in $\frac{9}{10}B\backslash \overline{\frac{7}{8}B}$, by Proposition \ref{prop_super_meas}, we have
$$\mu[u_2]\left(\frac{9}{10}B\right)=\mu[u_1]\left(\frac{9}{10}B\right)=\mu[u]\left(\frac{9}{10}B\right)\ge \mu[u]\left(\frac{1}{2}B\right).$$
Since $u_2$ is non-negative harmonic in $\frac{7}{8}B$, by \ref{eq_VD}, \ref{eq_EHI}, there exists $C_2>0$ depending only on $C_{VD}, C_H, A_H$ such that
$$\einf_{\frac{1}{2}B}u_2\le \esup_{\frac{3}{4}B}u_2\le C_2\einf_{\frac{3}{4}B}u_2\le C_2\einf_{\frac{3}{4}B}u_1.$$
Hence
$$\einf_{\frac{3}{4}B}u-\einf_{\frac{3}{2}B}u\ge \frac{1}{C_1C_2}\left(\frac{\mu[u]\left(\frac{1}{2}B\right)}{\mathrm{cap}\left(\frac{1}{2}B,X\backslash B\right)}\right)^{1/(p-1)}.$$
By \ref{eq_VD}, \hyperlink{eq_cap}{$\text{cap}(\Psi)$}, there exists $C_3>0$ depending only on $C_{\Psi}, C_{VD}, A_{cap}, C_{cap}$ such that
$${\mathrm{cap}\left(\frac{1}{2}B,X\backslash B\right)}\le C_3{\mathrm{cap}\left(\frac{1}{4}B,X\backslash \frac{1}{2}B\right)},$$
hence
$$\einf_{\frac{3}{4}B}u-\einf_{\frac{3}{2}B}u\ge\frac{1}{C_1C_2C_3^{1/(p-1)}}\left(\frac{\mu[u]\left(\frac{1}{2}B\right)}{\mathrm{cap}\left(\frac{1}{4}B,X\backslash \frac{1}{2}B\right)}\right)^{1/(p-1)}.$$
Then for any $k=0,1,\ldots$, we have
\begin{align*}
&\einf_{B(x_0,\frac{3}{2^{k+1}}R)}u-\einf_{B(x_0,\frac{3}{2^{k}}R)}u\\
&\ge\frac{1}{C_1C_2C_3^{1/(p-1)}}\left(\frac{\mu[u]\left(B(x_0,\frac{1}{2^{k}}R)\right)}{\mathrm{cap}\left(B(x_0,\frac{1}{2^{k+1}}R),X\backslash B(x_0,\frac{1}{2^{k}}R)\right)}\right)^{1/(p-1)}.
\end{align*}
Summing over $k=0,1,\ldots,l$, we have
\begin{align*}
&\einf_{B(x_0,\frac{3}{2^{l+1}}R)}u\ge \einf_{B(x_0,\frac{3}{2^{l+1}}R)}u-\einf_{B(x_0,{3}R)}u\\
&\ge\frac{1}{C_1C_2C_3^{1/(p-1)}}\sum_{k=0}^l\left(\frac{\mu[u]\left(B(x_0,\frac{1}{2^{k}}R)\right)}{\mathrm{cap}\left(B(x_0,\frac{1}{2^{k+1}}R),X\backslash B(x_0,\frac{1}{2^{k}}R)\right)}\right)^{1/(p-1)}.
\end{align*}
Since $x_0$ is a Lebesgue point of $u$, we have
\begin{align*}
&u(x_0)=\lim_{l\to+\infty}\frac{1}{m({B(x_0,\frac{3}{2^{l+1}}R)})}\int_{{B(x_0,\frac{3}{2^{l+1}}R)}}u\mathrm{d} m\ge\varlimsup_{l\to+\infty}\einf_{B(x_0,\frac{3}{2^{l+1}}R)}u\\
&\ge\frac{1}{C_1C_2C_3^{1/(p-1)}}\sum_{k=0}^{+\infty}\left(\frac{\mu[u]\left(B(x_0,\frac{1}{2^{k}}R)\right)}{\mathrm{cap}\left(B(x_0,\frac{1}{2^{k+1}}R),X\backslash B(x_0,\frac{1}{2^{k}}R)\right)}\right)^{1/(p-1)}.
\end{align*}

Secondly, we prove the upper bound. Let $\omega=\cup_{k=0}^{+\infty}(B(x_0,\frac{5}{4}\frac{1}{2^k}R)\backslash\overline{B(x_0,\frac{3}{4}\frac{1}{2^k}R)})$ and $v=H^\omega u$, then by Proposition \ref{prop_Pois_modi}, we have $v$ is non-negative superharmonic in $B(x_0,4R)$, harmonic in $\omega$, $\widetilde{v}=\widetilde{u}$ q.e. on $B(x_0,4R)\backslash\omega$, and $v\le u$ in $B(x_0,4R)$. By Proposition \ref{prop_super_meas}, for any $r\in(0,{R}]$, we have $\mu[v](B(x_0,r))\le\mu[u](B(x_0,2r))$. By (\ref{eq_Wolff2}), for any $k=0,1,\ldots$,
\begin{align*}
&\esup_{B(x_0,\frac{11}{10}\frac{1}{2^{k+1}}R)\backslash {B(x_0,\frac{9}{10}\frac{1}{2^{k+1}}R)}}v\\
&\le\esup_{B(x_0,\frac{11}{10}\frac{1}{2^{k}}R)\backslash {B(x_0,\frac{9}{10}\frac{1}{2^{k}}R)}}v+C_1 \left(\frac{\mu[v](B(x_0,\frac{1}{2^k}R))}{\mathrm{cap}\left(B(x_0,\frac{1}{2^{k+1}}R),X\backslash B(x_0,\frac{1}{2^k}R)\right)}\right)^{1/(p-1)}\\
&\le\esup_{B(x_0,\frac{11}{10}\frac{1}{2^{k}}R)\backslash {B(x_0,\frac{9}{10}\frac{1}{2^{k}}R)}}v+C_1 \left(\frac{\mu[u](B(x_0,\frac{1}{2^{k-1}}R))}{\mathrm{cap}\left(B(x_0,\frac{1}{2^{k+1}}R),X\backslash B(x_0,\frac{1}{2^k}R)\right)}\right)^{1/(p-1)}.
\end{align*}
Summing over $k=0,1,\ldots,l$, we have
\begin{align*}
&\esup_{B(x_0,\frac{11}{10}\frac{1}{2^{l+1}}R)\backslash {B(x_0,\frac{9}{10}\frac{1}{2^{l+1}}R)}}v\\
&\le\esup_{B(x_0,\frac{11}{10}R)\backslash {B(x_0,\frac{9}{10}R)}}v+C_1\sum_{k=0}^l\left(\frac{\mu[u](B(x_0,\frac{1}{2^{k-1}}R))}{\mathrm{cap}\left(B(x_0,\frac{1}{2^{k+1}}R),X\backslash B(x_0,\frac{1}{2^k}R)\right)}\right)^{1/(p-1)}.
\end{align*}
Since $v$ is non-negative harmonic in $B(x_0,\frac{5}{4}R)\backslash \overline{B(x_0,\frac{3}{4}R)}$, by \ref{eq_EHIann}, there exists $C_4>0$ such that
\begin{equation*}
\esup_{B(x_0,\frac{11}{10}R)\backslash {B(x_0,\frac{9}{10}R)}}v\le C_4\einf_{B(x_0,\frac{11}{10}R)\backslash {B(x_0,\frac{9}{10}R)}}v\le C_4\einf_{B(x_0,\frac{11}{10}R)\backslash {B(x_0,R)}}v,
\end{equation*}
where
$$\einf_{B(x_0,\frac{11}{10}R)\backslash {B(x_0,R)}}v\le\einf_{B(x_0,R)}v\le\einf_{B(x_0,R)}u.$$
By the lower semi-continuity, we have $u(x_0)\le\lim_{l\to+\infty}\einf_{B(x_0,\frac{9}{10}\frac{1}{2^{l}}R)}u$, where
\begin{align*}
&\einf_{B(x_0,\frac{9}{10}\frac{1}{2^{l}}R)}u= \einf_{B(x_0,\frac{9}{10}\frac{1}{2^{l}}R)\backslash B(x_0,\frac{11}{10}\frac{1}{2^{l+1}}R)}u=\einf_{B(x_0,\frac{9}{10}\frac{1}{2^{l}}R)\backslash B(x_0,\frac{11}{10}\frac{1}{2^{l+1}}R)}v\\
&\le \einf_{B(x_0,\frac{11}{10}\frac{1}{2^{l+1}}R)\backslash B(x_0,\frac{9}{10}\frac{1}{2^{l+1}}R)}v\le \esup_{B(x_0,\frac{11}{10}\frac{1}{2^{l+1}}R)\backslash B(x_0,\frac{9}{10}\frac{1}{2^{l+1}}R)}v,
\end{align*}
hence
$$u(x_0)\le C_4\einf_{B(x_0,R)}u+C_1\sum_{k=0}^{+\infty}\left(\frac{\mu[u](B(x_0,\frac{1}{2^{k-1}}R))}{\mathrm{cap}\left(B(x_0,\frac{1}{2^{k+1}}R),X\backslash B(x_0,\frac{1}{2^k}R)\right)}\right)^{1/(p-1)}.$$
By \ref{eq_VD}, \hyperlink{eq_cap}{$\text{cap}(\Psi)$}, there exists $C_5>0$ depending only on $C_\Psi,C_{VD},C_{cap},A_{cap}$ such that for any $k=0,1,\ldots$, we have
$$\mathrm{cap}\left(B(x_0,\frac{1}{2^{k+1}}R),X\backslash B(x_0,\frac{1}{2^k}R)\right)\ge \frac{1}{C_5}\mathrm{cap}\left(B(x_0,\frac{1}{2^{k}}R),X\backslash B(x_0,\frac{1}{2^{k-1}}R)\right).$$
Hence
$$u(x_0)\le C_4\einf_{B(x_0,R)}u+C_1C_5^{1/(p-1)}\sum_{k=0}^{+\infty}\left(\frac{\mu[u](B(x_0,\frac{1}{2^{k}}2R))}{\mathrm{cap}\left(B(x_0,\frac{1}{2^{k+1}}2R),X\backslash B(x_0,\frac{1}{2^k}2R)\right)}\right)^{1/(p-1)}.$$
\end{proof}

\begin{proposition}\label{prop_EF}
Assume \ref{eq_VD}, \hyperlink{eq_LSC}{\text{LSC}}, \ref{eq_BSP}, \ref{eq_EHI}, \ref{eq_EHIann}, \hyperlink{eq_cap}{$\text{cap}(\Psi)$}. Then there exists $C>0$ such that for any ball $B=B(x_0,R)$, for any $u\in\mathcal{F}(B)$ satisfying that $-\Delta_pu=1$ in $B$, that is, $u\in\mathcal{F}(B)$ is superharmonic in $B$ with $\mu[u]=m$ in $B$, we have
\begin{align}
\esup_{B}u&\le C\Psi(R)^{\frac{1}{p-1}},\label{eq_EF1}\\
\einf_{\frac{1}{2}B}u&\ge \frac{1}{C}\Psi(R)^{\frac{1}{p-1}}.\label{eq_EF2}
\end{align}
\end{proposition}

\begin{remark}
This result plays a role analogous to the mean exit time estimate in the case $p=2$, as formulated in condition $(\text{E}_{\text{F}})$ of \cite[Definition 3.10]{GH14}.
\end{remark}

\begin{proof}
By Proposition \ref{prop_comparison}, we have $u$ is non-negative in $B$. By \ref{eq_VD}, \hyperlink{eq_cap}{$\text{cap}(\Psi)$}, for any $x\in X$, for any $r>0$, we have
$$\mathcal{W}^{m}(x,r)\asymp \sum_{n=0}^{+\infty}\left(\frac{V(x,\frac{1}{2^n}r)}{\frac{V(x,\frac{1}{2^n}r)}{\Psi(\frac{1}{2^n}r)}}\right)^{\frac{1}{p-1}}=\sum_{n=0}^{+\infty}\Psi\left(\frac{1}{2^n}r\right)^{\frac{1}{p-1}}\asymp\Psi(r)^{\frac{1}{p-1}}.$$

Firstly, we prove (\ref{eq_EF2}). For any Lebesgue point $x\in \frac{1}{2}B$ of $u$, by the lower bound in \ref{eq_Wolff}, we have
$$u(x)\ge \frac{1}{C_1}\mathcal{W}^{\mu[u]}\left(x,\frac{1}{8}R\right)=\frac{1}{C_1}\mathcal{W}^{m}\left(x,\frac{1}{8}R\right)\asymp\Psi\left(\frac{1}{8}R\right)^{\frac{1}{p-1}}\asymp \Psi(R)^{\frac{1}{p-1}},$$
where $C_1$ is the positive constant appearing therein, which gives $\einf_{\frac{1}{2}B}u\gtrsim \Psi(R)^{\frac{1}{p-1}}$.

Secondly, we prove (\ref{eq_EF1}). By \hyperlink{eq_RBA}{\text{RBA}}, there exists a regular bounded open set $U$ such that $15B\subseteq U\subseteq 16B$. For any $\phi\in\mathcal{F}(U)$, by \cite[Lemma 4.1]{Yan25c}, we have
$$\int_{U}|\phi|^p\mathrm{d} m\le C_2\mathcal{E}(\phi),$$
where $C_2$ is the positive constant appearing therein, hence
$$\lvert \int_{B}\phi \mathrm{d}m\rvert\le\int_{B}\lvert \phi\rvert \mathrm{d}m\le\int_U \lvert \phi\rvert \mathrm{d}m\le m(U)^{1-\frac{1}{p}}\left(\int_{U}|\phi|^p \mathrm{d} m\right)^{\frac{1}{p}}\le C_2^{\frac{1}{p}}m(U)^{1-\frac{1}{p}}\mathcal{E}(\phi)^{\frac{1}{p}},$$
which gives $\phi\mapsto \int_{B}\phi \mathrm{d} m$ is a bounded linear functional on the Banach space $(\mathcal{F}(U),\mathcal{E})$. By \cite[Theorem 3.24]{KS24a}, there exists $v\in\mathcal{F}(U)$ such that $\mathcal{E}(v;\phi)=\int_{B}\phi \mathrm{d} m$ for any $\phi\in\mathcal{F}(U)$, that is, $v\in \mathcal{F}(U)$ is superharmonic in $U$ and $\mu[v]=1_{B}m$ in $U$, in particular, $v$ is harmonic in $U\backslash \overline{B}$. By \ref{eq_EHI}, we have $v\in C(U\backslash \overline{B})$. By Proposition \ref{prop_comparison}, we have $v$ is non-negative in $U$ and $u\le v$ in $B$. For any Lebesgue point $x\in B$ of $v$, by the upper bound in \ref{eq_Wolff}, we have
$$v(x)\le C_1\left(\einf_{B(x,3R)}v+\mathcal{W}^{\mu[v]}(x,6R)\right),$$
where, $\mathcal{W}^{\mu[v]}(x,6R)\le \mathcal{W}^m(x,6R)\asymp \Psi(R)^{\frac{1}{p-1}}$. Let $t=\einf_{B(x,3R)}v$, without loss of generality, we may assume that $t>0$, then $\widetilde{v}\ge t$ q.e. on $B(x,3R)\supseteq 2B\supseteq \overline{B}$. Since $U$ is regular and $v\in \mathcal{F}(U)$ is bounded and harmonic in $U\backslash \overline{14B}$, by (\ref{eq_cap_infchar}), we have
$$\mathrm{cap}\left(B,X\backslash 16B\right)\le\mathrm{cap}\left(\overline{B},X\backslash U\right)\le \mu \left[\frac{v}{t}\right](U)=\frac{1}{t^{p-1}}m(B).$$
Hence
$$t\le \left(\frac{m(B)}{\mathrm{cap}\left(B,X\backslash 16B\right)}\right)^{\frac{1}{p-1}}\overset{(*)}{\scalebox{2}[1]{$\lesssim$}}\left(\frac{m(B)}{\frac{m(B)}{\Psi(R)}}\right)^{\frac{1}{p-1}}=\Psi(R)^{\frac{1}{p-1}},$$
where $(*)$ follows from \ref{eq_VD}, \hyperlink{eq_cap}{$\text{cap}(\Psi)$}. Therefore, we have $\esup_{B}u\le\esup_{B}v\lesssim\Psi(R)^{\frac{1}{p-1}}$.
\end{proof}

\begin{proposition}\label{prop_Poi_lambda}
Assume \ref{eq_VD}, \hyperlink{eq_LSC}{\text{LSC}}, \ref{eq_BSP}, \ref{eq_EHI}, \ref{eq_EHIann}, \hyperlink{eq_cap}{$\text{cap}(\Psi)$}. Then there exist $C_1,C_2>0$ such that for any $\lambda>0$ and any ball $B=B(x_0,R)$, for any $u\in\mathcal{F}(B)$ satisfying that $-\Delta_pu+{\lambda}|u|^{p-2}u=1$ in $B$, that is,
$$\mathcal{E}(u;\varphi)+\lambda\int_X|u|^{p-2}u\varphi\mathrm{d} m=\int_X\varphi\mathrm{d} m\text{ for any }\varphi\in\mathcal{F}(B),$$
we have
\begin{align}
\esup_Bu&\le\frac{1}{\lambda^{\frac{1}{p-1}}},\nonumber\\
\einf_{\frac{1}{2}B}u&\ge\frac{C_1}{(C_2+\lambda\Psi(R))^{\frac{1}{p-1}}}{\Psi(R)^{\frac{1}{p-1}}}.\label{eq_lambda_lbd}
\end{align}
\end{proposition}

\begin{remark}
This result plays a role analogous to that of the estimate in \cite[LEMMA 3.1]{GHL15} in the case $p=2$.
\end{remark}

\begin{proof}
By Proposition \ref{prop_comparison}, we have $0\le u\le \frac{1}{\lambda^{\frac{1}{p-1}}}$ in $B$, thus the upper bound holds. It remains to prove the lower bound. By Proposition \ref{prop_EF} and its proof, there exists non-negative $v\in\mathcal{F}(B)$ such that $-\Delta_pv=1$ in $B$, $\esup_Bv\le C_3\Psi(R)^{\frac{1}{p-1}}$, and $\einf_{\frac{1}{2}B}v\ge \frac{1}{C_3}\Psi(R)^{\frac{1}{p-1}}$, where $C_3$ is the positive constant appearing therein. For any $a>0$, we have
\begin{align*}
&-\Delta_p(av)+\lambda|av|^{p-2}(av)=a^{p-1}\left(-\Delta_pv+\lambda|v|^{p-2}v\right)\le a^{p-1}(1+\lambda C_3^{p-1}\Psi(R))\text{ in }B.
\end{align*}
Let $a=\frac{1}{\left(1+ \lambda C_3^{p-1}\Psi(R)\right)^{\frac{1}{p-1}}}$, then
$$-\Delta_p(av)+\lambda|av|^{p-2}(av)\le1=-\Delta_p u+\lambda|u|^{p-2}u\text{ in }B.$$
Since $u,v\in\mathcal{F}(B)$, by Proposition \ref{prop_comparison}, we have $av \le u$ in $B$, in particular, we have
$$\einf_{\frac{1}{2}B}u\ge a\einf_{\frac{1}{2}B}v\ge \frac{1}{\left(1+ \lambda C_3^{p-1}\Psi(R)\right)^{\frac{1}{p-1}}}\frac{1}{C_3}\Psi(R)^{\frac{1}{p-1}}=\frac{C_1}{\left(C_2+\lambda\Psi(R)\right)^{\frac{1}{p-1}}}\Psi(R)^{\frac{1}{p-1}},$$
where $C_1=\frac{1}{C_3^2}$, $C_2=\frac{1}{C_3^{p-1}}$, which gives (\ref{eq_lambda_lbd}).
\end{proof}

We give the proof of Proposition \ref{prop_main_CS} as follows.

\begin{proof}[Proof of Proposition \ref{prop_main_CS}]
The argument is based on the proof of \cite[THEOREM 3.2]{GHL15}. For notational convenience, we may assume that all functions in $\mathcal{F}$ are quasi-continuous. For any ball $B=B(x_0,R)$, we construct a desired cutoff function for $B\subseteq 2B$.

Let $\Omega=2B\backslash\overline{B}$. For any $\varphi\in\mathcal{F}(\Omega)$, we have
$$\lvert\int_X\varphi\mathrm{d} m\rvert\le m(\Omega)^{1-\frac{1}{p}}\left(\int_X|\varphi|^p\mathrm{d} m\right)^{\frac{1}{p}}\le m(\Omega)^{1-\frac{1}{p}}\Psi(R)^{\frac{1}{p}}\mathcal{E}_{\frac{1}{\Psi(R)}}(\varphi)^{\frac{1}{p}},$$
hence $\varphi\mapsto\int_X\varphi\mathrm{d} m$ is a bounded linear functional on $(\mathcal{F}(\Omega),\mathcal{E}_{\frac{1}{\Psi(R)}}^{1/p})$. Since $(\mathcal{F}(\Omega),\mathcal{E}_{\frac{1}{\Psi(R)}}^{1/p})$ is a Banach space, by \cite[Corollary 3.25]{KS24a}, there exists $u_\Omega\in\mathcal{F}(\Omega)$ such that $-\Delta_p u_\Omega+\frac{1}{\Psi(R)}|u_\Omega|^{p-2}u_\Omega=1$ in $\Omega$. By Proposition \ref{prop_comparison}, we have $0\le u_\Omega\le {\Psi(R)^{\frac{1}{p-1}}}$ in $\Omega$. Similarly, for any $x\in \frac{3}{2}B\backslash \frac{5}{4}B$, there exists $v\in\mathcal{F}(B(x,\frac{1}{8}R))$ such that $-\Delta_p v+\frac{1}{\Psi(R)}|v|^{p-2}v=1$ in $B(x,\frac{1}{8}R)$. Since $B(x,\frac{1}{8}R)\subseteq 2B\backslash\overline{B}$, by Proposition \ref{prop_comparison}, we have $u_\Omega\ge v$ in $B(x,\frac{1}{8}R)$, by Proposition \ref{prop_Poi_lambda}, we have
$$\einf_{B(x,\frac{1}{16}R)}u_\Omega\ge \einf_{B(x,\frac{1}{16}R)}v\ge \frac{C_1}{\left(C_2+\frac{1}{\Psi(R)}\Psi\left(\frac{R}{8}\right)\right)^{\frac{1}{p-1}}}\Psi\left(\frac{R}{8}\right)^{\frac{1}{p-1}}\ge C_3 \Psi(R)^{\frac{1}{p-1}},$$
where $C_1,C_2$ are the positive constants given in Proposition \ref{prop_Poi_lambda}, and $C_3>0$ depends only on $p,C_\Psi,C_1,C_2$. Hence
$$\einf_{\frac{3}{2}B\backslash \frac{5}{4}B}u_\Omega\ge C_3 \Psi(R)^{\frac{1}{p-1}}.$$

Let $v_\Omega=\frac{1}{C_3 \Psi(R)^{\frac{1}{p-1}}}u_\Omega$, then $v_\Omega\in\mathcal{F}(\Omega)$ satisfies that $0\le v_\Omega\le \frac{1}{C_3}$ in $\Omega$ and $v_\Omega\ge1$ in $\frac{3}{2}B\backslash \frac{5}{4}B$. Let
$$\phi=
\begin{cases}
1&\text{in }\frac{3}{2}B,\\
v_\Omega\wedge1&\text{on }X\backslash \frac{3}{2}B,
\end{cases}
$$
then $\phi\in\mathcal{F}$ is a cutoff function for $B\subseteq 2B$.

Note that by the Markovian property and the chain rule, for any $f,g\in \mathcal{F}\cap L^\infty(X;m)$, for any locally Lipschitz function $\Phi:\mathbb{R}\to \mathbb{R}$ with $\Phi(0)=0$, we have $\Phi(f),fg\in \mathcal{F}\cap L^\infty(X;m)$.

Firstly, for any $f\in \mathcal{F}\cap C_c(X)$, by \cite[Proposition 4.16]{KS24a}, we have
$$\int_{X}|f|^p\mathrm{d}\Gamma(v_\Omega)=\mathcal{E}(v_\Omega;v_\Omega|f|^p)-\left(\frac{p-1}{p}\right)^{p-1}\mathcal{E}(v_\Omega^{\frac{p}{p-1}};|f|^p).$$
Since $v_\Omega|f|^p\in \mathcal{F}(\Omega)$, we have
\begin{align*}
&\mathcal{E}(v_\Omega;v_\Omega|f|^p)\\
&=\frac{1}{C_3^{p-1}\Psi(R)}\left(\mathcal{E}(u_\Omega;v_\Omega|f|^p)+\frac{1}{\Psi(R)}\int_X|u_\Omega|^{p-2}u_\Omega(v_\Omega|f|^p)\mathrm{d} m\right.\\
&\hspace{80pt}\left.-\frac{1}{\Psi(R)}\int_X|u_\Omega|^{p-2}u_\Omega(v_\Omega|f|^p)\mathrm{d} m\right)\\
&\le\frac{1}{C_3^{p-1}\Psi(R)}\int_Xv_\Omega|f|^p\mathrm{d} m\le\frac{1}{C_3^{p}\Psi(R)}\int_{2B}|f|^p\mathrm{d} m.
\end{align*}
Since
\begin{align*}
&\lvert\mathcal{E}(v_\Omega^{\frac{p}{p-1}};|f|^p)\rvert\overset{(*)}{\scalebox{2}[1]{$=$}}\lvert\int_X \left(\frac{p}{p-1}v_\Omega^{\frac{p}{p-1}-1}\right)^{p-1}\left(p|f|^{p-1}\right)\mathrm{d}\Gamma(v_\Omega;|f|)\rvert\\
&=p \left(\frac{p}{p-1}\right)^{p-1}\lvert\int_Xv_\Omega|f|^{p-1}\mathrm{d}\Gamma(v_\Omega;|f|)\rvert\\
&\overset{(**)}{\scalebox{2}[1]{$\le$}}p \left(\frac{p}{p-1}\right)^{p-1}\left(\int_X|f|^{(p-1)\frac{p}{p-1}}\mathrm{d}\Gamma(v_\Omega)\right)^{\frac{p-1}{p}}\left(\int_Xv_\Omega^p\mathrm{d}\Gamma(|f|)\right)^{\frac{1}{p}}\\
&\le p \left(\frac{p}{p-1}\right)^{p-1}\left(\int_X|f|^{p}\mathrm{d}\Gamma(v_\Omega)\right)^{\frac{p-1}{p}}\frac{1}{C_3}\left(\int_{2B}\mathrm{d}\Gamma(f)\right)^{\frac{1}{p}},
\end{align*}
where $(*)$ follows from \cite[Theorem 7.6]{Sas26} and $(**)$ follows from \cite[Proposition 4.8]{KS24a}, we have
\begin{align*}
&\int_X|f|^p\mathrm{d}\Gamma(v_\Omega)\\
&\le\frac{1}{C_3^{p}\Psi(R)}\int_{2B}|f|^p\mathrm{d} m+\frac{p}{C_3} \left(\int_X|f|^{p}\mathrm{d}\Gamma(v_\Omega)\right)^{\frac{p-1}{p}}\left(\int_{2B}\mathrm{d}\Gamma(f)\right)^{\frac{1}{p}}\\
&\le\frac{1}{C_3^{p}\Psi(R)}\int_{2B}|f|^p\mathrm{d} m+\frac{1}{2}\int_X|f|^{p}\mathrm{d}\Gamma(v_\Omega)+C_4\int_{2B}\mathrm{d}\Gamma(f),
\end{align*}
where in the last inequality, we use Young's inequality, and $C_4>0$ depends only on $p,C_3$. In summary
$$\int_{2B}|f|^p\mathrm{d}\Gamma(\phi)\le\int_{X}|f|^p\mathrm{d}\Gamma(v_\Omega)\le 2C_4\int_{2B}\mathrm{d}\Gamma(f)+\frac{2}{C_3^{p}\Psi(R)}\int_{2B}|f|^p\mathrm{d} m.$$

Secondly, for any $f\in \mathcal{F}$, by \cite[Proposition 8.5, Proposition 8.12]{Yan25a}, there exists $\{f_n\}\subseteq \mathcal{F}\cap C_c(X)$ such that $\{f_n\}$ is $\mathcal{E}_1$-convergent to $f$ and $\{f_n\}$ converges to $f$ q.e. on $X$, which is also $\Gamma(\phi)$-a.e. on $X$, then by Fatou's lemma, we have
\begin{align*}
&\int_{2B}|f|^p\mathrm{d}\Gamma(\phi)\le\varliminf_{n\to+\infty}\int_{2B}|f_n|^p\mathrm{d}\Gamma(\phi)\\
&\le\varliminf_{n\to+\infty}\left(2C_4\int_{2B}\mathrm{d}\Gamma(f_n)+\frac{2}{C_3^{p}\Psi(R)}\int_{2B}|f_n|^p\mathrm{d} m\right)\\
&=2C_4\int_{2B}\mathrm{d}\Gamma(f)+\frac{2}{C_3^{p}\Psi(R)}\int_{2B}|f|^p\mathrm{d} m.
\end{align*}
\end{proof}

To conclude this section, we make the following remark. By Theorem \ref{thm_main_conv}, whose proof will be given in Section \ref{sec_EHI}, the conditions \ref{eq_VD}, \ref{eq_PI}, \ref{eq_CS} ensure that the assumptions of all the results in this section are satisfied, and hence all these results hold. In particular, we obtain the following corollary.

\begin{corollary}
Assume \ref{eq_VD}, \ref{eq_PI}, \ref{eq_CS}. Then the \hyperlink{eq_RBA}{\text{RBA}} condition of Proposition \ref{prop_RBA} and the Wolff potential estimate \ref{eq_Wolff} of Theorem \ref{thm_Wolff} hold.
\end{corollary}

\section{Proof of Theorem \ref{thm_main_conv}}\label{sec_EHI}

The proof of \ref{eq_EHI} is classical and relies on the standard Nash--Moser--De Giorgi iteration technique, with \ref{eq_CS} supplying the necessary cutoff functions. While \hyperlink{eq_LSC}{\text{LSC}}, \ref{eq_EHIann} follow as byproducts of the mean value inequality and the weak Harnack inequality. We divide the proof into the following steps.

For the first step, we have the following Sobolev inequality. The proof is standard; see \cite{BCLS95,BM95,HK95,MS95,CG98,CHS01,Sal02}, and is therefore omitted.

\begin{proposition}[Sobolev inequality]\label{prop_Sobolev}
Assume \ref{eq_VD}, \ref{eq_PI}. Then there exist $\kappa>1$, $C>0$ such that for any ball $B(x_0,R)$, for any $f\in\mathcal{F}(B(x_0,R))$, we have
$$\left(\int_{B(x_0,R)}|f|^{p\kappa}\mathrm{d} m\right)^{\frac{1}{\kappa}}\le C\frac{\Psi(R)}{V(x_0,R)^{\frac{\kappa-1}{\kappa}}}\mathcal{E}(f).$$
Indeed, we can take $\kappa=\frac{\nu}{\nu-\beta_*}$, where $\nu=\max\{\beta_*+1,\log_2C_{VD}\}$.
\end{proposition}

For the second step, we have the following Caccioppoli inequality.

\begin{proposition}[Caccioppoli inequality]\label{prop_Cacci}
Assume \ref{eq_VD}, \ref{eq_CS}. Then there exists $C>0$ such that for any $x_0\in X$, for any $R,r>0$, there exists a cutoff function $\phi\in\mathcal{F}$ for $B(x_0,R)\subseteq B(x_0,R+r)$ such that for any $u\in\mathcal{F}$ which is non-negative bounded subharmonic in $B(x_0,R+r)$, for any $\theta\ge0$, we have
$$\int_{B(x_0,R+r)}\mathrm{d}\Gamma(\phi(u-\theta)_+)\le \frac{C}{\Psi(r)}\int_{B(x_0,R+r)}u^p\mathrm{d} m.$$
\end{proposition}

\begin{proof}
For notational convenience, we may assume that all functions in $\mathcal{F}$ are quasi-contin-uous. Let $v=(u-\theta)_+\in\mathcal{F}$. Let $\phi\in\mathcal{F}$ be a cutoff function for $B(x_0,R)\subseteq B(x_0,R+r)$ chosen later. Since $\phi\in \mathcal{F}(B(x_0,R+r))$ and $v\in \mathcal{F}$ are both non-negative bounded in $B(x_0,R+r)$, we have $\phi^pv\in\mathcal{F}(B(x_0,R+r))$ is non-negative. Since $u\in\mathcal{F}$ is subharmonic in $B(x_0,R+r)$, we have $\mathcal{E}(u;\phi^pv)\le0$, that is,
$$0\ge\int_X\mathrm{d}\Gamma(u;\phi^pv)=\int_X\mathrm{d}\Gamma(v;\phi^pv)=\int_X\phi^p\mathrm{d}\Gamma(v)+p\int_X\phi^{p-1}v\mathrm{d}\Gamma(v;\phi).$$
Hence
$$\int_X\phi^p\mathrm{d}\Gamma(v)\le p\lvert\int_X\phi^{p-1}v\mathrm{d}\Gamma(v;\phi)\rvert\le p \left(\int_X\phi^p\mathrm{d}\Gamma(v)\right)^{\frac{p-1}{p}}\left(\int_Xv^p\mathrm{d}\Gamma(\phi)\right)^{\frac{1}{p}},$$
which gives
$$\int_X\phi^p\mathrm{d}\Gamma(v)\le p^p\int_Xv^p\mathrm{d}\Gamma(\phi).$$
Let $\delta>0$ be chosen later, by the self-improvement property of \ref{eq_CS} (see \cite[Proposition 3.1]{Yan25c}), there exist $C_\delta>0$ and a cutoff function $\phi\in\mathcal{F}$ for $B(x_0,R)\subseteq B(x_0,R+r)$ such that
\begin{align*}
&\int_Xv^p\mathrm{d}\Gamma(\phi)=\int_{B(x_0,R+r)\backslash\overline{B(x_0,R)}}v^p\mathrm{d}\Gamma(\phi)\\
&\le\delta\int_{B(x_0,R+r)\backslash\overline{B(x_0,R)}}\phi^p\mathrm{d}\Gamma(v)+\frac{C_\delta}{\Psi(r)}\int_{B(x_0,R+r)\backslash\overline{B(x_0,R)}}v^p\mathrm{d} m\\
&\le\delta\int_X\phi^p\mathrm{d}\Gamma(v)+\frac{C_\delta}{\Psi(r)}\int_{B(x_0,R+r)}u^p\mathrm{d} m,
\end{align*}
hence
$$\int_X\phi^p\mathrm{d}\Gamma(v)\le p^p \left(\delta\int_X\phi^p\mathrm{d}\Gamma(v)+\frac{C_\delta}{\Psi(r)}\int_{B(x_0,R+r)}u^p\mathrm{d} m\right).$$
Let $\delta=\frac{1}{2p^p}$, then
$$\int_{B(x_0,R+r)}\phi^p\mathrm{d}\Gamma(v)=\int_X\phi^p\mathrm{d}\Gamma(v)\le \frac{2p^pC_\delta}{\Psi(r)}\int_{B(x_0,R+r)}u^p\mathrm{d} m.$$
By \cite[Lemma 6.1]{Yan25c}, we have
\begin{align*}
&\int_{B(x_0,R+r)}\mathrm{d}\Gamma(\phi v)\\
&\le2^{p-1}\left(\int_{B(x_0,R+r)}\phi^p\mathrm{d}\Gamma(v)+\int_{B(x_0,R+r)}v^p\mathrm{d}\Gamma(\phi)\right)\le\frac{C}{\Psi(r)}\int_{B(x_0,R+r)}u^p\mathrm{d} m,
\end{align*}
where $C=2^{p-1}\left(2(1+\delta)p^p+1\right)C_\delta$.
\end{proof}

We will need the following elementary iteration lemma for the subsequent analysis.

\begin{lemma}[{\cite[LEMMA 4.9]{CKW21}}]\label{lem_iteration}
Let $\beta>0$ and $\{A_j\}_{j\ge0}$ a sequence of positive numbers satisfying that
$$A_{j+1}\le c_0 b^jA_j^{1+\beta}\text{ for any }j\ge0,$$
where $c_0>0$ and $b>1$ are some constants. If
$$A_0\le c_0^{-1/\beta}b^{-1/\beta^2},$$
then
$$A_j\le b^{-j/\beta}A_0\text{ for any }j\ge0,$$
which implies $\lim_{j\to+\infty}A_j=0$.
\end{lemma}

We have the following property of subharmonic functions.

\begin{lemma}\label{lem_subhar}
Let $U$ be a bounded open set and $u\in\mathcal{F}$ bounded from above subharmonic in $U$. Then for any $\theta\ge0$, we have $(u-\theta)_+\in\mathcal{F}$ is bounded subharmonic in $U$.
\end{lemma}

\begin{proof}
The proof is essentially the same as that of \cite[PROPOSITION 2.1]{GHL15}, which was established for the case $p=2$. For general $p>1$, the argument remains valid upon noting that for any $f\in C^2(\mathbb{R})$ with $f''\ge0$, $1\ge f'\ge0$, $f|_{(-\infty,0]}=0$, for any non-negative bounded $\phi\in\mathcal{F}(U)$, we have $f(u)$, $f'(u)\in \mathcal{F}$ and
\begin{align*}
&\mathcal{E}(f(u);\phi)=\int_X\mathrm{d}\Gamma(f(u);\phi)=\int_Xf'(u)^{p-1}\mathrm{d}\Gamma(u;\phi)\\
&=\int_X\mathrm{d}\Gamma(u;\phi f'(u)^{p-1})-\int_X\phi\mathrm{d}\Gamma(u;f'(u)^{p-1})\\
&=\mathcal{E}(u;\phi f'(u)^{p-1})-(p-1)\int_X\phi f'(u)^{p-2}f''(u)\mathrm{d}\Gamma(u)\\
&\le\mathcal{E}(u;\phi f'(u)^{p-1})\le0,
\end{align*}
after which the remainder of the argument proceeds identically as in \cite{GHL15}.
\end{proof}

For the third step, we have the $L^p$-mean value inequality for subharmonic functions as follows. The proof follows the same argument as in \cite[THEOREM 6.2]{GHL15}.

\begin{proposition}[$L^p$-mean value inequality]\label{prop_MV}
Assume \ref{eq_VD},\ref{eq_PI}, \ref{eq_CS}. Then there exists $C>0$ such that for any ball $B(x_0,R)$, for any $u\in\mathcal{F}$ which is non-negative bounded subharmonic in $B(x_0,R)$, we have
$$\esup_{B(x_0,\frac{1}{2}R)}u^p\le \frac{C}{V(x_0,R)}\int_{B(x_0,R)}u^p\mathrm{d} m.$$
\end{proposition}

\begin{proof}
Let $\theta>0$ be chosen later. For any $k\ge0$, let $R_k=\left(\frac{1}{2}+\frac{1}{2^{k+1}}\right)R$ and $\theta_k=(1-\frac{1}{2^k})\theta$, then $R=R_0>R_1>\ldots>R_k\downarrow \frac{1}{2}R$ and $0=\theta_0<\theta_1<\ldots<\theta_k\uparrow\theta$. For any $k\ge0$, by Lemma \ref{lem_subhar}, $(u-\theta_k)_+\in\mathcal{F}$ is subharmonic in $B(x_0,R)$, since $(u-\theta_{k+1})_+=((u-\theta_k)_+-(\theta_{k+1}-\theta_k))_+$, by Proposition \ref{prop_Cacci}, there exists a cutoff function $\phi_k\in\mathcal{F}$ for $B(x_0,R_{k+1})\subseteq B(x_0,R_k)$ such that
$$\int_{B(x_0,R_k)}\mathrm{d}\Gamma(\phi_k(u-\theta_{k+1})_+)\le \frac{C_1}{\Psi(R_k-R_{k+1})}\int_{B(x_0,R_k)}(u-\theta_k)_+^p\mathrm{d} m,$$
where $C_1$ is the positive constant appearing therein, since $\phi_k(u-\theta_{k+1})_+\in\mathcal{F}(B(x_0,R_k))$, by Proposition \ref{prop_Sobolev}, we have
\begin{align*}
&\left(\int_{B(x_0,R_k)}|\phi_k(u-\theta_{k+1})_+|^{p\kappa}\mathrm{d} m\right)^{\frac{1}{\kappa}}\le\frac{C_2}{V(x_0,R_k)^{\frac{\kappa-1}{\kappa}}}\Psi(R_k)\int_{B(x_0,R_k)}\mathrm{d}\Gamma(\phi_k(u-\theta_{k+1})_+),
\end{align*}
where $C_2$ is the positive constant appearing therein. Hence
\begin{align*}
&\int_{B(x_0,R_{k+1})}(u-\theta_{k+1})_+^p\mathrm{d} m\le\int_{B(x_0,R_k)}|\phi_k(u-\theta_{k+1})_+|^{p}1_{\{u>\theta_{k+1}\}}\mathrm{d} m\\
&\le\left(\int_{B(x_0,R_k)}|\phi_k(u-\theta_{k+1})_+|^{p\kappa}\mathrm{d} m\right)^{\frac{1}{\kappa}}\left(\int_{B(x_0,R_k)}1_{\{u>\theta_{k+1}\}}\mathrm{d} m\right)^{1-\frac{1}{\kappa}}\\
&\le \frac{C_1C_2}{V(x_0,R_k)^{\frac{\kappa-1}{\kappa}}}\frac{\Psi(R_k)}{\Psi(R_k-R_{k+1})}\frac{1}{(\theta_{k+1}-\theta_k)^{p(1-\frac{1}{\kappa})}}\left(\int_{B(x_0,R_k)}(u-\theta_k)_+^p\mathrm{d} m\right)^{1+1-\frac{1}{\kappa}}.
\end{align*}
For any $k\ge0$, let
$$I_k=\int_{B(x_0,R_k)}(u-\theta_k)_+^p\mathrm{d} m,$$
recall that $\kappa=\frac{\nu}{\nu-\beta_*}$, then
\begin{align*}
&I_{k+1}\le \frac{C_3}{V(x_0,R)^{\frac{\beta_*}{\nu}}\theta^{\frac{p\beta_*}{\nu}}}2^{(\beta^*+\frac{p\beta_*}{\nu})k}I_{k}^{1+\frac{\beta_*}{\nu}},
\end{align*}
where $C_3=2^{2\beta^*+\frac{p\beta_*}{\nu}}C_1C_2C_\Psi C_{VD}^{\frac{\beta_*}{\nu}}$. By Lemma \ref{lem_iteration}, if
$$\int_{B(x_0,R)}u^p\mathrm{d} m=I_0\le \left(\frac{C_3}{V(x_0,R)^{{\frac{\beta_*}{\nu}}}\theta^{\frac{p\beta_*}{\nu}}}\right)^{-{1}/\left({\frac{\beta_*}{\nu}}\right)}2^{-(\beta^*+\frac{p\beta_*}{\nu})/(\frac{\beta_*}{\nu})^2}=\frac{1}{C_4}V(x_0,R)\theta^p,$$
where $C_4=C_3^{{1}/\left({\frac{\beta_*}{\nu}}\right)}2^{(\beta^*+\frac{p\beta_*}{\nu})/(\frac{\beta_*}{\nu})^2}$, then
$$\int_{B(x_0,\frac{1}{2}R)}(u-\theta)_+^p\mathrm{d} m\le\lim_{k\to+\infty}I_k=0.$$
Hence, let $\theta^p=\frac{C_4}{V(x_0,R)}I_0$, then $\esup_{B(x_0,\frac{1}{2}R)}u\le \theta$, that is,
$$\esup_{B(x_0,\frac{1}{2}R)}u^p\le \frac{C_4}{V(x_0,R)}\int_{B(x_0,R)}u^p\mathrm{d} m.$$
\end{proof}

For the fourth step, as a corollary, we have the $L^q$-mean value inequality for subharmonic functions for \emph{all} $q>0$ as follows. The proof follows the same argument as that of $(MV)_1$ in \cite{GHL15}, using $(MV)_2$ therein, as in the proofs of \cite[LEMMAS 9.1 and 9.2]{GHL15}, and is therefore omitted.

\begin{corollary}[$L^q$-mean value inequality]\label{cor_MVq}
Assume \ref{eq_VD},\ref{eq_PI}, \ref{eq_CS}. Then for any $q>0$, there exists $C>0$ such that for any ball $B(x_0,R)$, for any $u\in\mathcal{F}$ which is non-negative bounded subharmonic in $B(x_0,R)$, we have
$$\esup_{B(x_0,\frac{1}{4}R)}u^q\le \frac{C}{V(x_0,R)}\int_{B(x_0,R)}u^q\mathrm{d} m.$$
\end{corollary}

As a byproduct, we present the proof of \hyperlink{eq_LSC}{\text{LSC}} as follows; see also \cite[Theorem 8.22]{BB11} and \cite[LEMMA 3.4]{TW02AJM}.

\begin{proof}[Proof of \hyperlink{eq_LSC}{\text{LSC}} in Theorem \ref{thm_main_conv}]
Let $U$ be a bounded open subset of $X$ and $u\in \mathcal{F}$ which is bounded from below superharmonic in $U$. For any $x\in U$, let
$$u^*(x)=\lim_{r\downarrow0}\einf_{B(x,r)}u.$$
It is obvious that $u^*$ is lower semi-continuous in $U$. We only need to show that for any Lebesgue point $x_0\in U$ of $u$, we have $u^*(x_0)=u(x_0)$. It is obvious that $u^*(x_0)\le u(x_0)$. On the other hand, by Lemma \ref{lem_subhar}, we have $(u(x_0)-u)_+\in \mathcal{F}$ is bounded subharmonic in $U$. For notational convenience, we write $(u(x_0)-u)_+$ to denote $(u(x_0)\psi-u)_+$, where $\psi\in \mathcal{F}\cap C_c(X)$ satisfies $\psi=1$ on $\overline{U}$. By Corollary \ref{cor_MVq}, there exists $C>0$ such that for any $r>0$ with $B(x_0,r)\subseteq U$, we have
\begin{align*}
u(x_0)-\einf_{B(x_0,\frac{1}{4}r)}u\le\esup_{B(x_0,\frac{1}{4}r)}\left(u(x_0)-u\right)_+\le \frac{C}{V(x_0,r)}\int_{B(x_0,r)}\lvert u(x_0)-u\rvert\mathrm{d}m,
\end{align*}
hence
$$\einf_{B(x_0,\frac{1}{4}r)}u\ge u(x_0)-\frac{C}{V(x_0,r)}\int_{B(x_0,r)}\lvert u(x_0)-u\rvert\mathrm{d}m.$$
Letting $r\downarrow0$, we have $u^*(x_0)\ge u(x_0)$.
\end{proof}

For the fifth step, we have the following lemma of growth for superharmonic functions. We follow the idea of the proof given in \cite[LEMMA 3.5]{HY23} and \cite[Lemma 4.1]{GHH18}.

\begin{proposition}[Lemma of growth]\label{prop_LG}
Assume \ref{eq_VD}, \ref{eq_PI}, \ref{eq_CS}. Then there exists $\varepsilon_0\in(0,1)$ such that for any $\varepsilon\in(\varepsilon_0,1)$, there exists $\delta=\delta_\varepsilon\in(0,1)$ such that for any ball $B=B(x_0,R)$, for any $u\in\mathcal{F}$ which is non-negative superharmonic in $2B$, for any $a>0$, if
$$\frac{m(B\cap\{u\ge a\})}{m(B)}\ge\varepsilon,$$
then
$$\einf_{\frac{1}{2}B}u\ge\delta a.$$
\end{proposition}

\begin{proof}
For any $r>0$, let $B_r=B(x_0,r)$. For any $r_1,r_2\in(\frac{1}{2}R,R]$ with $r_1<r_2$, for any $a,b>0$ with $a<b$, let
$$I=\frac{m(B_{r_1}\cap\{u<a\})}{m(B_{r_1})}, J=\frac{m(B_{r_2}\cap\{u<b\})}{m(B_{r_2})}.$$
Since $u\in\mathcal{F}$ is non-negative superharmonic in $2B$, we have $b-u\in\mathcal{F}$ is bounded from above subharmonic in $2B$. For notational convenience, we write $b-u$ to denote $b\psi-u$, where $\psi\in\mathcal{F}\cap C_c(X)$ satisfies $\psi=1$ on $\overline{2B}$. By Lemma \ref{lem_subhar}, we have $(b-u)_+\in\mathcal{F}$ is bounded subharmonic in $2B$. By Proposition \ref{prop_Cacci}, there exists a cutoff function $\phi\in\mathcal{F}$ for $B_{r_1}\subseteq B_{r_2}$ such that
$$\int_{B_{r_2}}\mathrm{d}\Gamma(\phi(b-u)_+)\le \frac{C_1}{\Psi(r_2-r_1)}\int_{B_{r_2}}(b-u)_+^p\mathrm{d} m\le\frac{C_1}{\Psi(r_2-r_1)}b^pm(B_{r_2}\cap\{u<b\}),$$
where $C_1$ is the positive constant appearing therein. By Proposition \ref{prop_Sobolev}, we have
\begin{align*}
&\left(\int_{B_{r_2}}|\phi(b-u)_+|^{p\kappa}\mathrm{d} m\right)^{\frac{1}{\kappa}}\le \frac{C_2}{m(B_{r_2})^{\frac{\kappa-1}{\kappa}}}\Psi(r_2)\int_{B_{r_2}}\mathrm{d}\Gamma(\phi(b-u)_+),
\end{align*}
where $C_2$ is the positive constant appearing therein. Hence
\begin{align*}
&\int_{B_{r_1}}(b-u)_+^p\mathrm{d} m\le\int_{B_{r_2}}|\phi(b-u)_+|^p1_{\{u<b\}}\mathrm{d} m\\
&\le\left(\int_{B_{r_2}}|\phi(b-u)_+|^{p\kappa}\mathrm{d} m\right)^{\frac{1}{\kappa}}m(B_{r_2}\cap\{u<b\})^{1-\frac{1}{\kappa}}\\
&\le\frac{C_1C_2}{m(B_{r_2})^{\frac{\kappa-1}{\kappa}}}\frac{\Psi(r_2)}{\Psi(r_2-r_1)}b^pm(B_{r_2}\cap\{u<b\})^{1+1-\frac{1}{\kappa}},
\end{align*}
where
$$\int_{B_{r_1}}(b-u)_+^p\mathrm{d} m\ge \int_{B_{r_1}\cap\{u<a\}}(b-u)_+^p\mathrm{d} m\ge (b-a)^p m(B_{r_1}\cap\{u<a\}),$$
which gives
\begin{align*}
&I\le C_1C_2C_{VD}C_\Psi \left(\frac{r_2}{r_2-r_1}\right)^{\beta^*}\left(\frac{b}{b-a}\right)^pJ^{1+1-\frac{1}{\kappa}}.
\end{align*}

Let $\delta\in(0,1)$ be chosen later. For any $k\ge0$, let $R_k=\left(\frac{1}{2}+\frac{1}{2^{k+1}}\right)R$ and $a_k=(\delta+\frac{1-\delta}{2^k})a$, then $R=R_0>R_1>\ldots>R_k\downarrow \frac{R}{2}$ and $a=a_0>a_1>\ldots>a_k\downarrow\delta a$, let
$$I_k=\frac{m(B_{R_k}\cap\{u<a_k\})}{m(B_{R_k})},$$
recall that $\kappa=\frac{\nu}{\nu-\beta_*}$, then
\begin{align*}
&I_{k+1}\le \frac{C_3}{(1-\delta)^p}2^{(p+\beta^*)k}I_{k}^{1+\frac{\beta_*}{\nu}},
\end{align*}
where $C_3=\max\{2^{p+2\beta^*}C_1C_2C_{VD}C_\Psi,1\}$. By Lemma \ref{lem_iteration}, if
$$I_0\le \left(\frac{C_3}{(1-\delta)^p}\right)^{-1/(\frac{\beta_*}{\nu})}\left(2^{p+\beta^*}\right)^{-1/(\frac{\beta_*}{\nu})^2}=\frac{(1-\delta)^{\frac{p\nu}{\beta_*}}}{C_4},$$
where $C_4=2^{\frac{(p+\beta^*)\nu^2}{\beta_*^2}}C_3^{\frac{\nu}{\beta_*}}>1$, then
$\lim_{k\to+\infty}I_k=0$. Let $\varepsilon_0=1-\frac{1}{C_4}\in(0,1)$, then for any $\varepsilon\in(\varepsilon_0,1)$, there exists $\delta\in(0,1)$ depending only on $\frac{p\nu}{\beta_*}$, $C_4$, $\varepsilon$, such that $\frac{(1-\delta)^{\frac{p\nu}{\beta_*}}}{C_4}=1-\varepsilon$. Since $I_0=\frac{m(B\cap\{u<a\})}{m(B)}\le 1-\varepsilon$, we have
$$\frac{m (B_{\frac{1}{2}R}\cap \{u<\delta a\})}{m (B_{\frac{1}{2}R})}\le\lim_{k\to+\infty}I_k=0,$$
that is, $\einf_{\frac{1}{2}B}u\ge\delta a$. 
\end{proof}

\begin{corollary}\label{cor_infu}
Assume \ref{eq_VD}, \ref{eq_PI}, \ref{eq_CS}. Then for any $q>0$, there exists $C>0$ such that for any ball $B$, for any $u\in\mathcal{F}$ which is non-negative superharmonic in $2B$, we have
$$\einf_{\frac{1}{2}B}u\ge \frac{1}{C}\left(\frac{1}{m(B)}\int_Bu^{-q}\mathrm{d} m\right)^{-1/q}.$$
\end{corollary}

\begin{proof}
Without loss of generality, we may assume that $\int_Bu^{-q}\mathrm{d} m<+\infty$. For any $a>0$, we have
$$m(B\cap\{u<a\})=m(B\cap\left\{\frac{1}{u}>\frac{1}{a}\right\})\le \frac{1}{\left(\frac{1}{a}\right)^q}\int_{B}\left(\frac{1}{u}\right)^q\mathrm{d} m=a^q\int_B u^{-q}\mathrm{d} m,$$
then
$$\frac{m(B\cap\{u\ge a\})}{m(B)}=1-\frac{m(B\cap\{u<a\})}{m(B)}\ge 1-\frac{a^q}{m(B)}\int_B u^{-q}\mathrm{d} m.$$
Let $a=\left(\frac{2}{(1-\varepsilon_0)m(B)}\int_Bu^{-q}\mathrm{d} m\right)^{-1/q}$, where $\varepsilon_0\in(0,1)$ is the constant in Proposition \ref{prop_LG}, then $\frac{m(B\cap\{u\ge a\})}{m(B)}\ge \frac{\varepsilon_0+1}{2}>\varepsilon_0$, hence there exists $\delta\in(0,1)$ such that
$$\einf_{\frac{1}{2}B}u\ge\delta a=\frac{\delta(1-\varepsilon_0)^{1/q}}{2^{1/q}}\left(\frac{1}{m(B)}\int_Bu^{-q}\mathrm{d} m\right)^{-1/q}.$$
\end{proof}

For the sixth step, we present the following results concerning BMO spaces, similar techniques were also employed in \cite{HY23}. Let $U$ be an open set and $u$ a locally integrable function in $U$. We define the semi-norm $\lVert {u}\rVert_{\mathrm{BMO}(U)}$ as
$$\lVert {u}\rVert_{\mathrm{BMO}(U)}=\sup \left\{\dashint_B|u-u_B|\mathrm{d} m:B\subseteq U\text{ is a ball}\right\}.$$
Let $\mathrm{BMO}(U)$ be the family of all locally integrable functions $u$ in $U$ with $\lVert {u}\rVert_{\mathrm{BMO}(U)}<+\infty$.

\begin{lemma}[John-Nirenberg inequality; {\cite[THEOREM 5.2]{ABKY11}}]\label{lem_JohnNirenberg}
Assume \ref{eq_VD}. Then there exist $C_1,C_2>0$ such that for any open set $U$ and any $u\in \mathrm{BMO}(U)$, for any ball $B$ with $12B\subseteq U$ and any $\lambda>0$, we have
$$m \left(B\cap\{|u-u_B|>\lambda\}\right)\le C_1m(B)\exp \left(-C_2 \frac{\lambda}{\lVert {u}\rVert_{\mathrm{BMO}(U)}}\right).$$
\end{lemma}

\begin{lemma}[{\cite[COROLLARY 5.6]{BM95AMPA}}]\label{lem_crossover}
Assume \ref{eq_VD}. Then for any ball $B(x_0,R)$ and any $u\in\mathrm{BMO}(B(x_0,R))$, for any ball $B$ with $12B\subseteq B(x_0,R)$, for any $b\ge \lVert {u}\rVert_{\mathrm{BMO}(B(x_0,R))}$, we have
$$\left\{\dashint_{B}\exp \left(\frac{C_2}{2b}u\right)\mathrm{d} m\right\}\left\{\dashint_{B}\exp \left(-\frac{C_2}{2b}u\right)\mathrm{d} m\right\}\le (C_1+1)^2,$$
where $C_1,C_2$ are the positive constants appearing in Lemma \ref{lem_JohnNirenberg}.
\end{lemma}

\begin{lemma}\label{lem_loguBMO}
Assume \ref{eq_VD}, \ref{eq_PI}, \ref{eq_ucap}. Then there exist $C>0$, $A>1$ such that for any ball $B(x_0,R)$, for any $u\in\mathcal{F}$ which is non-negative superharmonic in $B(x_0,AR)$, for any $\varepsilon>0$, we have $\log(u+\varepsilon)\in\mathrm{BMO}(B(x_0,R))$ and $\lVert {\log (u+\varepsilon)}\rVert_{\mathrm{BMO}(B(x_0,R))}\le C$.
\end{lemma}

\begin{proof}
Let $u\in\mathcal{F}$ be non-negative superharmonic in $B(x_0,A_{cap}A_{PI}R)$. By replacing $u$ with $u\vee 0$, we may assume that $u$ is non-negative in $X$. For any ball $B\subseteq B(x_0,R)$ with radius $r$, since $\log \frac{u+\varepsilon}{\varepsilon}\in\mathcal{F}$, by \ref{eq_PI}, we have
\begin{align*}
&\dashint_B|\log (u+\varepsilon)-(\log (u+\varepsilon))_B|\mathrm{d} m\le \left(\dashint_B|\log \frac{u+\varepsilon}{\varepsilon}-\left(\log \frac{u+\varepsilon}{\varepsilon}\right)_B|^p\mathrm{d} m\right)^{1/p}\\
&\le \left(\frac{C_{PI}\Psi(r)}{m(B)}\int_{A_{PI}B}\mathrm{d}\Gamma(\log \frac{u+\varepsilon}{\varepsilon})\right)^{1/p}=\left(\frac{C_{PI}\Psi(r)}{m(B)}\int_{A_{PI}B}\mathrm{d}\Gamma(\log {(u+\varepsilon)})\right)^{1/p}.
\end{align*}
By \ref{eq_ucap}, there exists a cutoff function $\phi\in\mathcal{F}$ for $A_{PI}B\subseteq A_{cap}A_{PI}B$ such that
$$\mathcal{E}(\phi)\le 2\mathrm{cap}(A_{PI}B,X\backslash A_{cap}A_{PI}B)\le 2C_{cap}\frac{m(A_{PI}B)}{\Psi(A_{PI}r)},$$
hence
\begin{align}
&\int_{A_{PI}B}\mathrm{d}\Gamma(\log (u+\varepsilon))\le\int_{X}\phi^p(u+\varepsilon)^{-p}\mathrm{d}\Gamma(u)=-\frac{1}{p-1}\int_X\phi^p\mathrm{d}\Gamma(u;(u+\varepsilon)^{1-p})\nonumber\\
&=-\frac{1}{p-1}\left(\mathcal{E}(u;\phi^p(u+\varepsilon)^{1-p})-p\int_X\phi^{p-1}(u+\varepsilon)^{1-p}\mathrm{d}\Gamma(u;\phi)\right)\nonumber\\
&\overset{(*)}{\scalebox{2}[1]{$\le$}} \frac{p}{p-1}\int_X\phi^{p-1}(u+\varepsilon)^{1-p}\mathrm{d}\Gamma(u;\phi)\nonumber\\
&\le\frac{p}{p-1}\left(\int_X\phi^p(u+\varepsilon)^{-p}\mathrm{d}\Gamma(u)\right)^{(p-1)/p}\left(\int_X\mathrm{d}\Gamma(\phi)\right)^{1/p},\label{eq_energy_log}
\end{align}
where in $(*)$, we use the facts that $u$ is superharmonic in $B(x_0,A_{cap}A_{PI}R)\supseteq A_{cap}A_{PI}B$, and $\phi^p(u+\varepsilon)^{1-p}\in \mathcal{F}(A_{cap}A_{PI}B)$, which gives
\begin{align*}
&\int_{A_{PI}B}\mathrm{d}\Gamma(\log (u+\varepsilon))\le\int_X\phi^p(u+\varepsilon)^{-p}\mathrm{d}\Gamma(u)\le 2\left(\frac{p}{p-1}\right)^pC_{cap}\frac{m(A_{PI}B)}{\Psi(A_{PI}r)}.
\end{align*}
Hence by \ref{eq_VD}, we have
\begin{align*}
&\dashint_B|\log (u+\varepsilon)-(\log (u+\varepsilon))_B|\mathrm{d} m\le \frac{p}{p-1}\left(2C_{VD}C_{PI}C_{cap}A_{PI}^{\log_2C_{VD}}\right)^{1/p}=C,
\end{align*}
which implies $\log(u+\varepsilon)\in\mathrm{BMO}(B(x_0,R))$ and $\lVert {\log (u+\varepsilon)}\rVert_{\mathrm{BMO}(B(x_0,R))}\le C$.
\end{proof}

For the seventh step, we have the weak Harnack inequality for superharmonic functions as follows.

\begin{proposition}[weak Harnack inequality]\label{prop_wEHI}
Assume \ref{eq_VD}, \ref{eq_PI}, \ref{eq_CS}. Then there exist $q>0$, $C>0$, $A>1$ such that for any ball $B(x_0,R)$, for any $u\in\mathcal{F}$ which is non-negative superharmonic in $B(x_0,AR)$, we have
$$\einf_{B(x_0,\frac{1}{2}R)}u^q\ge \frac{1}{CV(x_0,R)}\int_{B(x_0,R)}u^q\mathrm{d} m.$$
\end{proposition}

\begin{proof}
Let $C_1,C_2$ be the positive constants appearing in Lemma \ref{lem_crossover}, and $A_1, C_3$ the positive constants appearing in Lemma \ref{lem_loguBMO}. Let $A=12A_1$ and $u\in\mathcal{F}$ non-negative superharmonic in $B(x_0,AR)$. For any $\varepsilon>0$, by Lemma \ref{lem_loguBMO}, we have $\log(u+\varepsilon)\in \mathrm{BMO}(B(x_0,12R))$ and $\lVert {\log(u+\varepsilon)}\rVert_{\mathrm{BMO}(B(x_0,12R))}\le C_3$, by Lemma \ref{lem_crossover}, we have
$$\left(\dashint_{B(x_0,R)}\exp \left(\frac{C_2}{2C_3}\log (u+\varepsilon)\right)\mathrm{d} m\right)\left(\dashint_{B(x_0,R)}\exp \left(-\frac{C_2}{2C_3}\log (u+\varepsilon)\right)\mathrm{d} m\right)\le(C_1+1)^2.$$
Let $q=\frac{C_2}{2C_3}$, then
$$\left(\dashint_{B(x_0,R)}(u+\varepsilon)^q\mathrm{d} m\right)\left(\dashint_{B(x_0,R)}(u+\varepsilon)^{-q}\mathrm{d} m\right)\le (C_1+1)^2.$$
Letting $\varepsilon\downarrow0$, we have
$$\left(\dashint_{B(x_0,R)}u^q\mathrm{d} m\right)\left(\dashint_{B(x_0,R)}u^{-q}\mathrm{d} m\right)\le(C_1+1)^2.$$
By Corollary \ref{cor_infu}, there exists $C_4>0$ such that
\begin{align*}
&\einf_{B(x_0,\frac{1}{2}R)}u^q\ge \frac{1}{C_4^q}\left(\dashint_{B(x_0,R)}u^{-q}\mathrm{d} m\right)^{-1}\ge \frac{1}{CV(x_0,R)}\int_{B(x_0,R)}u^q\mathrm{d} m,
\end{align*}
where $C=(C_1+1)^2C_4^q$.
\end{proof}

We have the following sharp weak Harnack inequality; see also \cite[Theorem 3.59]{HKM06}, \cite[Theorem 4.3]{KM03}, and \cite[Theorem 8.34]{BB11}.

\begin{proposition}[Sharp weak Harnack inequality]\label{prop_wEHI_sharp}
Assume \ref{eq_VD}, \ref{eq_PI}, \ref{eq_CS}. Let $\kappa>1$ be the constant in Proposition \ref{prop_Sobolev}. Then there exists $A>1$ such that for any $q\in(0,\kappa(p-1))$, there exists $C>0$ such that the following holds: for any ball $B(x_0,R)$ and any $u\in \mathcal{F}$ which is non-negative superharmonic in $B(x_0,AR)$, we have
\begin{equation}\label{eq_wEHI_sharp}
\einf_{B(x_0,R)}u\ge \frac{1}{C}\left(\dashint_{B(x_0,R)}u^q \mathrm{d}m\right)^{1/q}.
\end{equation}
\end{proposition}

We shall use the following Caccioppoli inequality for superharmonic functions. Its proof follows standard arguments; see, for example, \cite[Lemma 3.57]{HKM06}, \cite[Lemma 3.1]{KM03}, and \cite[Proposition 8.8]{BB11}. We therefore omit the details.

\begin{lemma}[Caccioppoli inequality for superharmonic functions]\label{lem_Cacci_super}
Assume \ref{eq_VD}, \ref{eq_PI}. Then for any $\varepsilon>0$ and any bounded open set $U$, for any $u\in \mathcal{F}$ which is non-negative superharmonic in $U$, for any non-negative bounded $\phi\in \mathcal{F}(U)$, we have
$$\int_{U}\phi^pu^{-1-\varepsilon} \mathrm{d}\Gamma(u)\le \left(\frac{p}{\varepsilon}\right)^p\int_{U}u^{p-1-\varepsilon}\mathrm{d}\Gamma(\phi).$$
\end{lemma}

\begin{proof}[Proof of Proposition \ref{prop_wEHI_sharp}]
For notational convenience, we may assume that all functions in $\mathcal{F}$ are quasi-continuous. Let $q_1,A_1,C_1$ be the constants in Proposition \ref{prop_wEHI}, and $A=2A_1$, then for any ball $B=B(x_0,R)$ and any $u\in \mathcal{F}$ which is non-negative superharmonic in $B(x_0,AR)$, we have
\begin{equation*}
\einf_{B}u\ge \frac{1}{C_1^{1/q_1}}\left(\dashint_{2B}u^{q_1} \mathrm{d}m\right)^{1/q_1}\ge\frac{1}{(C_1C_{VD})^{1/{q_1}}}\left(\dashint_{B}u^{q_1} \mathrm{d}m\right)^{1/q_1}.
\end{equation*}
By H\"older's inequality, we have (\ref{eq_wEHI_sharp}) for any $q\in(0,q_1]$. If $q_1\ge\kappa(p-1)$, then the desired result follows. Hence we may assume that $q_1<\kappa(p-1)$, and it remains to prove (\ref{eq_wEHI_sharp}) for any $q\in(q_1,\kappa(p-1))$.

Firstly, we prove the result for bounded $u$. Let $k\ge1$ be the integer satisfying that $\frac{q}{\kappa^k}\le q_1 <\frac{q}{\kappa^{k-1}}$. Let $B_j=(1+\frac{j}{k})B$ for any $j=0,\ldots,k$, then $B=B_0$ and $2B=B_k$. For any $j=0,\ldots,k-1$, let $\phi_j\in \mathcal{F}(B_{j+1})$ be a cutoff function for $B_j\subseteq B_{j+1}$ chosen later, then $\phi_ju^{\frac{q}{p\kappa^{j+1}}}\in \mathcal{F}(B_{j+1})$, by Proposition \ref{prop_Sobolev}, there exists $C_2>0$ such that
\begin{equation}\label{eq_wEHI_sharpA}
\left(\int_{B_j}u^{\frac{q}{\kappa^j}}\mathrm{d}m\right)^{\frac{1}{\kappa}}\le\left(\int_{B_{j+1}}\left(\phi_ju^{\frac{q}{p\kappa^{j+1}}}\right)^{p\kappa}\mathrm{d}m\right)^{\frac{1}{\kappa}}\le C_2 \frac{\Psi((1+\frac{j+1}{k})R)}{m(B_{j+1})^{\frac{\kappa-1}{\kappa}}}\mathcal{E}(\phi_ju^{\frac{q}{p\kappa^{j+1}}}),
\end{equation}
where by \cite[Lemma 6.1]{Yan25c}, we have
\begin{equation}\label{eq_wEHI_sharpB}
\mathcal{E}(\phi_ju^{\frac{q}{p\kappa^{j+1}}})\le 2^{p-1}\left(\int_{B_{j+1}}\phi_j^p \mathrm{d}\Gamma(u^{\frac{q}{p\kappa^{j+1}}})+\int_{B_{j+1}}u^{\frac{q}{\kappa^{j+1}}}\mathrm{d}\Gamma(\phi_j)\right).
\end{equation}
Let $\delta>0$ be chosen later, then by the self-improvement property of \ref{eq_CS} (see \cite[Proposition 3.1]{Yan25c}), there exist $C_\delta>0$ and a cutoff function $\phi_j\in \mathcal{F}(B_{j+1})$ for $B_j\subseteq B_{j+1}$ such that
\begin{align*}
\int_{B_{j+1}}u^{\frac{q}{\kappa^{j+1}}}\mathrm{d}\Gamma(\phi_j)\le \delta\int_{B_{j+1}}\phi_j^p \mathrm{d}\Gamma(u^{\frac{q}{p\kappa^{j+1}}})+\frac{C_\delta}{\Psi(\frac{1}{k}R)}\int_{B_{j+1}}u^{\frac{q}{\kappa^{j+1}}} \mathrm{d}m,
\end{align*}
where
\begin{align}
&\int_{B_{j+1}}\phi_j^p \mathrm{d}\Gamma(u^{\frac{q}{p\kappa^{j+1}}})\le \left(\frac{q}{p\kappa}\right)^p\int_{B_{j+1}}\phi_j^pu^{\frac{q}{\kappa^{j+1}}-p}\mathrm{d}\Gamma(u)\nonumber\\
&\overset{(*)}{\scalebox{2}[1]{$\le$}}\left(\frac{q}{p\kappa}\right)^p \left(\frac{p}{p-1-\frac{q}{\kappa}}\right)^p\int_{B_{j+1}}u^{\frac{q}{\kappa^{j+1}}}\mathrm{d}\Gamma(\phi_j)=C_{p,q,\kappa}\int_{B_{j+1}}u^{\frac{q}{\kappa^{j+1}}}\mathrm{d}\Gamma(\phi_j),\label{eq_wEHI_sharpC}
\end{align}
where $(*)$ follows from the fact that $\frac{q}{\kappa^{j+1}}-p\le \frac{q}{\kappa}-p<-1$ and Lemma \ref{lem_Cacci_super}, and $C_{p,q,\kappa}>0$ depends only on $p,q,\kappa$. Let $\delta=\frac{1}{2C_{p,q,\kappa}}$, then
$$\int_{B_{j+1}}u^{\frac{q}{\kappa^{j+1}}}\mathrm{d}\Gamma(\phi_j)\le \frac{2C_\delta}{\Psi(\frac{1}{k}R)}\int_{B_{j+1}}u^{\frac{q}{\kappa^{j+1}}} \mathrm{d}m.$$
Combining this with (\ref{eq_wEHI_sharpA})--(\ref{eq_wEHI_sharpC}), we obtain that there exists $C_3>0$ such that
$$\dashint_{B_j}u^{\frac{q}{\kappa^j}}\mathrm{d}m\le C_3 \left(\dashint_{B_{j+1}}u^{\frac{q}{\kappa^{j+1}}}\mathrm{d}m\right)^\kappa.$$
Therefore
\begin{align*}
&\dashint_{B}u^q \mathrm{d}m=\dashint_{B_0}u^{\frac{q}{\kappa^0}} \mathrm{d}m\le C_3\left(\dashint_{B_{1}}u^{\frac{q}{\kappa^{1}}}\mathrm{d}m\right)^\kappa\le\ldots\le C_3^{1+\ldots+\kappa^{k-1}}\left(\dashint_{B_{k}}u^{\frac{q}{\kappa^{k}}}\mathrm{d}m\right)^{\kappa^k}\\
&\overset{(\blacktriangle)}{\scalebox{2}[1]{$\le$}}C_3^{1+\ldots+\kappa^{k-1}} \left(\dashint_{2B}u^{q_1}\mathrm{d}m\right)^{\frac{q}{q_1}}\overset{(\diamond)}{\scalebox{2}[1]{$\le$}} C_1^{\frac{q}{q_1}}C_3^{1+\ldots+\kappa^{k-1}}\einf_{B}u^q=C^q\einf_{B}u^q,
\end{align*}
where $(\blacktriangle)$ follows from the fact $\frac{q}{\kappa^k}\le q_1$ and H\"older's inequality, $(\diamond)$ follows from Proposition \ref{prop_wEHI}, and $C^q=C_1^{\frac{q}{q_1}}C_3^{1+\ldots+\kappa^{k-1}}$.

Secondly, we prove the result for general $u$. For any $n\ge1$, by Lemma \ref{lem_pasting}, we have $u\wedge n\in \mathcal{F}$ is non-negative bounded superharmonic in $B(x_0,AR)$. By above, we have
$$\einf_{B(x_0,R)}u\ge\einf_{B(x_0,R)}(u\wedge n)\ge\frac{1}{C}\left(\dashint_{B(x_0,R)}(u\wedge n)^q \mathrm{d}m\right)^{\frac{1}{q}}.$$
Letting $n\to+\infty$, we have the desired result.
\end{proof}

For the final step, combining the $L^q$-mean value inequality for subharmonic functions from Corollary \ref{cor_MVq} with the weak Harnack inequality for superharmonic functions from Proposition \ref{prop_wEHI}, we have the elliptic Harnack inequality for harmonic functions as follows.

\begin{proof}[Proof of \ref{eq_EHI} in Theorem \ref{thm_main_conv}]
Let $q,C_1,A_1$  be the positive constants from Proposition \ref{prop_wEHI}. Let $A_H=4A_1$. For such $q$, let $C_2$ be the positive constant from Corollary \ref{cor_MVq}. For any ball $B(x_0,R)$, for any $u\in\mathcal{F}$ which is non-negative harmonic in $B(x_0,A_HR)$, by replacing $u$ with $u\vee 0$, we may assume that $u$ is non-negative in $X$.

Firstly, we prove the result for bounded $u$. Indeed, it follows directly from Proposition \ref{prop_wEHI} and Corollary \ref{cor_MVq} that
\begin{align*}
&\einf_{B(x_0,R)}u^q\ge\einf_{B(x_0,2R)}u^q\ge\frac{1}{C_1V(x_0,4R)}\int_{B(x_0,4R)}u^q\mathrm{d} m\ge \frac{1}{C_1C_2}\esup_{B(x_0,R)}u^q,
\end{align*}
that is, $\esup_{B(x_0,R)}u\le C_H\einf_{B(x_0,R)}u$, where $C_H=(C_1C_2)^{1/q}$.

Secondly, we prove the result for general $u$. For any $n\ge1$, let $u_n=H^{B(x_0,A_HR)}(u\wedge n)$, then by Proposition \ref{prop_harm_exist}, we have $u_n$ is non-negative bounded harmonic in $B(x_0,A_HR)$, $\widetilde{u}_n=\widetilde{u}\wedge n\le \widetilde{u}$ q.e. on $X\backslash B(x_0,A_HR)$, and $\mathcal{E}(u_n)\le \mathcal{E}(u\wedge n)\le \mathcal{E}(u)$. By above, we have
\begin{equation}\label{eq_EHI1}
\esup_{B(x_0,R)}u_n\le C_H\einf_{B(x_0,R)}u_n.
\end{equation}
Since $u_n,u$ are both harmonic in $B(x_0,A_HR)$, by Proposition \ref{prop_comparison}, we have $u_n\le u$ in $B(x_0,A_HR)$, which gives $u_n\le u$ in $X$, $\lVert u_n\rVert_{L^p(X;m)}\le \lVert u\rVert_{L^p(X;m)}$. Similarly, since $\widetilde{u}_n=\widetilde{u}\wedge n\le\widetilde{u}\wedge (n+1)=\widetilde{u}_{n+1}$ q.e. on $X\backslash B(x_0,A_HR)$, we have $u_n\le u_{n+1}$ in $X$. Hence there exists non-negative $v\in L^p(X;m)$ with $v\le u$ in $X$ such that $u_n\uparrow v$ in $X$. Letting $n\to+\infty$ in (\ref{eq_EHI1}), we have
$$\esup_{B(x_0,R)}v\le C_H\einf_{B(x_0,R)}v.$$

We only need to show that $v=u$. Indeed, since $\{u_n\}_n\subseteq \mathcal{F}$ is $\mathcal{E}_1$-bounded, by the Banach-Alaoglu theorem (see \cite[Theorem 3 in Chapter12]{Lax02}), $\{u_n\}_n$ is $\mathcal{E}_1$-weakly convergent to some function $w\in \mathcal{F}$. By Mazur's lemma (see \cite[Theorem 2 in Section V.1]{Yos95}), for any $n\ge1$, there exist $I_n\ge n$, $\lambda^{(n)}_k\ge0$ for $k=n,\ldots,I_n$ with $\sum_{k=n}^{I_n}\lambda^{(n)}_k=1$ such that $\{\sum_{k=n}^{I_n}\lambda^{(n)}_ku_k\}_n$ is $\mathcal{E}_1$-convergent to $w$, hence $L^p(X;m)$-convergent to $w$. Since $\{\sum_{k=n}^{I_n}\lambda^{(n)}_ku_k\}_n$ is also $L^p(X;m)$-convergent to $v$, we have $v=w\in \mathcal{F}$. Since $\widetilde{u}\wedge n=\widetilde{u}_n\le \widetilde{v}\le \widetilde{u}$ q.e. on $X\backslash B(x_0,A_HR)$, letting $n\to+\infty$, we have $\widetilde{v}=\widetilde{u}$ q.e. on $X\backslash B(x_0,A_HR)$. For any $\phi\in \mathcal{F}\cap C_c(B(x_0,A_HR))$, since $\{u_n\}_n$ is $\mathcal{E}_1$-weakly convergent to $v$, we have
$$\lim_{n\to+\infty}\left(\mathcal{E}(u_n;\phi)+\int_Xu_n^{p-1}\phi \mathrm{d}m\right)=\mathcal{E}(v;\phi)+\int_Xv^{p-1}\phi \mathrm{d}m.$$
Since $u_n$ is harmonic in $B(x_0,A_HR)$, we have $\mathcal{E}(u_n;\phi)=0$. By the dominated convergence theorem, we have $\lim_{n\to+\infty}\int_Xu_n^{p-1}\phi \mathrm{d}m=\int_Xv^{p-1}\phi \mathrm{d}m$. Hence $\mathcal{E}(v;\phi)=0$ for any $\phi\in \mathcal{F}\cap C_c(B(x_0,A_HR))$, that is, $v$ is harmonic in $B(x_0,A_HR)$.

In summary, $u,v\in \mathcal{F}$ are both harmonic in $B(x_0,A_HR)$ and satisfy $\widetilde{u}=\widetilde{v}$ q.e. on $X\backslash B(x_0,A_HR)$. By Proposition \ref{prop_comparison}, we therefore have $u=v$.
\end{proof}

We give the proof of \ref{eq_EHIann} as follows.

\begin{proof}[Proof of \ref{eq_EHIann} in Theorem \ref{thm_main_conv}]
Let $A>1$ be the constant in Proposition \ref{prop_wEHI_sharp}. Let $B=B(x_0,R)$ be a ball and $u\in \mathcal{F}$ which is non-negative superharmonic in $4AB$ and harmonic in $4AB\backslash \overline{B}$. By \ref{eq_EHI}, we have $u\in C(4AB\backslash \overline{B})$. For any $x\in 3B\backslash 2B$, we have $u\in \mathcal{F}$ is non-negative bounded (sub)harmonic in $B(x,\frac{1}{2}R)$, then by Corollary \ref{cor_MVq}, for any $q>0$, there exists $C_1>0$ such that
$$\esup_{B(x,\frac{1}{8}R)}u^q\le C_1\dashint_{B(x,\frac{1}{2}R)}u^q \mathrm{d}m\le C_1C_{VD}^5\dashint_{4B}u^q \mathrm{d}m.$$
By Proposition \ref{prop_wEHI_sharp}, for any $q\in(0,\kappa(p-1))$, there exists $C_2>0$ such that
$$\dashint_{4B}u^q \mathrm{d}m\le C_2\einf_{4B}u^q.$$
By Proposition \ref{prop_comparison}, we have
$$\einf_{4B}u=\einf_{4B\backslash 3B}u\le\einf_{3B\backslash 2B}u.$$
Hence
$$\esup_{B(x,\frac{1}{8}R)}u^q\le C_1C_2C_{VD}^5\einf_{3B\backslash 2B}u.$$
Taking the supremum with respect to $x\in 3B\backslash 2B$, we have
$$\esup_{3B\backslash 2B}u\le C\einf_{3B\backslash 2B}u,$$
where $C=(C_1C_2C_{VD}^5)^{1/q}$.
\end{proof}

\section{Proofs of Theorems \ref{thm_main_RSVR} and \ref{thm_HK}}\label{sec_RSVR}

First, we establish the following result, which yields the \hyperlink{eq_LLC}{\text{LLC}} condition and hence allows us to apply Theorem \ref{thm_main}.

\begin{proposition}\label{prop_LLC}
Assume \ref{eq_FVR}, \ref{eq_PI}, \ref{eq_ucap}. Then \hyperlink{eq_LLC}{\text{LLC}} holds.
\end{proposition}

We follow the argument used in the proof of \cite[Theorem 3.3]{Kor07}; see also \cite[THEOREM 3.13]{HK98}, \cite[Theorem 9.4.1]{HKST15}, \cite[Proposition 3.4]{Mur19} for related results. The outline of the proof is as follows. For any $\varepsilon>0$, pick up an $\varepsilon$-net $V^{(\varepsilon)}$ of $X$, we construct an infinite graph $(V^{(\varepsilon)},E^{(\varepsilon)})$, where $V^{(\varepsilon)}$ is the set of vertices and $E^{(\varepsilon)}$ is the set of edges, which consists of all pairs of vertices whose distance is ``approximately" $\varepsilon$. We introduce a ``naturally" defined metric $d^{(\varepsilon)}$, measure $m^{(\varepsilon)}$, and $p$-energy $(\mathbf{E}^{(\varepsilon)},\mathbf{F}^{(\varepsilon)})$ on $V^{(\varepsilon)}$ such that at large scales the volume regular condition, the Poincar\'e inequality and the capacity upper bound still hold. Using the argument in \cite{Kor07}, we prove the \hyperlink{eq_LLC}{\text{LLC}} condition at large scales on $V^{(\varepsilon)}$, which in turn implies the \hyperlink{eq_LLC}{\text{LLC}} condition on $X$.

We now start the formal proof. An important remark is that the constants appearing in the following results will be \emph{independent} of $\varepsilon$. For any $\varepsilon>0$, let $V^{(\varepsilon)}$ be an $\varepsilon$-net and
$$E^{(\varepsilon)}=\left\{(z_1,z_2):z_1,z_2\in V^{(\varepsilon)},z_1\ne z_2, B(z_1,\frac{5}{4}\varepsilon)\cap B(z_2,\frac{5}{4}\varepsilon)\ne\emptyset\right\},$$
then for any $(z_1,z_2)\in E^{(\varepsilon)}$, we have $\varepsilon<d(z_1,z_2)\le \frac{5}{2}\varepsilon$. Assume \ref{eq_VD}, then there exists some positive integer $N$ depending only on $C_{VD}$ such that
\begin{equation}\label{eq_bdddeg}
\sup_{\overline{z}\in V^{(\varepsilon)}}\# \left\{z\in V^{(\varepsilon)}:(z,\overline{z})\in E^{(\varepsilon)}\right\}\le N.
\end{equation}
For any distinct $x,y\in V^{{(\varepsilon)}}$, let $d^{(\varepsilon)}(x,x)=0$ and
$$d^{(\varepsilon)}(x,y)=\inf \left\{\sum_{k=1}^Kd(z_{k-1},z_k):z_0=x,z_K=y,(z_{k-1},z_k)\in E^{(\varepsilon)}\text{ for any }k=1,\ldots,K\right\}.$$
Obviously, $d^{(\varepsilon)}$ defines a metric on $V^{{(\varepsilon)}}$ and $d(z_1,z_2)\le d^{(\varepsilon)}(z_1,z_2)$ for any $z_1,z_2\in V^{(\varepsilon)}$. Moreover, by \ref{eq_CC}, there exists $L\ge1$ depending only on $C_{cc}$ such that $d^{(\varepsilon)}(z_1,z_2)\le Ld(z_1,z_2)$ for any $z_1,z_2\in V^{(\varepsilon)}$, hence
\begin{equation}\label{eq_eps_ball}
B(z,\frac{1}{L}r)\cap V^{(\varepsilon)}\subseteq B^{(\varepsilon)}(z,r)\subseteq B(z,r)\cap V^{(\varepsilon)}\text{ for any }z\in V^{(\varepsilon)},r>0,
\end{equation}
where $B^{(\varepsilon)}(z,r)=\{y\in V^{(\varepsilon)}:d^{(\varepsilon)}(y,z)<r\}$. See also \cite[Subsection 6.1]{GT12} for another similar metric based on $\varepsilon$-chains. For any subset $A$ of $V^{(\varepsilon)}$, let $m^{(\varepsilon)}(A)=\Phi(\varepsilon)(\# A)$. It is obvious that
\begin{equation*}
m^{(\varepsilon)}(B^{(\varepsilon)}(z,r))\asymp \Phi(r)\text{ for any }z\in V^{(\varepsilon)}, r>\varepsilon.
\end{equation*}

For any subset $A\subseteq X$ or $A\subseteq V^{(\varepsilon)}$, for any $\mathbf{u}\in l(A\cap V^{(\varepsilon)})$, let
$$\mathbf{E}_A^{(\varepsilon)}(\mathbf{u})=\mathbf{E}_{A\cap V^{(\varepsilon)}}^{(\varepsilon)}(\mathbf{u})=\frac{\Phi(\varepsilon)}{\Psi(\varepsilon)}\sum_{(z_1,z_2)\in E^{(\varepsilon)},z_1,z_2\in A\cap V^{(\varepsilon)}}|\mathbf{u}(z_1)-\mathbf{u}(z_2)|^p.$$
Let $\mathbf{E}^{(\varepsilon)}=\mathbf{E}^{(\varepsilon)}_{X}=\mathbf{E}^{(\varepsilon)}_{V^{(\varepsilon)}}$, and
$$\mathbf{F}^{(\varepsilon)}=\left\{\mathbf{u}\in l(V^{(\varepsilon)}):\sum_{(z_1,z_2)\in E^{(\varepsilon)}}|\mathbf{u}(z_1)-\mathbf{u}(z_2)|^p<+\infty\right\}.$$

For any subsets $A_0,A_1$ of $V^{(\varepsilon)}$, we define the capacity between $A_0,A_1$ as
$$\mathrm{cap}^{(\varepsilon)}(A_0,A_1)=\inf \left\{\mathbf{E}^{(\varepsilon)}(\mathbf{u}):\mathbf{u}\in l(V^{(\varepsilon)}),\mathbf{u}=0\text{ on }A_0,\mathbf{u}=1\text{ on }A_1\right\}.$$

To carry out the discretization and approximation procedures, we make use of the following partition of unity with controlled energy.

\begin{lemma}[{\cite[Lemma 3.2]{Yan25a}}]\label{lem_partition}
Assume \ref{eq_VD}, \ref{eq_ucap}. Then we have the following controlled cutoff condition. There exists $C_{cut}>0$ depending only on $p, C_{\Psi}, C_{VD}, C_{cap}, A_{cap}$ such that for any $\varepsilon>0$, for any $\varepsilon$-net $V$, there exists a family of functions $\{\psi_z\in\mathcal{F}:z\in V\}$ satisfying the following conditions.
\begin{enumerate}[label=(CO\arabic*)]
\item\label{item_COspt} For any $z\in V$, $0\le\psi_z\le 1$ in $X$, $\psi_z=1$ in $B(z,\varepsilon/4)$, and $\psi_z=0$ on $X\backslash B(z,5\varepsilon/4)$.
\item\label{item_COunit} $\sum_{z\in V}\psi_z=1$.
\item\label{item_COenergy} For any $z\in V$, $\mathcal{E}(\psi_z)\le C_{cut}\frac{V(z,\varepsilon)}{\Psi(\varepsilon)}$.
\end{enumerate}
\end{lemma}

Let $\{\psi_z^{(\varepsilon)}\in\mathcal{F}:z\in V^{(\varepsilon)}\}$ be the partition of unity given by Lemma \ref{lem_partition}. For any $\mathbf{u}\in l(V^{(\varepsilon)})$, let $u=\sum_{z\in V^{(\varepsilon)}}\mathbf{u}(z)\psi_z^{(\varepsilon)}$, then $u\in\mathcal{F}_{loc}$, hence $\Gamma(u)$ is well-defined. We have the following estimate for $\Gamma(u)(\cdot)$ in terms of $\mathbf{E}^{(\varepsilon)}_{\bullet}(\mathbf{u})$.

\begin{lemma}\label{lem_Gammabfu}
Assume \ref{eq_VPhi}, \ref{eq_ucap}. Then there exists $C>0$ such that for any $\mathbf{u}\in l(V^{(\varepsilon)})$, let $u=\sum_{z\in V^{(\varepsilon)}}\mathbf{u}(z)\psi_z^{(\varepsilon)}$, for any ball $B(x_0,R)$ in $X$, we have
$$\Gamma(u)(B(x_0,R))\le C\mathbf{E}^{(\varepsilon)}_{B(x_0,R+\frac{5}{2}\varepsilon)}(\mathbf{u}).$$
\end{lemma}

\begin{proof}
For any $\overline{z}\in V^{(\varepsilon)}$, let
$$N_{\overline{z}}=\{z\in V^{(\varepsilon)}:B(z,\frac{5}{4}\varepsilon)\cap B(\overline{z},\frac{5}{4}\varepsilon)\ne\emptyset\}=\{\overline{z}\}\bigcup\{z\in V^{(\varepsilon)}:(z,\overline{z})\in E^{(\varepsilon)}\},$$
then by (\ref{eq_bdddeg}), we have $\# N_{\overline{z}}\le N+1$, and by Lemma \ref{lem_partition}, we have
\begin{align*}
&\Gamma(u)\left(B(\overline{z},\frac{5}{4}\varepsilon)\right)=\Gamma \left(\sum_{z\in N_{\overline{z}}}\left(\mathbf{u}(z)-\mathbf{u}(\overline{z})\right)\psi_z^{(\varepsilon)}\right)\left(B(\overline{z},\frac{5}{4}\varepsilon)\right)\\
&\le C_1\sum_{z\in N_{\overline{z}}}\Gamma \left(\left(\mathbf{u}(z)-\mathbf{u}(\overline{z})\right)\psi_z^{(\varepsilon)}\right)\left(B(\overline{z},\frac{5}{4}\varepsilon)\right)\le C_1C_{VR}C_{cut}\frac{\Phi(\varepsilon)}{\Psi(\varepsilon)}\sum_{{z:(z,\overline{z})\in E^{(\varepsilon)}}}|\mathbf{u}(z)-\mathbf{u}(\overline{z})|^p,
\end{align*}
where $C_1>0$ depends only on $p, N$. Hence
\begin{align*}
&\Gamma(u)(B(x_0,R))\le\sum_{{\overline{z}\in B(x_0,R+\frac{5}{4}\varepsilon)\cap V^{(\varepsilon)}}}\Gamma(u)(B(\overline{z},\frac{5}{4}\varepsilon))\\
&\le NC_1C_{VR}C_{cut}\frac{\Phi(\varepsilon)}{\Psi(\varepsilon)}\sum_{\substack{(z_1,z_2)\in E^{(\varepsilon)}\\z_1,z_2\in B(x_0,R+\frac{5}{2}\varepsilon)\cap V^{(\varepsilon)}}}|\mathbf{u}(z_1)-\mathbf{u}(z_2)|^p=C\mathbf{E}^{(\varepsilon)}_{B(x_0,R+\frac{5}{2}\varepsilon)}(\mathbf{u}),
\end{align*}
which gives the desired result with $C=NC_1C_{VR}C_{cut}$.
\end{proof}

We have the Poincar\'e inequality on $V^{(\varepsilon)}$ as follows.

\begin{proposition}\label{prop_PI_Veps}
Assume \ref{eq_VPhi}, \ref{eq_PI}, \ref{eq_ucap}. Then there exist $C>0$, $A>1$ such that for any ball $B^{(\varepsilon)}(z,R)$ in $V^{(\varepsilon)}$ with $R>\varepsilon$, for any $\mathbf{u}\in l(V^{(\varepsilon)})$, we have
\begin{align*}
&\int_{B^{(\varepsilon)}(z,R)}|\mathbf{u}-\frac{1}{m^{(\varepsilon)}(B^{(\varepsilon)}(z,R))}\int_{B^{(\varepsilon)}(z,R)}\mathbf{u}\mathrm{d} m^{(\varepsilon)}|^p\mathrm{d} m^{(\varepsilon)}\\
&\le C\Psi(R)\mathbf{E}^{(\varepsilon)}_{B^{(\varepsilon)}(z,AR)}(\mathbf{u}).\label{eq_PI_Veps}\tag*{$\mathbf{PI}^{(\varepsilon)}(\Psi)$}
\end{align*}

\end{proposition}

\begin{proof}
Let $u=\sum_{z\in V^{(\varepsilon)}}\mathbf{u}(z)\psi_z^{(\varepsilon)}$, then for any $c\in \mathbb{R}$, we have
\begin{align*}
\int_{B^{(\varepsilon)}(z,R)}|\mathbf{u}-c|^p\mathrm{d} m^{(\varepsilon)}\lesssim \sum_{x\in B^{(\varepsilon)}(z,R)}\int_{B(x,\frac{1}{4}\varepsilon)}|u-c|^p\mathrm{d} m\le \int_{B(z,2R)}|u-c|^p\mathrm{d} m.
\end{align*}
Let $c=u_{B(z,2R)}$, then by \ref{eq_PI} and Lemma \ref{lem_Gammabfu}, we have
\begin{align*}
&\int_{B^{(\varepsilon)}(z,R)}|\mathbf{u}-\frac{1}{m^{(\varepsilon)}(B^{(\varepsilon)}(z,R))}\int_{B^{(\varepsilon)}(z,R)}\mathbf{u}\mathrm{d} m^{(\varepsilon)}|^p\mathrm{d} m^{(\varepsilon)}\le2^p\int_{B^{(\varepsilon)}(z,R)}|\mathbf{u}-c|^p\mathrm{d} m^{(\varepsilon)}\\
&\lesssim \int_{B(z,2R)}|u-u_{B(z,2R)}|^p\mathrm{d} m\lesssim \Psi(R)\int_{B(z,2A_{PI}R)}\mathrm{d}\Gamma(u)\lesssim \Psi(R)\mathbf{E}^{(\varepsilon)}_{B^{(\varepsilon)}(z,8A_{PI}LR)}(\mathbf{u}).
\end{align*}
\end{proof}

We have the capacity upper bound on $V^{(\varepsilon)}$ as follows.

\begin{proposition}\label{prop_ucap_Veps}
Assume \ref{eq_VPhi}, \ref{eq_PI}, \ref{eq_ucap}. Then there exist $C>0$, $A>1$ such that for any ball $B^{(\varepsilon)}(x_0,R)$ in $V^{(\varepsilon)}$ with $R>\varepsilon$, we have
\begin{equation*}\label{eq_ucap_Veps}\tag*{$\mathbf{cap}^{(\varepsilon)}(\Psi)$}
\mathrm{cap}^{(\varepsilon)}(B^{(\varepsilon)}(x_0,R),V^{(\varepsilon)}\backslash B^{(\varepsilon)}(x_0,AR))\le C \frac{\Phi(R)}{\Psi(R)}.
\end{equation*}
\end{proposition}

\begin{proof}
By \ref{eq_VPhi}, \ref{eq_ucap}, there exists a cutoff function $\phi\in\mathcal{F}$ for\\
\noindent $B(x_0,2R)\subseteq B(x_0,2A_{cap}R)$ such that $\mathcal{E}(\phi)\lesssim \frac{\Phi(R)}{\Psi(R)}$. For any $z\in V^{(\varepsilon)}$, let $\psi(z)=\phi_{B(z,\varepsilon)}$, then $\psi=1$ in $B(x_0,R)$ and $\psi=0$ on $V^{(\varepsilon)}\backslash B(x_0,4A_{cap}R)$. By (\ref{eq_eps_ball}), we have $\psi=1$ in $B^{(\varepsilon)}(x_0,R)$ and $\psi=0$ on $V^{(\varepsilon)}\backslash B^{(\varepsilon)}(x_0,4A_{cap}LR)$. Moreover
$$\mathbf{E}^{(\varepsilon)}(\psi)=\frac{\Phi(\varepsilon)}{\Psi(\varepsilon)}\sum_{(z_1,z_2)\in E^{(\varepsilon)}}|\phi_{B(z_1,\varepsilon)}-\phi_{B(z_2,\varepsilon)}|^p,$$
where by \ref{eq_PI}, we have
\begin{align*}
&|\phi_{B(z_1,\varepsilon)}-\phi_{B(z_2,\varepsilon)}|^p\lesssim\frac{\Psi(\varepsilon)}{\Phi(\varepsilon)}\int_{B(z_1,4A_{PI}\varepsilon)}\mathrm{d}\Gamma(\phi),
\end{align*}
then by (\ref{eq_bdddeg}), we have
\begin{align*}
&\mathbf{E}^{(\varepsilon)}(\psi)\lesssim\sum_{(z_1,z_2)\in E^{(\varepsilon)}}\int_{B(z_1,4A_{PI}\varepsilon)}\mathrm{d}\Gamma(\phi)\lesssim\int_X \left(\sum_{z\in V^{(\varepsilon)}}1_{B(z,4A_{PI}\varepsilon)}\right)\mathrm{d}\Gamma(\phi).
\end{align*}
By \ref{eq_VD}, there exists some positive integer $M$ depending only on $C_{VD},A_{PI}$ such that
$$\sum_{z\in V^{(\varepsilon)}}1_{B(z,4A_{PI}\varepsilon)}\le M1_{\cup_{z\in V^{(\varepsilon)}}B(z,4A_{PI}\varepsilon)},$$
hence
$$\mathrm{cap}^{(\varepsilon)}(B^{(\varepsilon)}(x_0,R),V^{(\varepsilon)}\backslash B^{(\varepsilon)}(x_0,4A_{cap}LR))\le\mathbf{E}^{(\varepsilon)}(\psi)\lesssim \mathcal{E}(\phi)\lesssim \frac{\Phi(R)}{\Psi(R)}.$$
\end{proof}

The following version of the capacity lower bound will play a crucial role. The proof technique has appeared previously in the literature, notably in \cite[THEOREM 5.9]{HK98}, \cite[Lemma 3.2]{Kor07}, and \cite[PROPOSITION 3.2]{Mur23a}. A similar argument will also be applied in the proof of Lemma \ref{lem_lcap_Xeps}.

\begin{lemma}\label{lem_lcap_Veps}
Assume \ref{eq_VPhi}, \ref{eq_PI}, \ref{eq_ucap}. Then there exist $C>0$, $A>1$ such that for any $x_0\in V^{(\varepsilon)}$ and any $R>256\varepsilon$, for any balls $B^{(\varepsilon)}(y_0,\frac{1}{16}R)$, $B^{(\varepsilon)}(z_0,\frac{1}{16}R)$ contained in $B^{(\varepsilon)}(x_0,R)$, for any $\mathbf{u}\in l(V^{(\varepsilon)})$ with $\mathbf{u}\ge1$ in $B^{(\varepsilon)}(y_0,\frac{1}{16}R)$ and $\mathbf{u}\le 0$ in $B^{(\varepsilon)}(z_0,\frac{1}{16}R)$,
$$\mathbf{E}^{(\varepsilon)}_{B^{(\varepsilon)}(x_0,AR)}(\mathbf{u})\ge \frac{1}{C}\frac{\Phi(R)}{\Psi(R)}.$$
\end{lemma}

\begin{proof}
For notational convenience, we write $\mathbf{u}_{W}=\frac{1}{m^{(\varepsilon)}(W)}\int_{W}\mathbf{u}\mathrm{d} m^{(\varepsilon)}$ for any non-empty finite subset $W$ of $V^{(\varepsilon)}$. By assumption, we have $\mathbf{u}_{B^{(\varepsilon)}(y,2\varepsilon)}\ge1$ for any $y\in B^{(\varepsilon)}(y_0,\frac{1}{32}R)$, and $\mathbf{u}_{B^{(\varepsilon)}(z,2\varepsilon)}\le0$ for any $z\in B^{(\varepsilon)}(z_0,\frac{1}{32}R)$.

If there exist $y\in B^{(\varepsilon)}(y_0,\frac{1}{32}R)$ and $z\in B^{(\varepsilon)}(z_0,\frac{1}{32}R)$ such that $|\mathbf{u}_{B^{(\varepsilon)}(y,2\varepsilon)}-\mathbf{u}_{B^{(\varepsilon)}(y,R)}|\le \frac{1}{5}$ and $|\mathbf{u}_{B^{(\varepsilon)}(z,2\varepsilon)}-\mathbf{u}_{B^{(\varepsilon)}(z,R)}|\le \frac{1}{5}$, then $|\mathbf{u}_{B^{(\varepsilon)}(y,R)}-\mathbf{u}_{B^{(\varepsilon)}(z,R)}|\ge \frac{3}{5}$. By \ref{eq_PI_Veps}, we have
\begin{align*}
|\mathbf{u}_{B^{(\varepsilon)}(y,R)}-\mathbf{u}_{B^{(\varepsilon)}(z,R)}|^p\lesssim \frac{\Psi(R)}{\Phi(R)}\mathbf{E}^{(\varepsilon)}_{B^{(\varepsilon)}(x_0,2A_1R)}(\mathbf{u}),
\end{align*}
where $A_1$ is the positive constant appearing therein, hence
$$\mathbf{E}^{(\varepsilon)}_{B^{(\varepsilon)}(x_0,2A_1R)}(\mathbf{u})\gtrsim \frac{\Phi(R)}{\Psi(R)}.$$

Assume that no such $y,z$ as above exist, then without loss of generality, we may assume that $|\mathbf{u}_{B^{(\varepsilon)}(y,2\varepsilon)}-\mathbf{u}_{B^{(\varepsilon)}(y,R)}|>\frac{1}{5}$ for any $y\in B^{(\varepsilon)}(y_0,\frac{1}{32}R)$. Let $n$ be the integer satisfying that $2^n(2\varepsilon)\le R<2^{n+1}(2\varepsilon)$, then by \ref{eq_PI_Veps} again, we have
\begin{align*}
\frac{1}{5}&<|\mathbf{u}_{B^{(\varepsilon)}(y,2\varepsilon)}-\mathbf{u}_{B^{(\varepsilon)}(y,R)}|\\
&\le\sum_{k=0}^{n-1}|\mathbf{u}_{B^{(\varepsilon)}(y,2^k(2\varepsilon))}-\mathbf{u}_{B^{(\varepsilon)}(y,2^{k+1}(2\varepsilon))}|+|\mathbf{u}_{B^{(\varepsilon)}(y,2^{n}(2\varepsilon))}-\mathbf{u}_{B^{(\varepsilon)}(y,R)}|\\
&\lesssim\sum_{k=0}^n \left(\frac{\Psi(2^k\varepsilon)}{\Phi(2^k\varepsilon)}\mathbf{E}^{(\varepsilon)}_{B^{(\varepsilon)}(y,2^{k+2}A_1\varepsilon)}(\mathbf{u})\right)^{1/p},
\end{align*}
that is,
$$\sum_{k=0}^n \left(\frac{\Psi(2^k\varepsilon)}{\Phi(2^k\varepsilon)}\mathbf{E}^{(\varepsilon)}_{B^{(\varepsilon)}(y,2^{k+2}A_1\varepsilon)}(\mathbf{u})\right)^{1/p}\ge C_1,$$
where $C_1>0$ is some constant. Hence there exists $k_y=0,\ldots,n$ such that
$$\frac{\Psi(2^{k_y}\varepsilon)}{\Phi(2^{k_y}\varepsilon)}\mathbf{E}^{(\varepsilon)}_{B^{(\varepsilon)}(y,2^{k_y+2}A_1\varepsilon)}(\mathbf{u})\ge \left(\frac{C_1}{2}\right)^p2^{p(k_y-n)}.$$
Otherwise
$$C_1\le\sum_{k=0}^n \left(\frac{\Psi(2^k\varepsilon)}{\Phi(2^k\varepsilon)}\mathbf{E}^{(\varepsilon)}_{B^{(\varepsilon)}(y,2^{k+2}A_1\varepsilon)}(\mathbf{u})\right)^{1/p}<\sum_{k=0}^n \frac{C_1}{2} 2^{k-n}\le C_1,$$
contradiction. By the $5B$-covering lemma (see \cite[Theorem 1.2]{Hei01}), there exists a countable family of disjoint balls $\{B^{(\varepsilon)}(y_l,2^{k_l+2}A_1\varepsilon)\}_l$ such that
$$B^{(\varepsilon)}(y_0,\frac{1}{32}R)\subseteq \bigcup_{l}B^{(\varepsilon)}(y_l,5\cdot2^{k_l+2}A_1\varepsilon).$$
Hence
$$\mathbf{E}^{(\varepsilon)}_{B^{(\varepsilon)}(x_0,4A_1R)} (\mathbf{u})\ge\sum_{l}\mathbf{E}^{(\varepsilon)}_{B^{(\varepsilon)}(y_l,2^{k_l+2}A_1\varepsilon)}(\mathbf{u})\ge\sum_{l}\left(\frac{C_1}{2}\right)^p2^{p(k_l-n)}\frac{\Phi(2^{k_l}\varepsilon)}{\Psi(2^{k_l}\varepsilon)},$$
where by \ref{eq_VPhi}, we have
$$\Phi(2^{k_l}\varepsilon)\asymp{m^{(\varepsilon)}(B^{(\varepsilon)}(y_l,5\cdot2^{k_l+2}A_1\varepsilon))},$$
and by \cite[Proposition 2.1, Remark 2.2]{Yan25a}, under \ref{eq_VD}, \ref{eq_PI}, \ref{eq_ucap}, we have
$$\frac{\Psi(R)}{\Psi(2^{k_l}\varepsilon)}\gtrsim \left(\frac{R}{2^{k_l}\varepsilon}\right)^p\asymp 2^{p(n-k_l)},$$
hence
\begin{align*}
&\mathbf{E}^{(\varepsilon)}_{B^{(\varepsilon)}(x_0,4A_1R)} (\mathbf{u})\gtrsim\frac{1}{\Psi(R)}\sum_{l}{m^{(\varepsilon)}(B^{(\varepsilon)}(y_l,5\cdot2^{k_l+2}A_1\varepsilon))}\\
&\ge \frac{1}{\Psi(R)}m^{(\varepsilon)}(B^{(\varepsilon)}(y_0,\frac{1}{32}R))\asymp \frac{\Phi(R)}{\Psi(R)}.
\end{align*}
\end{proof}

We introduce the notion of combinatorial $p$-modulus following \cite[Section 2]{MS25}. A finite sequence of vertices $\theta=\{z_0,\ldots,z_n\}$ is called a path in $(V^{(\varepsilon)},E^{(\varepsilon)})$ if $(z_{i},z_{i+1})\in E^{(\varepsilon)}$ for any $i=0,\ldots,n-1$. For any non-negative function $\rho\in l(V^{(\varepsilon)})$, let $L_\rho(\theta)$ be the $\rho$-length of the path $\theta$ given by $L_\rho(\theta)=\sum_{i=0}^n\rho(z_i)$. For a family $\Theta$ of paths in $(V^{(\varepsilon)},E^{(\varepsilon)})$, we define the set $\mathrm{Adm}(\Theta)$ of admissible functions for $\Theta$ as
$$\mathrm{Adm}(\Theta)=\left\{\rho\in l(V^{(\varepsilon)}):\text{non-negative},L_\rho(\theta)\ge1\text{ for any }\theta\in\Theta\right\},$$
and the combinatorial $p$-modulus $\mathrm{mod}^{(\varepsilon)}(\Theta)$ of $\Theta$ as
$$\mathrm{mod}^{(\varepsilon)}(\Theta)=\inf \left\{\sum_{v\in V^{(\varepsilon)}}\rho(v)^p:\rho\in\mathrm{Adm}(\Theta)\right\}.$$
For any subsets $A_0, A_1, A_2$ of $V^{(\varepsilon)}$ with $A_0\cup A_1\subseteq A_2$, let
\begin{align*}
&\mathrm{Path}(A_0,A_1;A_2)=\left\{\{z_0,\ldots,z_n\}\middle|
\begin{array}{l}
z_0\in A_0,z_n\in A_1,z_i\in A_2\text{ for any }i=0,\ldots,n,\\
(z_i,z_{i+1})\in E^{(\varepsilon)}\text{ for any }i=0,\ldots,n-1
\end{array}
\right\}
\end{align*}
be the family of all paths contained in $A_2$ that connect $A_0$ and $A_1$, write 
$$\mathrm{mod}^{(\varepsilon)}(A_0,A_1;A_2)=\mathrm{mod}^{(\varepsilon)}\left(\mathrm{Path}(A_0,A_1;A_2)\right),$$
let
$$\mathrm{cap}^{(\varepsilon)}(A_0,A_1;A_2)=\inf\left\{\mathbf{E}^{(\varepsilon)}_{A_2}(\mathbf{u}):\mathbf{u}\in l(A_2),\mathbf{u}=0\text{ on }A_0,\mathbf{u}=1\text{ on }A_1\right\}.$$
Indeed, $\mathrm{cap}^{(\varepsilon)}(A_0,A_1;V^{(\varepsilon)})=\mathrm{cap}^{(\varepsilon)}(A_0,A_1)$ for any subsets $A_0, A_1$ of $V^{(\varepsilon)}$. By \cite[Lemma 2.12]{MS25} and (\ref{eq_bdddeg}), we have the modulus $\mathrm{mod}^{(\varepsilon)}$ and the capacity $\mathrm{cap}^{(\varepsilon)}$ are comparable in the following sense: there exists $C_{mod}>0$ depending only on $p,C_{VD}$ such that
$$\frac{1}{C_{mod}}\mathrm{mod}^{(\varepsilon)}(A_0,A_1;A_2)\le \frac{\Psi(\varepsilon)}{\Phi(\varepsilon)}\mathrm{cap}^{(\varepsilon)}(A_0,A_1;A_2)\le {C_{mod}}\mathrm{mod}^{(\varepsilon)}(A_0,A_1;A_2)$$
for any subsets $A_0, A_1, A_2$ of $V^{(\varepsilon)}$ with $A_0\cup A_1\subseteq A_2$.

We give the proof of the \hyperlink{eq_LLC}{\text{LLC}} condition as follows.

\begin{proof}[Proof of Proposition \ref{prop_LLC}]
By (\ref{eq_eps_ball}), we only need to show that there exists $A>1$ such that for any $R>0$ and any $\varepsilon\in(0,\frac{1}{256A}R)$, for any $\varepsilon$-net $V^{(\varepsilon)}$ and any $x_0\in V^{(\varepsilon)}$, for any balls $B^{(\varepsilon)}(y_0,\frac{1}{16}R)$, $B^{(\varepsilon)}(z_0,\frac{1}{16}R)$ contained in $B^{(\varepsilon)}(x_0,R)\backslash B^{(\varepsilon)}(x_0,\frac{1}{2}R)$, we have
$$\mathrm{mod}^{(\varepsilon)}\left(B^{(\varepsilon)}(y_0,\frac{1}{16}R),B^{(\varepsilon)}(z_0,\frac{1}{16}R);B^{(\varepsilon)}(x_0,AR)\backslash B^{(\varepsilon)}(x_0,\frac{1}{2A}R)\right)>0.$$

Let $C_1,A_1$ be the constants appearing in Lemma \ref{lem_lcap_Veps}, then
$$\mathrm{cap}^{(\varepsilon)}\left(B^{(\varepsilon)}(y_0,\frac{1}{16}R),B^{(\varepsilon)}(z_0,\frac{1}{16}R);B^{(\varepsilon)}(x_0,A_1R)\right)\ge \frac{1}{C_1}\frac{\Phi(R)}{\Psi(R)}.$$
Let $C_2, A_2$ be the constants appearing in Proposition \ref{prop_ucap_Veps}. Let $M\ge1$ be some integer chosen later. Let
\begin{align*}
\Gamma&=\mathrm{Path}\left(B^{(\varepsilon)}(y_0,\frac{1}{16}R),B^{(\varepsilon)}(z_0,\frac{1}{16}R);B^{(\varepsilon)}(x_0,A_1R)\right),\\
\Gamma_1&=\mathrm{Path}\left(B^{(\varepsilon)}(y_0,\frac{1}{16}R),B^{(\varepsilon)}(z_0,\frac{1}{16}R);B^{(\varepsilon)}(x_0,A_1R)\backslash B^{(\varepsilon)}(x_0,\frac{1}{2A_2^M}R)\right),\\
\Gamma_2&=\mathrm{Path}\left(B^{(\varepsilon)}(x_0,\frac{1}{2A_2^{M}}R),V^{(\varepsilon)}\backslash B^{(\varepsilon)}(x_0,\frac{1}{2}R);V^{(\varepsilon)}\right),
\end{align*}
then for any $\theta=\{z_0,\ldots,z_n\}\in\Gamma$, either $\theta\in\Gamma_1$, or there exist $k,l$ with $0\le k\le l\le n$ such that $\{z_k,\ldots,z_l\}\in\Gamma_2$. By \cite[Lemma 2.3 (iii) (iv)]{MS25}, we have
\begin{equation}\label{eq_LLC1}
\mathrm{mod}^{(\varepsilon)}(\Gamma)\le\mathrm{mod}^{(\varepsilon)}(\Gamma_1)+\mathrm{mod}^{(\varepsilon)}(\Gamma_2),
\end{equation}
where
\begin{align}
&\mathrm{mod}^{(\varepsilon)}(\Gamma)\nonumber\\
&\ge \frac{1}{C_{mod}}\frac{\Psi(\varepsilon)}{\Phi(\varepsilon)}\mathrm{cap}^{(\varepsilon)}\left(B^{(\varepsilon)}(y_0,\frac{1}{16}R),B^{(\varepsilon)}(z_0,\frac{1}{16}R);B^{(\varepsilon)}(x_0,A_1R)\right)\nonumber\\
&\ge \frac{1}{C_1C_{mod}}\frac{\Psi(\varepsilon)}{\Phi(\varepsilon)}\frac{\Phi(R)}{\Psi(R)}.\label{eq_LLC2}
\end{align}
For any $n\ge0$, by Proposition \ref{prop_ucap_Veps}, we have
$$\mathrm{cap}^{(\varepsilon)}\left(B^{(\varepsilon)}(x_0,\frac{1}{2A_{2}^{n+1}}R),V^{(\varepsilon)}\backslash B^{(\varepsilon)}(x_0,\frac{1}{2A_{2}^{n}}R)\right)\le C_2 \frac{\Phi \left(\frac{1}{2A_{2}^{n+1}}R\right)}{\Psi \left(\frac{1}{2A_{2}^{n+1}}R\right)},$$
then there exists $\phi_n\in l(V^{(\varepsilon)})$ satisfying that $\phi_n=1$ in $B^{(\varepsilon)}(x_0,\frac{1}{2A_{2}^{n+1}}R)$, $\phi_n=0$ on $V^{(\varepsilon)}\backslash B^{(\varepsilon)}(x_0,\frac{1}{2A_{2}^{n}}R)$, and
$$\mathbf{E}^{(\varepsilon)}(\phi_n)\le 2C_2 \frac{\Phi \left(\frac{1}{2A_{2}^{n+1}}R\right)}{\Psi \left(\frac{1}{2A_{2}^{n+1}}R\right)}.$$
Let $\phi=\frac{1}{M}\sum_{n=0}^{M-1}\phi_n$, then $\phi=1$ in $B^{(\varepsilon)}(x_0,\frac{1}{2A_{2}^{M}}R)$, $\phi=0$ on $V^{(\varepsilon)}\backslash B^{(\varepsilon)}(x_0,\frac{1}{2}R)$, and
\begin{align*}
&\mathrm{cap}^{(\varepsilon)}\left(B^{(\varepsilon)}(x_0, \frac{1}{2A_{2}^{M}}R),V^{(\varepsilon)}\backslash B^{(\varepsilon)}(x_0,\frac{1}{2}R)\right)\le\mathbf{E}^{(\varepsilon)}(\phi)\le \frac{1}{M^p}\sum_{n=0}^{M-1}\mathbf{E}^{(\varepsilon)}(\phi_n)\\
&\le \frac{1}{M^p}\sum_{n=0}^{M-1}2C_2 \frac{\Phi \left(\frac{1}{2A_{2}^{n+1}}R\right)}{\Psi \left(\frac{1}{2A_{2}^{n+1}}R\right)}\overset{(*)}{\scalebox{2}[1]{$\le$}}\frac{1}{M^p}\sum_{n=0}^{M-1}2C_2 C_{FVR}\frac{\Phi(R)}{\Psi(R)}=\frac{C_3}{M^{p-1}}\frac{\Phi(R)}{\Psi(R)},
\end{align*}
where $(*)$ follows from \ref{eq_FVR}, and $C_3=2C_2C_{FVR}$, which gives
\begin{align}
&\mathrm{mod}^{(\varepsilon)}(\Gamma_2)\nonumber\\
&\le C_{mod}\frac{\Psi(\varepsilon)}{\Phi(\varepsilon)}\mathrm{cap}^{(\varepsilon)}\left(B^{(\varepsilon)}(x_0, \frac{1}{2A_{2}^{M}}R),V^{(\varepsilon)}\backslash B^{(\varepsilon)}(x_0,\frac{1}{2}R)\right)\nonumber\\
&\le \frac{C_3C_{mod}}{M^{p-1}}\frac{\Psi(\varepsilon)}{\Phi(\varepsilon)}\frac{\Phi(R)}{\Psi(R)}.\label{eq_LLC3}
\end{align}
Combining (\ref{eq_LLC1})--(\ref{eq_LLC3}), we have
\begin{align*}
&\frac{1}{C_1C_{mod}}\frac{\Psi(\varepsilon)}{\Phi(\varepsilon)}\frac{\Phi(R)}{\Psi(R)}\le\mathrm{mod}^{(\varepsilon)}(\Gamma_1)+\frac{C_3C_{mod}}{M^{p-1}}\frac{\Psi(\varepsilon)}{\Phi(\varepsilon)}\frac{\Phi(R)}{\Psi(R)}.
\end{align*}
Take $M\ge1$ such that $M^{p-1}\ge2C_1C_3C_{mod}^2$, then
$$\mathrm{mod}^{(\varepsilon)}(\Gamma_1)\ge\frac{1}{2C_1C_{mod}}\frac{\Psi(\varepsilon)}{\Phi(\varepsilon)}\frac{\Phi(R)}{\Psi(R)}>0.$$

In summary, we have the desired result with $A=\max\{A_1,A_2^M\}$.
\end{proof}

Second, the proofs of Theorems \ref{thm_main_RSVR} and \ref{thm_HK} can be reduced to the following result.

\begin{proposition}\label{prop_CS_X}
Assume \ref{eq_RSVR}, \ref{eq_PI}, \ref{eq_ucap}. If either
\begin{enumerate}[label=(\arabic*)]
\item $p=2$, or 
\item \ref{eq_FVR},
\end{enumerate}
holds, then \ref{eq_CS} holds.
\end{proposition}

Given an $\varepsilon$-net $V^{(\varepsilon)}$, we associate the corresponding edge set $E^{(\varepsilon)}$, metric $d^{(\varepsilon)}$, measure $m^{(\varepsilon)}$, and $p$-energy $(\mathbf{E}^{(\varepsilon)},\mathbf{F}^{(\varepsilon)})$ as in the proof of Proposition \ref{prop_LLC}. Following \cite[Section 3]{DRY23}, we construct a cable system (also known as a metric graph) $X^{(\varepsilon)}$, endowed with a metric $d^{(\varepsilon)}$, a measure $\lambda^{(\varepsilon)}$, and a $p$-energy $(\mathsf{E}^{(\varepsilon)},\mathsf{F}^{(\varepsilon)})$ as follows. For each edge $(z_1,z_2)\in E^{(\varepsilon)}$, we ``replace" it with a closed cable, denoted by $[z_1,z_2]$, which is identified with the closed interval $[0,d^{(\varepsilon)}(z_1,z_2)]$. We also denote by $]z_1,z_2[=[z_1,z_2]\backslash\{z_1,z_2\}$ the corresponding open cable. Notice that $V^{(\varepsilon)}\subseteq X^{(\varepsilon)}$.

We extend the metric $d^{(\varepsilon)}$ on $V^{(\varepsilon)}$ to a metric $d^{(\varepsilon)}$ on $X^{(\varepsilon)}$ as follows. For any $x,y\in X^{(\varepsilon)}$, if $x, y$ lie on the same cable $[z_1,z_2]$, then let $d^{(\varepsilon)}(x,y)=|x-y|$, where $x,y$ on the RHS are interpreted as points in the interval $[0,d^{(\varepsilon)}(x,y)]$, otherwise, there exist distinct cables $[z_1,z_2]$, $[z_3,z_4]$ such that $x$ lies on $[z_1,z_2]$ and $y$ lies on $[z_3,z_4]$, let
\begin{align*}
&d^{(\varepsilon)}(x,y)=\min \left\{
\begin{array}{l}
d^{(\varepsilon)}(x,z_1)+d^{(\varepsilon)}(z_1,z_3)+d^{(\varepsilon)}(z_3,y),\\
d^{(\varepsilon)}(x,z_1)+d^{(\varepsilon)}(z_1,z_4)+d^{(\varepsilon)}(z_4,y),\\
d^{(\varepsilon)}(x,z_2)+d^{(\varepsilon)}(z_2,z_3)+d^{(\varepsilon)}(z_3,y),\\
d^{(\varepsilon)}(x,z_2)+d^{(\varepsilon)}(z_2,z_4)+d^{(\varepsilon)}(z_4,y)
\end{array}
\right\}.
\end{align*}
By (\ref{eq_bdddeg}), $d^{(\varepsilon)}$ on $X^{(\varepsilon)}$ is well-defined and $(X^{(\varepsilon)},d^{(\varepsilon)})$ is a complete unbounded locally compact separable geodesic metric space. We also use $B^{(\varepsilon)}(x,r)$ to denote the open ball in $X^{(\varepsilon)}$ centered at $x\in X^{(\varepsilon)}$ with radius $r>0$.

Let $\mathcal{H}^{(\varepsilon)}$ and $\lambda^{(\varepsilon)}$ be the unique positive Radon measures on $X^{(\varepsilon)}$ such that, for any $(z_1,z_2)\in E^{(\varepsilon)}$, the restriction of $\mathcal{H}^{(\varepsilon)}$ to the closed cable $[z_1,z_2]$ coincides with the one-dimensional Lebesgue measure on the closed interval $[0,d^{(\varepsilon)}(z_1,z_2)]$, and the restriction of $\lambda^{(\varepsilon)}$ to $[z_1,z_2]$ coincides with the one-dimensional Lebesgue measure on the same interval, scaled by the factor $\frac{\Phi(\varepsilon)}{d^{(\varepsilon)}(z_1,z_2)}$, in particular, $\mathcal{H}^{(\varepsilon)}([z_1,z_2])=d^{(\varepsilon)}(z_1,z_2)$ and $\lambda^{(\varepsilon)}([z_1,z_2])=\Phi(\varepsilon)$. Let
$$\Phi^{(\varepsilon)}(r)=
\begin{cases}
\Phi(\varepsilon)\frac{r}{\varepsilon}&\text{if }r\le\varepsilon,\\
\Phi(r)&\text{if }r>\varepsilon,
\end{cases}
$$
then by (\ref{eq_bdddeg}) and \ref{eq_VPhi}, we have
\begin{equation*}\label{eq_VPhieps}\tag*{$\mathsf{V}^{(\varepsilon)}(\Phi^{(\varepsilon)})$}
\lambda^{(\varepsilon)}(B^{(\varepsilon)}(x,r))\asymp \Phi^{(\varepsilon)}(r)\text{ for any }x\in X^{(\varepsilon)},r>0.
\end{equation*}
Let $\mathsf{u},\mathsf{v}\in l(X^{(\varepsilon)})$, and $(z_1,z_2)\in E^{(\varepsilon)}$. For any $x$ in the open cable $]z_1,z_2[$, we define
$$\nabla^{(\varepsilon)}\mathsf{u}(x)=\lim_{]z_1,z_2[\ni y\to x}\frac{\mathsf{u}(y)-\mathsf{u}(x)}{d^{(\varepsilon)}(y,z_1)-d^{(\varepsilon)}(x,z_1)}.$$
We define the directional derivative at $z_1$ in the direction of $z_2$ as
$$\nabla^{(\varepsilon)}_{z_2}\mathsf{u}(z_1)=\lim_{]z_1,z_2[\ni y\to z_1}\frac{\mathsf{u}(y)-\mathsf{u}(z_1)}{d^{(\varepsilon)}(y,z_1)}.$$
Note that the choice of the roles of $z_1$ and $z_2$ determines the sign of $\nabla^{(\varepsilon)}\mathsf{u}(x)$, but does not affect the quantities $\lvert\nabla^{(\varepsilon)}\mathsf{u}(x)\rvert$ and $\nabla^{(\varepsilon)}\mathsf{u}(x)\nabla^{(\varepsilon)}\mathsf{v}(x)$. For any measurable subset $A$ of $X^{(\varepsilon)}$, let
$$\lVert {|\nabla^{(\varepsilon)}\mathsf{u}|}\rVert_{L^\infty(X^{(\varepsilon)};\lambda^{(\varepsilon)})}=\esup_{x\in A\backslash V^{(\varepsilon)}}|\nabla^{(\varepsilon)}\mathsf{u}(x)|.$$
Since $\lambda^{(\varepsilon)}(V^{(\varepsilon)})=0$, the above definition makes sense even $\nabla^{(\varepsilon)}\mathsf{u}(x)$ is not defined for all $x\in V^{(\varepsilon)}$.

Let
\begin{align*}
&\mathsf{K}^{(\varepsilon)}=\left\{\mathsf{u}\in C_c(X^{(\varepsilon)})\middle|
\begin{array}{l}
\nabla^{(\varepsilon)}\mathsf{u}(x),\nabla^{(\varepsilon)}_{z_2}\mathsf{u}(z_1)\text{ exist for any }x\in ]z_1,z_2[,\\
\text{for any }(z_1,z_2)\in E^{(\varepsilon)},\lVert {|\nabla^{(\varepsilon)}\mathsf{u}|}\rVert_{L^\infty(X^{(\varepsilon)};\lambda^{(\varepsilon)})}<+\infty
\end{array}
\right\}.
\end{align*}
Let
\begin{align*}
\mathsf{E}^{(\varepsilon)}(\mathsf{u})&=\frac{\Phi(\varepsilon)}{\Psi(\varepsilon)}\sum_{(z_1,z_2)\in E^{(\varepsilon)}}d^{(\varepsilon)}(z_1,z_2)^{p-1}\int_{]z_1,z_2[}|\nabla^{(\varepsilon)} \mathsf{u}|^p\mathrm{d}\mathcal{H}^{(\varepsilon)},\\
\mathsf{F}^{(\varepsilon)}&=\text{the }\left(\mathsf{E}^{(\varepsilon)}(\cdot)^{1/p}+\lVert {\cdot}\rVert_{L^p(X^{(\varepsilon)};\lambda^{(\varepsilon)})}\right)\text{-closure of }\mathsf{K}^{(\varepsilon)}.
\end{align*}
By the classical Sobolev space theory on $\mathbb{R}$, we have $\mathsf{F}^{(\varepsilon)}\subseteq C(X^{(\varepsilon)})$. It is obvious that $(\mathsf{E}^{(\varepsilon)},\mathsf{F}^{(\varepsilon)})$ is a $p$-energy on $(X^{(\varepsilon)},d^{(\varepsilon)},\lambda^{(\varepsilon)})$ with a $p$-energy measure $\Gamma^{(\varepsilon)}$ given by
$$\Gamma^{(\varepsilon)}(\mathsf{u})(A)=\mathsf{E}^{(\varepsilon)}_A(\mathsf{u})=\frac{\Phi(\varepsilon)}{\Psi(\varepsilon)}\sum_{(z_1,z_2)\in E^{(\varepsilon)}}d^{(\varepsilon)}(z_1,z_2)^{p-1}\int_{A\cap ]z_1,z_2[}|\nabla^{(\varepsilon)} \mathsf{u}|^p\mathrm{d}\mathcal{H}^{(\varepsilon)}$$
for any $A\in\mathcal{B}(X^{(\varepsilon)})$, $\mathsf{u}\in\mathsf{F}^{(\varepsilon)}$. By definition, we have $\mathsf{E}^{(\varepsilon)}_{X^{(\varepsilon)}}=\mathsf{E}^{(\varepsilon)}$.

Let $A_1,A_2\in\mathcal{B}(X^{(\varepsilon)})$. We define the capacity between $A_1, A_2$ as
\begin{align*}
&\mathsf{cap}^{(\varepsilon)}(A_1,A_2)=\inf\left\{\mathsf{E}^{(\varepsilon)}(\mathsf{u}):\mathsf{u}\in\mathsf{F}^{(\varepsilon)},
\begin{array}{l}
\mathsf{u}=1\text{ in an open neighborhood of }A_1,\\
\mathsf{u}=0\text{ in an open neighborhood of }A_2
\end{array}
\right\}.
\end{align*}

Let
$$\Psi^{(\varepsilon)}(r)=
\begin{cases}
\Psi(\varepsilon)\left(\frac{r}{\varepsilon}\right)^{p}&\text{if }r\le\varepsilon,\\
\Psi(r)&\text{if }r>\varepsilon.
\end{cases}
$$
We show that the corresponding functional inequalities with scaling function $\Psi^{(\varepsilon)}$ also hold on $X^{(\varepsilon)}$. Firstly, we have the Poincar\'e inequality on $X^{(\varepsilon)}$ as follows.

\begin{proposition}\label{prop_PI_Xeps}
Assume \ref{eq_VPhi}, \ref{eq_PI}, \ref{eq_ucap}. Then there exist $C>0$, $A>1$ such that for any ball $B^{(\varepsilon)}(x_0,R)$ in $X^{(\varepsilon)}$, for any $\mathsf{u}\in\mathsf{F}^{(\varepsilon)}$, we have
\begin{align*}
&\int_{B^{(\varepsilon)}(x_0,R)}\lvert \mathsf{u}-\frac{1}{\lambda^{(\varepsilon)}(B^{(\varepsilon)}(x_0,R))}\int_{B^{(\varepsilon)}(x_0,R)}\mathsf{u}\mathrm{d}\lambda^{(\varepsilon)}\rvert^p\mathrm{d}\lambda^{(\varepsilon)}\\
&\le C\Psi^{(\varepsilon)}(R)\mathsf{E}_{B^{(\varepsilon)}(x_0,AR)}^{(\varepsilon)}(\mathsf{u}).\label{eq_PI_Xeps}\tag*{$\mathsf{PI}^{(\varepsilon)}(\Psi^{(\varepsilon)})$}
\end{align*}
\end{proposition}

\begin{proof}
If $R\le\varepsilon$, then the result follows directly from the Poincar\'e inequality on $\mathbb{R}$. We may assume that $R>\varepsilon$. For notational convenience, we write $\mathsf{u}_{W}=\frac{1}{\lambda^{(\varepsilon)}(W)}\int_{W}\mathsf{u}\mathrm{d}\lambda^{(\varepsilon)}$ for any measurable subset $W$ of $X^{(\varepsilon)}$ with $\lambda^{(\varepsilon)}(W)\in(0,+\infty)$. For any $c\in \mathbb{R}$, we have
\begin{align*}
&\int_{B^{(\varepsilon)}(x_0,R)}\lvert \mathsf{u}-c\rvert^p\mathrm{d}\lambda^{(\varepsilon)}\le\sum_{z\in B^{(\varepsilon)}(x_0,R+\frac{5}{2}\varepsilon)\cap V^{(\varepsilon)}}\int_{B^{(\varepsilon)}(z,\frac{5}{4}\varepsilon)}|\mathsf{u}-c|^p\mathrm{d}\lambda^{(\varepsilon)}\\
&\le2^{p-1}\sum_{z\in B^{(\varepsilon)}(x_0,R+\frac{5}{2}\varepsilon)\cap V^{(\varepsilon)}}\int_{B^{(\varepsilon)}(z,\frac{5}{4}\varepsilon)}|\mathsf{u}-\mathsf{u}_{B^{(\varepsilon)}(z,\frac{5}{4}\varepsilon)}|^p\mathrm{d}\lambda^{(\varepsilon)}\\
&+2^{p-1}\sum_{z\in B^{(\varepsilon)}(x_0,R+\frac{5}{2}\varepsilon)\cap V^{(\varepsilon)}}\int_{B^{(\varepsilon)}(z,\frac{5}{4}\varepsilon)}|\mathsf{u}_{B^{(\varepsilon)}(z,\frac{5}{4}\varepsilon)}-c|^p\mathrm{d}\lambda^{(\varepsilon)}\\
&=I_1+I_2.
\end{align*}
By the Poincar\'e inequality on $\mathbb{R}$, there exists $A_1\ge1$ such that
\begin{align*}
&\int_{B^{(\varepsilon)}(z,\frac{5}{4}\varepsilon)}|\mathsf{u}-\mathsf{u}_{B^{(\varepsilon)}(z,\frac{5}{4}\varepsilon)}|^p\mathrm{d}\lambda^{(\varepsilon)}\lesssim \Psi(\varepsilon)\mathsf{E}^{(\varepsilon)}_{{B^{(\varepsilon)}(z,\frac{5}{4}A_1\varepsilon)}}(\mathsf{u}),
\end{align*}
hence
\begin{align*}
&I_1\lesssim \Psi(\varepsilon)\sum_{z\in B^{(\varepsilon)}(x_0,R+\frac{5}{2}\varepsilon)\cap V^{(\varepsilon)}}\mathsf{E}^{(\varepsilon)}_{{B^{(\varepsilon)}(z,\frac{5}{4}A_1\varepsilon)}}(\mathsf{u})\\
&\overset{(*)}{\scalebox{2}[1]{$\lesssim$}}\Psi(\varepsilon)\mathsf{E}^{(\varepsilon)}_{\cup_{z\in B^{(\varepsilon)}(x_0,R+\frac{5}{2}\varepsilon)\cap V^{(\varepsilon)}}{B^{(\varepsilon)}(z,\frac{5}{4}A_1\varepsilon)}}(\mathsf{u})\le \Psi(R)\mathsf{E}^{(\varepsilon)}_{{B^{(\varepsilon)}(x_0,5A_1R)}}(\mathsf{u}),
\end{align*}
where $(*)$ follows from \ref{eq_VPhieps}. For any $z\in V^{(\varepsilon)}$, let $\mathbf{u}(z)=\mathsf{u}_{B^{(\varepsilon)}(z,\frac{5}{4}\varepsilon)}$. Take $y_0\in V^{(\varepsilon)}$ such that $d^{(\varepsilon)}(x_0,y_0)<\frac{5}{2}\varepsilon$, then
\begin{align*}
I_2\asymp\int_{B^{(\varepsilon)}(x_0,R+\frac{5}{2}\varepsilon)\cap V^{(\varepsilon)}} \lvert \mathbf{u}-c\rvert^p \mathrm{d}m^{(\varepsilon)}\le\int_{B^{(\varepsilon)}(y_0,6R)} \lvert \mathbf{u}-c\rvert^p \mathrm{d}m^{(\varepsilon)}.
\end{align*}
Let $c=\frac{1}{m^{(\varepsilon)}(B^{(\varepsilon)}(y_0,6R))}\int_{B^{(\varepsilon)}(y_0,6R)}\mathbf{u}\mathrm{d} m^{(\varepsilon)}$, then by \ref{eq_PI_Veps} on $V^{(\varepsilon)}$, we have
$$I_2\lesssim \Psi(R)\mathbf{E}^{(\varepsilon)}_{B^{(\varepsilon)}(y_0,6A_2R)}(\mathbf{u}),$$
where $A_2$ is the positive constant appearing therein. Recall that
$$\mathbf{E}^{(\varepsilon)}_{B^{(\varepsilon)}(y_0,6A_2R)}(\mathbf{u})=\frac{\Phi(\varepsilon)}{\Psi(\varepsilon)}\sum_{\substack{(z_1,z_2)\in E^{(\varepsilon)}\\z_1,z_2\in B^{(\varepsilon)}(y_0,6A_2R)\cap V^{(\varepsilon)}}}\lvert\mathbf{u}(z_1)-\mathbf{u}(z_2)\rvert^p,$$
where by the Poincar\'e inequality on $\mathbb{R}$, we have
\begin{align*}
\lvert \mathbf{u}(z_1)-\mathbf{u}(z_2)\rvert^p=\lvert \mathsf{u}_{B^{(\varepsilon)}(z_1,\frac{5}{4}\varepsilon)}-\mathsf{u}_{B^{(\varepsilon)}(z_2,\frac{5}{4}\varepsilon)}\rvert^p\lesssim \frac{\Psi(\varepsilon)}{\Phi(\varepsilon)}\mathsf{E}^{(\varepsilon)}_{B^{(\varepsilon)}(z_1,4A_1\varepsilon)}(\mathsf{u}),
\end{align*}
hence
\begin{align*}
&\mathbf{E}^{(\varepsilon)}_{B^{(\varepsilon)}(y_0,6A_2R)}(\mathbf{u})\lesssim \sum_{z\in B^{(\varepsilon)}(y_0,6A_2R+4A_1\varepsilon)}\mathsf{E}^{(\varepsilon)}_{B^{(\varepsilon)}(z,4A_1\varepsilon)}(\mathsf{u})\\
&\overset{(**)}{\scalebox{2}[1]{$\lesssim$}}\mathsf{E}^{(\varepsilon)}_{\cup_{{z\in B^{(\varepsilon)}(y_0,6A_2R+4A_1\varepsilon)}}B^{(\varepsilon)}(z,4A_1\varepsilon)}(\mathsf{u})\le\mathsf{E}^{(\varepsilon)}_{B^{(\varepsilon)}(x_0,6(A_1+A_2+1)R)}(\mathsf{u}),
\end{align*}
where $(**)$ follows from \ref{eq_VPhieps}. Therefore
$$I_2\lesssim \Psi(R)\mathsf{E}^{(\varepsilon)}_{B^{(\varepsilon)}(x_0,6(A_1+A_2+1)R)}(\mathsf{u}).$$
In summary, we have
\begin{align*}
&\int_{B^{(\varepsilon)}(x_0,R)}\lvert \mathsf{u}-\frac{1}{\lambda^{(\varepsilon)}(B^{(\varepsilon)}(x_0,R))}\int_{B^{(\varepsilon)}(x_0,R)}\mathsf{u}\mathrm{d}\lambda^{(\varepsilon)}\rvert^p\mathrm{d}\lambda^{(\varepsilon)}\\
&\le2^p\int_{B^{(\varepsilon)}(x_0,R)}\lvert \mathsf{u}-c\rvert^p\mathrm{d}\lambda^{(\varepsilon)}\le2^p \left(I_1+I_2\right)\lesssim \Psi(R)\mathsf{E}^{(\varepsilon)}_{B^{(\varepsilon)}(x_0,6(A_1+A_2+1)R)}(\mathsf{u}).
\end{align*}
\end{proof}

Secondly, we have the capacity upper bound on $X^{(\varepsilon)}$ as follows.

\begin{proposition}\label{prop_ucap_Xeps}
Assume \ref{eq_VPhi}, \ref{eq_PI}, \ref{eq_ucap}. Then there exist $C>0$, $A>1$ such that for any ball $B^{(\varepsilon)}(x_0,R)$ in $X^{(\varepsilon)}$, we have
\begin{equation*}\label{eq_ucap_Xeps}\tag*{$\mathsf{cap}^{(\varepsilon)}(\Psi^{(\varepsilon)})_\le$}
\mathsf{cap}^{(\varepsilon)}(B^{(\varepsilon)}(x_0,R),X^{(\varepsilon)}\backslash B^{(\varepsilon)}(x_0,AR))\le C \frac{\Phi^{(\varepsilon)}(R)}{\Psi^{(\varepsilon)}(R)}.
\end{equation*}
\end{proposition}

\begin{proof}
If $R\le\varepsilon$, then the result follows directly from the capacity upper bound on $\mathbb{R}$. We may assume that $R>\varepsilon$. For any $\mathbf{u}\in l(V^{(\varepsilon)})$, we define $\mathsf{u}\in l(X^{(\varepsilon)})$ by linear interpolation on each cable as follows: for any $(z_1,z_2)\in E^{(\varepsilon)}$ and any $x\in[z_1,z_2]$, let
\begin{equation}\label{eq_linear_inter}
\mathsf{u}(x)=\frac{d^{(\varepsilon)}(x,z_1)\mathbf{u}(z_2)+d^{(\varepsilon)}(x,z_2)\mathbf{u}(z_1)}{d^{(\varepsilon)}(z_1,z_2)}.
\end{equation}
Then $\mathsf{u}$ is well-defined, $\mathsf{u}|_{V^{(\varepsilon)}}=\mathbf{u}$, and for any $(z_1,z_2)\in E^{(\varepsilon)}$, we have
$$d^{(\varepsilon)}(z_1,z_2)^{p-1}\int_{]z_1,z_2[}|\nabla^{(\varepsilon)}\mathsf{u}|^p\mathrm{d}\mathcal{H}^{(\varepsilon)}=|\mathbf{u}(z_1)-\mathbf{u}(z_2)|^p,$$
which gives $\mathsf{E}^{(\varepsilon)}(\mathsf{u})=\mathbf{E}^{(\varepsilon)}(\mathbf{u})\in[0,+\infty]$.

Let $A_1,C_1$ be the constants appearing in \ref{eq_ucap_Veps} on $V^{(\varepsilon)}$. Take $y_0\in V^{(\varepsilon)}$ such that $d^{(\varepsilon)}(x_0,y_0)<\frac{5}{2}\varepsilon$, then there exists $\mathbf{u}\in l(V^{(\varepsilon)})$ with $\mathbf{u}=1$ in $B^{(\varepsilon)}(y_0,R+5\varepsilon)$ and $\mathbf{u}=0$ on $V^{(\varepsilon)}\backslash B^{(\varepsilon)}(y_0,A_1(R+5\varepsilon))$ such that $\mathbf{E}^{(\varepsilon)}(\mathbf{u})\le 2C_1 \frac{\Phi(R+5\varepsilon)}{\Psi(R+5\varepsilon)}$. Let $\mathsf{u}\in l(X^{(\varepsilon)})$ be given by (\ref{eq_linear_inter}), then $\mathsf{u}=1$ in $B^{(\varepsilon)}(x_0,R)$, $\mathsf{u}=0$ on $X^{(\varepsilon)}\backslash B^{(\varepsilon)}(x_0,11A_1R)$, and
\begin{align*}
&\mathsf{cap}^{(\varepsilon)}(B^{(\varepsilon)}(x_0,R),X^{(\varepsilon)}\backslash B^{(\varepsilon)}(x_0,11A_1R))\le\mathsf{E}^{(\varepsilon)}(\mathsf{u})=\mathbf{E}^{(\varepsilon)}(\mathbf{u})\le2C_1C_\Phi^3 \frac{\Phi(R)}{\Psi(R)}.\\
\end{align*}
\end{proof}

By \cite[Proposition 2.1]{Yan25a}, under \ref{eq_VPhi}, \ref{eq_PI}, \ref{eq_ucap}, there exists $C>0$ such that
$$\frac{\Psi(R)}{\Psi(r)}\le C \left(\frac{R}{r}\right)^{p-1}\frac{\Phi(R)}{\Phi(r)}\text{ for any }r\le R.$$
If \ref{eq_RSVR} also holds, then
$$\frac{1}{CC_{RSVR}}\le \left(\frac{R}{r}\right)^{p-1+\tau}\text{ for any }r\le R,$$
hence $\tau\in[1-p,1)$. Without loss of generality, by assuming $C_{RSVR}\ge1$, it is easy to verify that $\Phi^{(\varepsilon)}, \Psi^{(\varepsilon)}$ also satisfy that
\begin{equation*}\label{eq_RSVR_Xeps}\tag*{$\mathsf{RSVR}^{(\varepsilon)}(\Phi^{(\varepsilon)},\Psi^{(\varepsilon)})$}
\frac{1}{C_{RSVR}}\left(\frac{r}{R}\right)^{\tau}\frac{\Phi^{(\varepsilon)}(R)}{\Phi^{(\varepsilon)}(r)}\le\frac{\Psi^{(\varepsilon)}(R)}{\Psi^{(\varepsilon)}(r)}\text{ for any }r\le R,
\end{equation*}
in particular, in the case $r\in(0,\varepsilon]$, we require the fact that $\tau\in[1-p,1)$.

Recall that the Hausdorff $1$-content $\mathcal{H}_1^\infty(A)$ of a set $A\subseteq X^{(\varepsilon)}$ is given by
$$\mathcal{H}_1^{\infty}(A)=\inf\sum_{i}r_i,$$
where the infimum is taken over all countable covers of the set $A$ by balls $B^{(\varepsilon)}_i$ in $X^{(\varepsilon)}$ with radius $r_i$. If $A$ is connected, then its Hausdorff $1$-content is comparable to its diameter, indeed,
\begin{equation}\label{eq_Haus_diam}
\frac{1}{2}\mathrm{diam}(A)\le\mathcal{H}_1^\infty(A)\le\mathrm{diam}(A),
\end{equation}
where the lower bound follows from \cite[Lemma 2.6.1]{BBI01} and the fact that $(X^{(\varepsilon)},d^{(\varepsilon)})$ is a geodesic metric space.

The following version of the capacity lower bound follows by an argument analogous to that of Lemma \ref{lem_lcap_Veps}.

\begin{lemma}\label{lem_lcap_Xeps}
Assume \ref{eq_RSVR_Xeps}, \ref{eq_PI_Xeps}. Then there exist $C>0$, $A>1$ such that for any ball $B^{(\varepsilon)}(x_0,R)$ in $X^{(\varepsilon)}$ and any subsets $E,F$ of $B^{(\varepsilon)}(x_0,R)$ satisfying that
$$\min\{\mathcal{H}_1^\infty(E),\mathcal{H}_1^\infty(F)\}\ge \frac{1}{4}R,$$
for any $\mathsf{u}\in\mathsf{F}^{(\varepsilon)}\subseteq C(X^{(\varepsilon)})$ with $\mathsf{u}\ge1$ on $E$ and $\mathsf{u}\le0$ on $F$, we have
$$\mathsf{E}^{(\varepsilon)}_{B^{(\varepsilon)}(x_0,AR)}(\mathsf{u})\ge \frac{1}{C}\frac{\Phi^{(\varepsilon)}(R)}{\Psi^{(\varepsilon)}(R)}.$$
\end{lemma}

\begin{proof}
For notational convenience, we write $\mathsf{u}_W=\frac{1}{\lambda^{(\varepsilon)}(W)}\int_W\mathsf{u}\mathrm{d}\lambda^{(\varepsilon)}$ for any measurable subset $W$ of $X^{(\varepsilon)}$ with $\lambda^{(\varepsilon)}(W)\in(0,+\infty)$.
If there exist $y\in E$ and $z\in F$ such that $|\mathsf{u}(y)-\mathsf{u}_{B^{(\varepsilon)}(y,R)}|\le \frac{1}{5}$ and $|\mathsf{u}(z)-\mathsf{u}_{B^{(\varepsilon)}(z,R)}|\le \frac{1}{5}$, then $\lvert \mathsf{u}_{B^{(\varepsilon)}(y,R)}-\mathsf{u}_{B^{(\varepsilon)}(z,R)}\rvert\ge \frac{3}{5}$. By \ref{eq_PI_Xeps}, we have 
\begin{align*}
\lvert \mathsf{u}_{B^{(\varepsilon)}(y,R)}-\mathsf{u}_{B^{(\varepsilon)}(z,R)}\rvert^p\lesssim \frac{\Psi^{(\varepsilon)}(R)}{\Phi^{(\varepsilon)}(R)}\mathsf{E}^{(\varepsilon)}_{B^{(\varepsilon)}(x_0,2A_1R)}(\mathsf{u}),
\end{align*}
where $A_1>1$ is the positive constant appearing therein, hence
$$\mathsf{E}^{(\varepsilon)}_{B^{(\varepsilon)}(x_0,2A_1R)}(\mathsf{u})\gtrsim \frac{\Phi^{(\varepsilon)}(R)}{\Phi^{(\varepsilon)}(R)}.$$

Assume that no such $y,z$ as above exist, then without loss of generality, we may assume that $|\mathsf{u}(y)-\mathsf{u}_{B^{(\varepsilon)}(y,R)}|>\frac{1}{5}$ for any $y\in E$, then by \ref{eq_PI_Xeps} again, we have
\begin{align*}
\frac{1}{5}&<|\mathsf{u}(y)-\mathsf{u}_{B^{(\varepsilon)}(y,R)}|\le\sum_{n=0}^{+\infty}|\mathsf{u}_{B^{(\varepsilon)}(y,\frac{1}{A_1^{n+1}}R)}-\mathsf{u}_{B^{(\varepsilon)}(y,\frac{1}{A_1^{n}}R)}|\\
&\lesssim\sum_{n=0}^{+\infty}\left(\frac{\Psi^{(\varepsilon)}(\frac{1}{A_1^{n}}R)}{\Phi^{(\varepsilon)}(\frac{1}{A_1^{n}}R)}\mathsf{E}^{(\varepsilon)}_{B^{(\varepsilon)}(y,\frac{1}{A_1^{n-1}}R)}(\mathsf{u})\right)^{1/p},
\end{align*}
that is,
$$\sum_{n=0}^{+\infty}\left(\frac{\Psi^{(\varepsilon)}(\frac{1}{A_1^{n}}R)}{\Phi^{(\varepsilon)}(\frac{1}{A_1^{n}}R)}\mathsf{E}^{(\varepsilon)}_{B^{(\varepsilon)}(y,\frac{1}{A_1^{n-1}}R)}(\mathsf{u})\right)^{1/p}\ge C_1,$$
where $C_1>0$ is some constant. Hence there exists $n_y\ge0$ such that
\begin{equation}\label{eq_lcap_RSVR}
\frac{\Psi^{(\varepsilon)}(\frac{1}{A_1^{n_y}}R)}{\Phi^{(\varepsilon)}(\frac{1}{A_1^{n_y}}R)}\mathsf{E}^{(\varepsilon)}_{B^{(\varepsilon)}(y,\frac{1}{A_1^{n_y-1}}R)}(\mathsf{u})\ge \left((1-A_1^{\frac{\tau-1}{p}})C_1\right)^pA_1^{(\tau-1)n_y}.
\end{equation}
Otherwise
$$C_1\le\sum_{n=0}^{+\infty}\left(\frac{\Psi^{(\varepsilon)}(\frac{1}{A_1^{n}}R)}{\Phi^{(\varepsilon)}(\frac{1}{A_1^{n}}R)}\mathsf{E}^{(\varepsilon)}_{B^{(\varepsilon)}(y,\frac{1}{A_1^{n-1}}R)}(\mathsf{u})\right)^{1/p}<(1-A_1^{\frac{\tau-1}{p}})C_1\sum_{n=0}^{+\infty}A_1^{\frac{\tau-1}{p}n}=C_1,$$
contradiction. By the $5B$-covering lemma (see \cite[Theorem 1.2]{Hei01}), there exists a countable family of disjoint balls $\{B^{(\varepsilon)}(y_k,\frac{1}{A_1^{n_k-1}}R)\}_k$ such that
$$E\subseteq\bigcup_kB^{(\varepsilon)}(y_k,5\frac{1}{A_1^{n_k-1}}R).$$
Hence
\begin{align*}
\frac{1}{4}R&\le\mathcal{H}_1^\infty(E)\le\sum_k5\frac{1}{A_1^{n_k-1}}R=5A_1R\sum_k\left(\frac{1}{A_1^{n_k}}\right)^{1-\tau}\left(\frac{1}{A_1^{n_k}}\right)^\tau\\
&\overset{(*)}{\scalebox{2}[1]{$\le$}}5A_1R\sum_kA_1^{(\tau-1)n_k}\left(C_{RSVR}\frac{\Psi^{(\varepsilon)}(R)}{\Phi^{(\varepsilon)}(R)}\frac{\Phi^{(\varepsilon)}(\frac{1}{A_1^{n_k}}R)}{\Psi^{(\varepsilon)}(\frac{1}{A_1^{n_k}}R)}\right)\\
&\overset{(**)}{\scalebox{2}[1]{$\le$}}5A_1C_{RSVR}R\frac{\Psi^{(\varepsilon)}(R)}{\Phi^{(\varepsilon)}(R)}\sum_k \frac{1}{\left((1-A_1^{{\frac{\tau-1}{p}}})C_1\right)^p}\mathsf{E}^{(\varepsilon)}_{B^{(\varepsilon)}(y_k,\frac{1}{A_1^{n_k-1}}R)}(\mathsf{u})\\
&=C_2R\frac{\Psi^{(\varepsilon)}(R)}{\Phi^{(\varepsilon)}(R)}\sum_k\mathsf{E}^{(\varepsilon)}_{B^{(\varepsilon)}(y_k,\frac{1}{A_1^{n_k-1}}R)}(\mathsf{u})\le C_2R\frac{\Psi^{(\varepsilon)}(R)}{\Phi^{(\varepsilon)}(R)}\mathsf{E}^{(\varepsilon)}_{B^{(\varepsilon)}(x_0,2A_1R)}(\mathsf{u}),
\end{align*}
where $(*)$ follows from \ref{eq_RSVR_Xeps}, $(**)$ follows from (\ref{eq_lcap_RSVR}), the last inequality follows from the disjointness of the balls and the fact that they are all contained in $B^{(\varepsilon)}(x_0,2A_1R)$, and $C_2=\frac{5A_1C_{RSVR}}{\left((1-A_1^{{\frac{\tau-1}{p}}})C_1\right)^p}$, then
$$\mathsf{E}^{(\varepsilon)}_{B^{(\varepsilon)}(x_0,2A_1R)}(\mathsf{u})\ge \frac{1}{4C_2}\frac{\Phi^{(\varepsilon)}(R)}{\Psi^{(\varepsilon)}(R)}.$$
\end{proof}

We have the elliptic Harnack inequality on $X^{(\varepsilon)}$ as follows. The proof follows the same argument as in \cite[THEOREM 4.3]{Hol03}, \cite[THEOREM 1.1]{Mur23a}, and \cite[Theorem 5.4]{MS25}.

\begin{proposition}\label{prop_EHI_Xeps}
Assume \ref{eq_RSVR_Xeps}, \ref{eq_PI_Xeps}, \ref{eq_ucap_Xeps}. Then there exist $C>0$, $A>1$ such that for any ball $B^{(\varepsilon)}(x_0,R)$ in $X^{(\varepsilon)}$, for any $\mathsf{h}\in\mathsf{F}^{(\varepsilon)}\subseteq C(X^{(\varepsilon)})$ which is non-negative harmonic in $B^{(\varepsilon)}(x_0,AR)$, we have
\begin{equation*}\label{eq_EHI_Xeps}\tag*{$\mathsf{EHI}^{(\varepsilon)}$}
\sup_{B^{(\varepsilon)}(x_0,\frac{1}{2}R)}\mathsf{h}\le C\inf_{B^{(\varepsilon)}(x_0,\frac{1}{2}R)}\mathsf{h}.
\end{equation*}
\end{proposition}

\begin{proof}
Let $A>1$ be chosen later. Let $\mathsf{h}\in \mathsf{F}^{(\varepsilon)}\subseteq C(X^{(\varepsilon)})$ be non-negative harmonic in $B^{(\varepsilon)}(x_0,AR)$. Let $M=\sup_{B^{(\varepsilon)}(x_0,\frac{1}{2}R)}\mathsf{h}$ and $m=\inf_{B^{(\varepsilon)}(x_0,\frac{1}{2}R)}\mathsf{h}$, then by the continuity of $\mathsf{h}$, we have $0\le m\le M<+\infty$. Without loss of generality, we may assume that $m<M$.

Let $E=\{\mathsf{h}\ge M\}$ and $F=\{\mathsf{h}\le m\}$, then by the classical maximum principle, the set $E$ contains a connected subset $E_0$ that intersects both $\partial B^{(\varepsilon)}(x_0,\frac{1}{2}R)$ and $\partial B^{(\varepsilon)}(x_0,R)$, hence $\mathrm{diam}(E_0)\ge \frac{1}{2}R$, by (\ref{eq_Haus_diam}), we have $\mathcal{H}_1^\infty(E)\ge \mathcal{H}_1^{\infty}(E_0)\ge \frac{1}{4}R$. Similarly, we also have $\mathcal{H}_1^\infty(F)\ge \frac{1}{4}R$. For any $\delta>0$, let
$$\mathsf{u}=\left(\frac{\log(\mathsf{h}+\delta)-\log(m+\delta)}{\log(M+\delta)-\log(m+\delta)}\vee0\right)\wedge1,$$
then $\mathsf{u}=\Upsilon(\mathsf{h})$, where $\Upsilon:\mathbb{R}\to \mathbb{R}$ is given by
$$\Upsilon(t)=\left(\frac{\log(t+\delta)-\log(m+\delta)}{\log(M+\delta)-\log(m+\delta)}\vee0\right)\wedge1=\frac{\log \left((t\vee m)\wedge M+\delta\right)-\log(m+\delta)}{\log(M+\delta)-\log(m+\delta)}.$$
Since $\Upsilon|_{(-\infty,m]}=0$, $\Upsilon|_{[M,+\infty)}=1$, and $\Upsilon'|_{(m,M)}\in[0,\frac{1}{\log(M+\delta)-\log(m+\delta)}\frac{1}{m+\delta}]$, by the Markovian property, we have $\mathsf{u}=\Upsilon(\mathsf{h})\in \mathsf{F}^{(\varepsilon)}$. It is obvious that $\mathsf{u}=1$ on $E$ and $\mathsf{u}=0$ on $F$. By Lemma \ref{lem_lcap_Xeps}, we have
$$\mathsf{E}^{(\varepsilon)}_{B^{(\varepsilon)}(x_0,A_1R)}(\mathsf{u})\ge \frac{1}{C_1}\frac{\Phi^{(\varepsilon)}(R)}{\Psi^{(\varepsilon)}(R)},$$
where $A_1,C_1$ are the positive constants appearing therein. By \ref{eq_ucap_Xeps}, there exists $\phi\in \mathsf{F}^{(\varepsilon)}$ with $\phi=1$ in $B^{(\varepsilon)}(x_0,A_1R)$ and $\phi=0$ on $X^{(\varepsilon)}\backslash B^{(\varepsilon)}(x_0,A_1A_2R)$ such that
$$\mathsf{E}^{(\varepsilon)}(\phi)\le 2C_2 \frac{\Phi^{(\varepsilon)}(R)}{\Psi^{(\varepsilon)}(R)},$$
where $A_2, C_2$ are the positive constants appearing therein. Let $A=A_1A_2$. By the same argument as in (\ref{eq_energy_log}) in the proof of Lemma \ref{lem_loguBMO}, and using the fact that $\mathsf{h}$ is harmonic in $B^{(\varepsilon)}(x_0,AR)$, we have
$$\mathsf{E}^{(\varepsilon)}_{B^{(\varepsilon)}(x_0,A_1R)}\left(\log(\mathsf{h}+\delta)\right)\le\int_{X^{(\varepsilon)}}\phi^p\mathrm{d}\Gamma^{(\varepsilon)}\left(\log (\mathsf{h}+\delta)\right)\le 2 \left(\frac{p}{p-1}\right)^pC_2 \frac{\Phi^{(\varepsilon)}(R)}{\Psi^{(\varepsilon)}(R)}.$$
Therefore, we have
\begin{align*}
&\frac{1}{C_1}\frac{\Phi^{(\varepsilon)}(R)}{\Psi^{(\varepsilon)}(R)}\le \mathsf{E}^{(\varepsilon)}_{B^{(\varepsilon)}(x_0,A_1R)}(\mathsf{u})\\
&\le \frac{1}{\left(\log \frac{M+\delta}{m+\delta}\right)^p}\mathsf{E}^{(\varepsilon)}_{B^{(\varepsilon)}(x_0,A_1R)}\left(\log(\mathsf{h}+\delta)\right)\\
&\le \frac{1}{\left(\log \frac{M+\delta}{m+\delta}\right)^p}2 \left(\frac{p}{p-1}\right)^pC_2 \frac{\Phi^{(\varepsilon)}(R)}{\Psi^{(\varepsilon)}(R)},
\end{align*}
which gives
$$\frac{M+\delta}{m+\delta}\le\exp \left(\frac{p}{p-1}\left(2C_1C_2\right)^{1/p}\right).$$
Letting $\delta\downarrow0$, we have the desired result.
\end{proof}

We have the cutoff Sobolev inequality on $X^{(\varepsilon)}$ as follows.

\begin{proposition}\label{prop_CS_Xeps}
Assume \ref{eq_RSVR}, \ref{eq_PI}, \ref{eq_ucap}. If either
\begin{enumerate}[label=(\arabic*)]
\item $p=2$, or
\item \ref{eq_FVR},
\end{enumerate}
holds, then there exist $C_{1}, C_{2}>0$, $A>1$ such that for any ball $B^{(\varepsilon)}(x_0,R)$ in $X^{(\varepsilon)}$, there exists a cutoff function $\phi\in\mathsf{F}^{(\varepsilon)}$ for $B^{(\varepsilon)}(x_0,R)\subseteq B^{(\varepsilon)}(x_0,AR)$ such that for any $\mathsf{f}\in\mathsf{F}^{(\varepsilon)}\subseteq C(X^{(\varepsilon)})$, we have
\begin{align*}
&\int_{B^{(\varepsilon)}(x_0,AR)}|{\mathsf{f}}|^p\mathrm{d}\Gamma^{(\varepsilon)}(\phi)\\
&\le C_{1}\int_{B^{(\varepsilon)}(x_0,AR)}\mathrm{d}\Gamma^{(\varepsilon)}(\mathsf{f})+\frac{C_{2}}{\Psi^{(\varepsilon)}(R)}\int_{B^{(\varepsilon)}(x_0,AR)}|\mathsf{f}|^p\mathrm{d} \lambda^{(\varepsilon)}.\label{eq_CS_Xeps}\tag*{$\mathsf{CS}^{(\varepsilon)}(\Psi^{(\varepsilon)})$}
\end{align*}
\end{proposition}

\begin{proof}
By Propositions \ref{prop_PI_Xeps}, \ref{prop_ucap_Xeps}, \ref{prop_EHI_Xeps}, we have \ref{eq_PI_Xeps}, \ref{eq_ucap_Xeps}, \ref{eq_EHI_Xeps}, respectively. Since $\mathsf{F}^{(\varepsilon)}\subseteq C(X^{(\varepsilon)})$, we have the corresponding \hyperlink{eq_LSC}{$\mathsf{LSC}^{(\varepsilon)}$} holds trivially on $X^{(\varepsilon)}$. Since any ball in $X^{(\varepsilon)}$ is a finite union of cable subintervals, \hyperref[eq_BSP]{$\mathsf{BSP}^{(\varepsilon)}$} on $X^{(\varepsilon)}$ follows from the corresponding one-dimensional result for bounded intervals. By \ref{eq_VPhieps} and \ref{eq_PI_Xeps}, applying the same line as the proof of \cite[Lemma 5.1 (a)]{BB04}, we have the capacity lower bound \hyperref[eq_lcap]{$\mathsf{cap}^{(\varepsilon)}(\Psi^{(\varepsilon)})_\ge$} on $X^{(\varepsilon)}$, then the two-sided capacity bounds \hyperlink{eq_cap}{$\mathsf{cap}^{(\varepsilon)}(\Psi^{(\varepsilon)})$} hold on $X^{(\varepsilon)}$.

Firstly, assume $p=2$. By \cite[Theorem 3.14]{GH14} and \cite[THEOREM 1.2]{GHL15}, under \ref{eq_VPhieps}, the conjunction of \ref{eq_EHI_Xeps} and \hyperlink{eq_cap}{$\mathsf{cap}^{(\varepsilon)}(\Psi^{(\varepsilon)})$} is equivalent to the conjunction of \ref{eq_PI_Xeps} and \ref{eq_CS_Xeps}, hence we have the desired result.

Secondly, assume \ref{eq_FVR}. By Proposition \ref{prop_LLC}, we have \hyperlink{eq_LLC}{\text{LLC}} on $X$, which implies the corresponding \hyperlink{eq_LLC}{$\mathsf{LLC}^{(\varepsilon)}$} on $X^{(\varepsilon)}$ for all balls with radii $\gtrsim\varepsilon$. Combining this with \ref{eq_VPhieps}, \ref{eq_EHI_Xeps}, by a standard chaining argument as in \cite[Theorem 3.1]{Bla01} and \cite[LEMMA 6.3]{GS05}, we have the corresponding \hyperref[eq_EHIann]{$\mathsf{EHI}_{\mathsf{ann}}^{(\varepsilon)}$} on $X^{(\varepsilon)}$ for all balls with radii $\gtrsim\varepsilon$. Since any ball in $X^{(\varepsilon)}$ with radius $\lesssim\varepsilon$ is the union of a bounded number of cable subintervals, where the bound depends only on \ref{eq_VPhieps}, we have \hyperref[eq_EHIann]{$\mathsf{EHI}_{\mathsf{ann}}^{(\varepsilon)}$} on $X^{(\varepsilon)}$ for such balls follows from the corresponding one-dimensional result on $\mathbb{R}$. By Theorem \ref{thm_main}, we have the desired result.
\end{proof}

We ``transfer" the cutoff Sobolev inequality from $X^{(\varepsilon)}$ to $V^{(\varepsilon)}$ as follows. The proof follows essentially the same idea as in \cite[Proposition 3.4]{BB04}.

\begin{proposition}\label{prop_CS_Veps}
Assume the assumptions of Proposition \ref{prop_CS_Xeps} are satisfied. Then there exist $C_1, C_2>0$, $A>1$ such that for any ball $B^{(\varepsilon)}(x_0,R)$ in $V^{(\varepsilon)}$ with $R>\varepsilon$, there exists a cutoff function $\phi\in\mathbf{F}^{(\varepsilon)}$ for $B^{(\varepsilon)}(x_0,R)\subseteq B^{(\varepsilon)}(x_0,AR)$ such that for any $\mathbf{f}\in\mathbf{F}^{(\varepsilon)}$, we have
\begin{align*}
&\sum_{z_1\in B^{(\varepsilon)}(x_0,AR)}|\mathbf{f}(z_1)|^p \frac{\Phi(\varepsilon)}{\Psi(\varepsilon)}\sum_{z_2:(z_1,z_2)\in E^{(\varepsilon)}}|\phi(z_1)-\phi(z_2)|^p\\
&\le C_1\frac{\Phi(\varepsilon)}{\Psi(\varepsilon)}\sum_{\substack{(z_1,z_2)\in E^{(\varepsilon)}\\z_1,z_2\in B^{(\varepsilon)}(x_0,AR)}}|\mathbf{f}(z_1)-\mathbf{f}(z_2)|^p+\frac{C_2}{\Psi(R)}\sum_{z\in B^{(\varepsilon)}(x_0,AR)}|\mathbf{f}(z)|^pm^{(\varepsilon)}(\{z\}).\label{eq_CS_Veps}\tag*{$\mathbf{CS}^{(\varepsilon)}(\Psi)$}
\end{align*}
\end{proposition}

\begin{proof}
Let $C_1,C_2,A_1$ be the constants from Proposition \ref{prop_CS_Xeps}. Consider $B^{(\varepsilon)}(x_0,R)$ also as a ball in $X^{(\varepsilon)}$, then there exists a cutoff function $\phi\in\mathsf{F}^{(\varepsilon)}\subseteq C(X^{(\varepsilon)})$ for $B^{(\varepsilon)}(x_0,4R)\subseteq B^{(\varepsilon)}(x_0,4A_1R)$ satisfying the condition stated therein. Since $\phi$ is continuous on $X^{(\varepsilon)}$, we consider its restriction to $V^{(\varepsilon)}$, still denoted by $\phi$.

For any $\mathbf{f}\in\mathbf{F}^{(\varepsilon)}$, without loss of generality, we may assume that $\mathbf{f}$ is non-negative. Let $\mathsf{f}\in l(X^{(\varepsilon)})$  be given by linear interpolation on each cable as in (\ref{eq_linear_inter}), then
$$\int_{]z_1,z_2[}\mathrm{d}\Gamma^{(\varepsilon)}(\mathsf{f})=\frac{\Phi(\varepsilon)}{\Psi(\varepsilon)}\lvert\mathbf{f}(z_1)-\mathbf{f}(z_2)\rvert^p\text{ for any }(z_1,z_2)\in E^{(\varepsilon)},$$
which gives $\mathsf{f}\in \mathsf{F}^{(\varepsilon)}$. By \ref{eq_CS_Xeps} on $X^{(\varepsilon)}$, we have
\begin{equation}\label{eq_CS_VepsA}
\int_{B^{(\varepsilon)}(x_0,4A_1R)}\mathsf{f}^p\mathrm{d}\Gamma^{(\varepsilon)}(\phi)\le C_{1}\int_{B^{(\varepsilon)}(x_0,4A_1R)}\mathrm{d}\Gamma^{(\varepsilon)}(\mathsf{f})+\frac{C_{2}}{\Psi(4R)}\int_{B^{(\varepsilon)}(x_0,4A_1R)}\mathsf{f}^p\mathrm{d} \lambda^{(\varepsilon)},
\end{equation}
where
\begin{align}
&\int_{B^{(\varepsilon)}(x_0,4A_1R)}\mathrm{d}\Gamma^{(\varepsilon)}(\mathsf{f})\le\frac{\Phi(\varepsilon)}{\Psi(\varepsilon)}\sum_{\substack{(z_1,z_2)\in E^{(\varepsilon)}\\z_1,z_2\in B^{(\varepsilon)}(x_0,8A_1R)}}|\mathbf{f}(z_1)-\mathbf{f}(z_2)|^p.\label{eq_CS_VepsB}
\end{align}
For any $(z_1,z_2)\in E^{(\varepsilon)}$ and any $z\in[z_1,z_2]$, we have $\mathsf{f}(z)^p\le\mathbf{f}(z_1)^p+\mathbf{f}(z_2)^p$, hence
\begin{align*}
&\int_{]z_1,z_2[}\mathsf{f}^p\mathrm{d}\lambda^{(\varepsilon)}\le\mathbf{f}(z_1)^pm^{(\varepsilon)}(\{z_1\})+\mathbf{f}(z_2)^pm^{(\varepsilon)}(\{z_2\}),
\end{align*}
which implies
\begin{align}
&\int_{B^{(\varepsilon)}(x_0,4A_1R)}\mathsf{f}^p\mathrm{d} \lambda^{(\varepsilon)}\lesssim\sum_{z\in B^{(\varepsilon)}(x_0,8A_1R)}\mathbf{f}(z)^pm^{(\varepsilon)}(\{z\})\label{eq_CS_VepsC}.
\end{align}
Moreover, if $\mathbf{f}(z_1)<\frac{1}{2}\mathbf{f}(z_2)$ or $\mathbf{f}(z_2)<\frac{1}{2}\mathbf{f}(z_1)$, then $\mathbf{f}(z_1)^p+\mathbf{f}(z_2)^p\le2^{p+1}\lvert\mathbf{f}(z_1)-\mathbf{f}(z_2)\rvert^p$. Otherwise, $\frac{1}{2}\mathbf{f}(z_2)\le\mathbf{f}(z_1)\le2\mathbf{f}(z_2)$, without loss of generality, we may assume that $\mathbf{f}(z_1)\le\mathbf{f}(z_2)$, then $\mathbf{f}(z_1)\le\mathsf{f}(z)\le\mathbf{f}(z_2)\le2\mathsf{f}(z)$ for any $z\in[z_1,z_2]$, hence
\begin{align*}
&\left(\mathbf{f}(z_1)^p+\mathbf{f}(z_2)^p\right)|\phi(z_1)-\phi(z_2)|^p\le2^{p+1}\frac{\Psi(\varepsilon)}{\Phi(\varepsilon)}\int_{]z_1,z_2[}\mathsf{f}^p\mathrm{d}\Gamma^{(\varepsilon)}(\phi).
\end{align*}
Therefore
\begin{align*}
&\sum_{z_1\in B^{(\varepsilon)}(x_0,8A_1R)}\mathbf{f}(z_1)^p \frac{\Phi(\varepsilon)}{\Psi(\varepsilon)}\sum_{z_2:(z_1,z_2)\in E^{(\varepsilon)}}|\phi(z_1)-\phi(z_2)|^p\\
&\lesssim\frac{\Phi(\varepsilon)}{\Psi(\varepsilon)}\sum_{\substack{(z_1,z_2)\in E^{(\varepsilon)}\\z_1\in B^{(\varepsilon)}(x_0,4A_1R)\text{ or }z_2\in B^{(\varepsilon)}(x_0,4A_1R)}}\left(\mathbf{f}(z_1)^p+\mathbf{f}(z_2)^p\right)|\phi(z_1)-\phi(z_2)|^p\\
&=\frac{\Phi(\varepsilon)}{\Psi(\varepsilon)}\left(\sum_{\mathbf{f}(z_1)<\frac{1}{2}\mathbf{f}(z_2)}+\sum_{\mathbf{f}(z_2)<\frac{1}{2}\mathbf{f}(z_1)}+\sum_{\frac{1}{2}\mathbf{f}(z_2)\le\mathbf{f}(z_1)\le2\mathbf{f}(z_2)}\right)\left(\mathbf{f}(z_1)^p+\mathbf{f}(z_2)^p\right)|\phi(z_1)-\phi(z_2)|^p\\
&\lesssim\frac{\Phi(\varepsilon)}{\Psi(\varepsilon)}\sum_{\substack{(z_1,z_2)\in E^{(\varepsilon)}\\z_1,z_2\in B^{(\varepsilon)}(x_0,4A_1R+\frac{5}{2}\varepsilon)}}|\mathbf{f}(z_1)-\mathbf{f}(z_2)|^p+\sum_{\substack{(z_1,z_2)\in E^{(\varepsilon)}\\z_1,z_2\in B^{(\varepsilon)}(x_0,4A_1R+\frac{5}{2}\varepsilon)}}\int_{]z_1,z_2[}\mathsf{f}^p\mathrm{d}\Gamma^{(\varepsilon)}(\phi)\\
&\le\frac{\Phi(\varepsilon)}{\Psi(\varepsilon)}\sum_{\substack{(z_1,z_2)\in E^{(\varepsilon)}\\z_1,z_2\in B^{(\varepsilon)}(x_0,8A_1R)}}|\mathbf{f}(z_1)-\mathbf{f}(z_2)|^p+\int_{B^{(\varepsilon)}(x_0,4A_1R)}\mathsf{f}^p\mathrm{d}\Gamma^{(\varepsilon)}(\phi)\\
&\overset{(*)}{\scalebox{2}[1]{$\le$}}\frac{\Phi(\varepsilon)}{\Psi(\varepsilon)}\sum_{\substack{(z_1,z_2)\in E^{(\varepsilon)}\\z_1,z_2\in B^{(\varepsilon)}(x_0,8A_1R)}}|\mathbf{f}(z_1)-\mathbf{f}(z_2)|^p\\
&\hspace{15pt}+C_{1}\int_{B^{(\varepsilon)}(x_0,4A_1R)}\mathrm{d}\Gamma^{(\varepsilon)}(\mathsf{f})+\frac{C_{2}}{\Psi(4R)}\int_{B^{(\varepsilon)}(x_0,4A_1R)}\mathsf{f}^p\mathrm{d} \lambda^{(\varepsilon)}\\
&\overset{(**)}{\scalebox{2}[1]{$\lesssim$}}\frac{\Phi(\varepsilon)}{\Psi(\varepsilon)}\sum_{\substack{(z_1,z_2)\in E^{(\varepsilon)}\\z_1,z_2\in B^{(\varepsilon)}(x_0,8A_1R)}}|\mathbf{f}(z_1)-\mathbf{f}(z_2)|^p+\frac{1}{\Psi(R)}\sum_{z\in B^{(\varepsilon)}(x_0,8A_1R)}\mathbf{f}(z)^pm^{(\varepsilon)}(\{z\}),
\end{align*}
where $(*)$ follows from (\ref{eq_CS_VepsA}) and $(**)$ follows from (\ref{eq_CS_VepsB}) and (\ref{eq_CS_VepsC}).
\end{proof}

We ``transfer" the cutoff Sobolev inequality from $V^{(\varepsilon)}$ to $X$ as follows, which yields Proposition \ref{prop_CS_X}.

\begin{proof}[Proof of Proposition \ref{prop_CS_X}]
Let $L$ be the constant in (\ref{eq_eps_ball}), and $A_1$ the constant appearing in Proposition \ref{prop_CS_Veps}. Fix a ball $B(x_0,R)$ in $X$, let $\varepsilon=\min\{\frac{1}{2}R,\frac{A_1}{4A_{PI}}R\}$. For any $k\ge0$, let $\varepsilon_k=\frac{1}{2^k}\varepsilon$ and $V^{(\varepsilon_k)}$ an $\varepsilon_k$-net containing $x_0$, let $\{\psi_z^{(\varepsilon_k)}:z\in V^{(\varepsilon_k)}\}$ be a partition of unity given by Lemma \ref{lem_partition}, $\phi_k\in\mathbf{F}^{(\varepsilon_k)}$ a cutoff function for $B^{(\varepsilon_k)}(x_0,2LR)\subseteq B^{(\varepsilon_k)}(x_0,2A_1LR)$ given by Proposition \ref{prop_CS_Veps}, and
$$\psi_k=\sum_{z\in V^{(\varepsilon_k)}}\phi_k(z)\psi_z^{(\varepsilon_k)}.$$
By (\ref{eq_eps_ball}), $\psi_k\in\mathcal{F}$ is a cutoff function for $B(x_0,R)\subseteq B(x_0,4A_1LR)$ for any $k$. By Lemma \ref{lem_Gammabfu}, $\{\psi_k\}_k$ is $\mathcal{E}$-bounded, since $\{\psi_k\}_k$ is obviously $L^p(X;m)$-bounded, $\{\psi_k\}_k$ is $\mathcal{E}_1$-bounded. By the Banach-Alaoglu theorem (see \cite[Theorem 3 in Chapter12]{Lax02}), $\{\psi_k\}_k$ is $\mathcal{E}_1$-weakly convergent to some function $\phi\in\mathcal{F}$. By Mazur's lemma (see \cite[Theorem 2 in Section V.1]{Yos95}), for any $k\ge0$, there exist $I_k\ge k$, $\lambda^{(k)}_i\ge0$ for $i=k,\ldots,I_k$ with $\sum_{i=k}^{I_k}\lambda^{(k)}_i=1$ such that $\{\sum_{i=k}^{I_k}\lambda^{(k)}_i\psi_i\}_k$ is $\mathcal{E}_1$-convergent to $\phi$.

It is obvious that $\phi\in\mathcal{F}$ is a cutoff function for $B(x_0,R)\subseteq B(x_0,4A_1LR)$.

Firstly, for any $f\in \mathcal{F}\cap C_c(X)$, we have
\begin{align*}
&\left(\int_{B(x_0,4A_1LR)}|f|^p\mathrm{d}\Gamma(\phi)\right)^{1/p}=\lim_{k\to+\infty}\left(\int_{B(x_0,4A_1LR)}|f|^p\mathrm{d}\Gamma \left(\sum_{i=k}^{I_k}\lambda^{(k)}_i\psi_i\right)\right)^{1/p}\\
&\le\varliminf_{k\to+\infty}\sum_{i=k}^{I_k}\lambda^{(k)}_i\left(\int_{B(x_0,4A_1LR)}|f|^p\mathrm{d}\Gamma \left(\psi_i\right)\right)^{1/p}\le\varlimsup_{k\to+\infty}\left(\int_{B(x_0,4A_1LR)}|f|^p\mathrm{d}\Gamma \left(\psi_k\right)\right)^{1/p}.
\end{align*}
For any $\delta>0$, by the uniform continuity of $f$, there exists $K\ge0$ such that for any $k>K$, for any $x,y\in X$ with $d(x,y)<\frac{5}{2}\varepsilon_k$, we have $|f(x)-f(y)|<\delta$, then for any $\overline{z}\in V^{(\varepsilon_k)}$, we have
\begin{equation*}
|f|\le|f_{B(\overline{z},\frac{5}{4}\varepsilon_k)}|+\delta\text{ in }B(\overline{z},\frac{5}{4}\varepsilon_k).
\end{equation*}
Thus we have
\begin{align*}
&\int_{B(x_0,4A_1LR)}|f|^p\mathrm{d}\Gamma \left(\psi_k\right)\\
&\le\sum_{\overline{z}\in V^{(\varepsilon_k)}}\int_{B(\overline{z},\frac{5}{4}\varepsilon_k)\cap B(x_0,4A_1LR)}|f|^p\mathrm{d}\Gamma \left(\sum_{z:(z,\overline{z})\in E^{(\varepsilon_k)}}(\phi_k(z)-\phi_k(\overline{z}))\psi_z^{(\varepsilon_k)}\right)\\
&\lesssim\sum_{\overline{z}\in B(x_0,4A_1LR+\frac{5}{4}\varepsilon_k)\cap V^{(\varepsilon_k)}}\left(|f_{B(\overline{z},\frac{5}{4}\varepsilon_k)}|+\delta\right)^p\sum_{z:(z,\overline{z})\in E^{(\varepsilon_k)}}|\phi_k(z)-\phi_k(\overline{z})|^p\mathcal{E} \left(\psi_z^{(\varepsilon_k)}\right)\\
&\lesssim\sum_{\overline{z}\in B^{(\varepsilon_k)}(x_0,8A_1L^2R)}\left(|f_{B(\overline{z},\frac{5}{4}\varepsilon_k)}|+\delta\right)^p\frac{\Phi(\varepsilon_k)}{\Psi(\varepsilon_k)}\sum_{z:(z,\overline{z})\in E^{(\varepsilon_k)}}|\phi_k(z)-\phi_k(\overline{z})|^p\\
&=\sum_{\overline{z}\in B^{(\varepsilon_k)}(x_0,2A_1LR)}\left(|f_{B(\overline{z},\frac{5}{4}\varepsilon_k)}|+\delta\right)^p\frac{\Phi(\varepsilon_k)}{\Psi(\varepsilon_k)}\sum_{z:(z,\overline{z})\in E^{(\varepsilon_k)}}|\phi_k(z)-\phi_k(\overline{z})|^p\\
&\lesssim\frac{\Phi(\varepsilon_k)}{\Psi(\varepsilon_k)}\sum_{\substack{(z,\overline{z})\in E^{(\varepsilon_k)}\\z,\overline{z}\in B^{(\varepsilon_k)}(x_0,2A_1LR)}}\left(|f_{B(z,\frac{5}{4}\varepsilon_k)}|-|f_{B(\overline{z},\frac{5}{4}\varepsilon_k)}|\right)^p\\
&\hspace{15pt}+\frac{1}{\Psi(R)}\sum_{\overline{z}\in B^{(\varepsilon_k)}(x_0,2A_1LR)}\left(|f_{B(\overline{z},\frac{5}{4}\varepsilon_k)}|+\delta\right)^pm^{(\varepsilon_k)}(\{\overline{z}\}),
\end{align*}
where the last inequality follows from \hyperref[eq_CS_Veps]{$\mathbf{CS}^{(\varepsilon_k)}(\Psi)$} on $V^{(\varepsilon_k)}$. By \ref{eq_PI} and \ref{eq_VPhi}, we have
\begin{align*}
&\sum_{\substack{(z,\overline{z})\in E^{(\varepsilon_k)}\\z,\overline{z}\in B^{(\varepsilon_k)}(x_0,2A_1LR)}}\left(|f_{B(z,\frac{5}{4}\varepsilon_k)}|-|f_{B(\overline{z},\frac{5}{4}\varepsilon_k)}|\right)^p\le\sum_{\substack{(z,\overline{z})\in E^{(\varepsilon_k)}\\z,\overline{z}\in B^{(\varepsilon_k)}(x_0,2A_1LR)}}\lvert f_{B(z,\frac{5}{4}\varepsilon_k)}-f_{B(\overline{z},\frac{5}{4}\varepsilon_k)}\rvert^p\\
&\lesssim\sum_{z\in B^{(\varepsilon_k)}(x_0,2A_1LR)}\frac{\Psi(\varepsilon_k)}{\Phi(\varepsilon_k)}\int_{B(z,\frac{15}{4}A_{PI}\varepsilon_k)}\mathrm{d}\Gamma(f)\lesssim \frac{\Psi(\varepsilon_k)}{\Phi(\varepsilon_k)}\int_{B(x_0,4A_1LR)}\mathrm{d}\Gamma(f).
\end{align*}
By H\"older's inequality and \ref{eq_VPhi}, we have
\begin{align*}
&\sum_{\overline{z}\in B^{(\varepsilon_k)}(x_0,2A_1LR)}\left(|f_{B(\overline{z},\frac{5}{4}\varepsilon_k)}|+\delta\right)^pm^{(\varepsilon_k)}(\{\overline{z}\})\\
&\lesssim\sum_{\overline{z}\in B^{(\varepsilon_k)}(x_0,2A_1LR)}\int_{B(\overline{z},\frac{5}{4}\varepsilon_k)}\left(|f|^p+\delta^p\right)\mathrm{d} m\lesssim\int_{B(x_0,4A_1LR)}|f|^p\mathrm{d} m+\delta^pV(x_0,4A_1LR).
\end{align*}
Hence
\begin{align*}
&\int_{B(x_0,4A_1LR)}|f|^p\mathrm{d}\Gamma(\phi)\le\varlimsup_{K<k\to+\infty}\int_{B(x_0,4A_1LR)}|f|^p\mathrm{d}\Gamma \left(\psi_k\right)\\
&\lesssim\int_{B(x_0,4A_1LR)}\mathrm{d}\Gamma(f)+\frac{1}{\Psi(R)}\left(\int_{B(x_0,4A_1LR)}|f|^p\mathrm{d} m+\delta^pV(x_0,4A_1LR)\right).
\end{align*}
Since $\delta>0$ is arbitrary, we have
$$\int_{B(x_0,4A_1LR)}|f|^p\mathrm{d}\Gamma(\phi)\lesssim\int_{B(x_0,4A_1LR)}\mathrm{d}\Gamma(f)+\frac{1}{\Psi(R)}\int_{B(x_0,4A_1LR)}|f|^p\mathrm{d} m.$$

Secondly, for any $f\in\mathcal{F}$, by \cite[Proposition 8.5, Proposition 8.12]{Yan25a}, there exists $\{f_n\}\subseteq\mathcal{F}\cap C_c(X)$ such that $\{f_n\}$ is $\mathcal{E}_1$-convergent to $f$ and $\{f_n\}$ converges to $\widetilde{f}$ q.e. on $X$, which is also $\Gamma(\phi)$-a.e. on $X$, then by Fatou's lemma, we have
\begin{align*}
&\int_{B(x_0,4A_1LR)}|\widetilde{f}|^p\mathrm{d}\Gamma(\phi)\le\varliminf_{n\to+\infty}\int_{B(x_0,4A_1LR)}|f_n|^p\mathrm{d}\Gamma(\phi)\\
&\lesssim\varliminf_{n\to+\infty}\left(\int_{B(x_0,4A_1LR)}\mathrm{d}\Gamma(f_n)+\frac{1}{\Psi(R)}\int_{B(x_0,4A_1LR)}|f_n|^p\mathrm{d} m\right)\\
&=\int_{B(x_0,4A_1LR)}\mathrm{d}\Gamma(f)+\frac{1}{\Psi(R)}\int_{B(x_0,4A_1LR)}|f|^p\mathrm{d} m.
\end{align*}
\end{proof}

\section{Proof of Proposition \ref{prop_RBA}}\label{sec_RBA}

The proof follows the same technique as that of the classical result stating that, under suitable functional inequalities, general bounded open sets can be approximated by bounded open sets satisfying the so-called corkscrew condition. The corkscrew condition, in turn, guarantees the validity of a certain Wiener test, which implies that the set is regular.

For simplicity, by appropriately adjusting the constants in the following results, we always assume that the annuli appearing therein satisfy the assumption as in the \ref{eq_EHIann} condition.

\begin{lemma}\label{lem_Wiener1}
Assume the setting of Proposition \ref{prop_RBA}. Then there exists $C>0$ such that for any ball $B$ and any compact set $K\subseteq B$, letting $u\in \mathcal{F}(4B)$ be the capacitory potential of $(K,X\backslash 4B)$ as given in Lemma \ref{lem_cap_dual} (1), we have
$$\einf_{2B}u\ge {C}\left(\frac{\mathrm{cap}\left(K,X\backslash 4B\right)}{\mathrm{cap}\left(B,X\backslash 4B\right)}\right)^{1/(p-1)}.$$
\end{lemma}

\begin{proof}
The proof follows essentially the same idea as in \cite[Lemma 11.20]{BB11} and in \cite[Lemma 6.21]{HKM06}. By Lemma \ref{lem_cap_dual} (1), we have $u\in \mathcal{F}(4B)$ is superharmonic in $4B$, harmonic in $4B\backslash K$, satisfies $0\le u\le 1$ in $X$, and $\widetilde{u}=1$ q.e. on $K$. By \hyperlink{eq_LSC}{\text{LSC}}, we may assume that $u$ is lower semi-continuous. Since $u\in \mathcal{F}(4B)$ is superharmonic in $4B$, by Proposition \ref{prop_comparison}, we have
$$\einf_{2B}u\ge\einf_{3B\backslash {2B}}u.$$
Since $u\in \mathcal{F}(4B)$ is non-negative superharmonic in $4B$ and harmonic in $4B\backslash \overline{B}$, by \ref{eq_EHIann}, there exists $C_1>0$ such that
\begin{equation}\label{eq_RBA1}
\esup_{4B\backslash {2B}}u=\esup_{3B\backslash {2B}}u\le C_1\einf_{3B\backslash {2B}}u\le C_1\einf_{2B}u.
\end{equation}

Let $M=\esup_{4B\backslash {2B}}u\in[0,1]$. Firstly, if $M<1$, then for any $\varepsilon\in(0,1-M)$, by the lower semi-continuity of $u$, we have $\{u>M+\varepsilon\}$ is an open set containing $K$, and by \ref{eq_EHI}, we have $u\in C(4B\backslash K)$, hence $\overline{\{u>M+\varepsilon\}}=\{u\ge M+\varepsilon\}\subseteq 2B$ is a compact set. Let
\begin{align*}
u_1&=\frac{u}{M+\varepsilon}\wedge 1=\frac{u}{M+\varepsilon}1_{\{u< M+\varepsilon\}}+1_{\{u\ge M+\varepsilon\}},\\
u_2&=\frac{u-(M+\varepsilon)u_1}{1-(M+\varepsilon)}=\frac{u-(M+\varepsilon)}{1-(M+\varepsilon)}1_{\{u>M+\varepsilon\}},
\end{align*}
then $u_1,u_2\in \mathcal{F}(4B)$, and $\widetilde{u}_1=\widetilde{u}_2=1$ q.e. on $K$. Hence for any $a\in[0,1]$, we have $au_1+(1-a)u_2\in \mathcal{F}(4B)$, $a \widetilde{u}_1+(1-a)\widetilde{u}_2=1$ q.e. on $K$, and
$$\Gamma \left(au_1+(1-a)u_2\right)=a^p1_{\{u< M+\varepsilon\}}\Gamma(u_1)+(1-a)^p1_{\{u>M+\varepsilon\}}\Gamma(u_2),$$
which gives
\begin{align*}
&\mathrm{cap}\left(K,X\backslash 4B\right)\le \mathcal{E}(au_1+(1-a)u_2)=\int_{4B}\mathrm{d}\Gamma\left(au_1+(1-a)u_2\right)\\
&=a^p\int_{4B\cap\{u< M+\varepsilon\}}\mathrm{d}\Gamma(u_1)+(1-a)^p\int_{4B\cap\{u>M+\varepsilon\}}\mathrm{d}\Gamma(u_2)\le a^p \mathcal{E}(u_1)+(1-a)^p \mathcal{E}(u_2).
\end{align*}
Since $a\in[0,1]$ is arbitrary, we have
$$\mathcal{E}(u)=\mathrm{cap}\left(K,X\backslash 4B\right)\le \inf_{a\in[0,1]}\left(a^p \mathcal{E}(u_1)+(1-a)^p \mathcal{E}(u_2)\right)=\frac{1}{\left(\frac{1}{\mathcal{E}(u_1)^{\frac{1}{p-1}}}+\frac{1}{\mathcal{E}(u_2)^{\frac{1}{p-1}}}\right)^{p-1}},$$
where
$$\mathcal{E}(u_2)=\int_{4B\cap \{u>M+\varepsilon\}}\frac{1}{(1-(M+\varepsilon))^p}\mathrm{d}\Gamma(u)\le \frac{1}{(1-(M+\varepsilon))^p}\mathcal{E}(u),$$
which gives
$$\mathcal{E}(u)\le\frac{1}{\left(\frac{1}{\mathcal{E}(u_1)^{\frac{1}{p-1}}}+\frac{(1-(M+\varepsilon))^{\frac{p}{p-1}}}{\mathcal{E}(u)^{\frac{1}{p-1}}}\right)^{p-1}},$$
that is,
$$\mathcal{E}(u)\le \left(1-\left(1-(M+\varepsilon)\right)^{\frac{p}{p-1}}\right)^{p-1}\mathcal{E}(u_1)\le \left(\frac{p}{p-1}(M+\varepsilon)\right)^{p-1}\mathcal{E}(u_1),$$
where in the second inequality, we use the fact that $(1-x)^\alpha\ge 1-\alpha x$ for any $x\in[0,1]$ and any $\alpha\ge1$. Since $u_1\in \mathcal{F}(4B)$ satisfies $\widetilde{u}_1=1$ q.e. on ${\{u\ge M+\varepsilon\}}$ and $u_1$ is harmonic in $4B\backslash {\{u\ge M+\varepsilon\}}$, we have
$$\mathcal{E}(u_1)=\mathrm{cap}\left({\{u\ge M+\varepsilon\}},X\backslash 4B\right)\le \mathrm{cap}\left(2B,X\backslash 4B\right)\le C_2 \mathrm{cap}\left(B,X\backslash 4B\right),$$
where in the last inequality, we use \ref{eq_VD}, \hyperlink{eq_cap}{$\text{cap}(\Psi)$}, and $C_2>0$ depends only on $C_\Psi, C_{VD}, A_{cap}$, $C_{cap}$. Hence
$$\mathrm{cap}\left(K,X\backslash 4B\right)=\mathcal{E}(u)\le \left(\frac{p}{p-1}(M+\varepsilon)\right)^{p-1} C_2 \mathrm{cap}\left(B,X\backslash 4B\right),$$
that is,
$$M+\varepsilon\ge \frac{p-1}{p}\frac{1}{C_2^{1/(p-1)}}\left(\frac{\mathrm{cap}\left(K,X\backslash 4B\right)}{\mathrm{cap}\left(B,X\backslash 4B\right)}\right)^{1/(p-1)}.$$
Since $\varepsilon\in(0,1-M)$ is arbitrary, we have
$$M\ge \frac{p-1}{p}\frac{1}{C_2^{1/(p-1)}}\left(\frac{\mathrm{cap}\left(K,X\backslash 4B\right)}{\mathrm{cap}\left(B,X\backslash 4B\right)}\right)^{1/(p-1)}.$$
Secondly, if $M=1$, then obviously
$$M=1\ge \left(\frac{\mathrm{cap}\left(K,X\backslash 4B\right)}{\mathrm{cap}\left(B,X\backslash 4B\right)}\right)^{1/(p-1)}.$$
Therefore, combining with (\ref{eq_RBA1}), we have
$$\einf_{2B}u\ge \frac{1}{C_1}M\ge \frac{1}{C_1}\left(\frac{p-1}{p}\frac{1}{C_2^{1/(p-1)}}\wedge 1\right)\left(\frac{\mathrm{cap}\left(K,X\backslash 4B\right)}{\mathrm{cap}\left(B,X\backslash 4B\right)}\right)^{1/(p-1)}.$$
\end{proof}

\begin{lemma}\label{lem_Wiener2}
Assume the setting of Proposition \ref{prop_RBA}. Let $C$ be the constant given in Lemma \ref{lem_Wiener1}. Then the following holds. For any ball $B$ and any compact set $K\subseteq B$, let $u\in \mathcal{F}(4B)$ be the capacitory potential of $(K,X\backslash 4B)$ as given in Lemma \ref{lem_cap_dual} (1). For any $n\ge-1$, let $B_n=\frac{1}{4^n}B$, then for any $n\ge0$, we have
$$\esup_{2B_n}\left(1-u\right)\le \exp \left(-C\sum_{k=0}^n \left(\frac{\mathrm{cap}\left(\overline{B}_k\cap K,X\backslash B_{k-1}\right)}{\mathrm{cap}(B_k,X\backslash B_{k-1})}\right)^{1/(p-1)}\right).$$
\end{lemma}

The proof follows essentially the same idea as in \cite[Theorem 11.21]{BB11} and in \cite[Lemma 6.25]{HKM06}, using the estimate in Lemma \ref{lem_Wiener1} together with an inductive argument, and is therefore omitted.

\begin{proposition}\label{prop_Wiener}
Assume the setting of Proposition \ref{prop_RBA}. Let $U$ be a bounded open subset of $X$ and $x_0\in \partial U$. Assume that there exists $r_0>0$ such that
\begin{equation}\label{eq_Wiener}
\sum_{n=0}^\infty \left(\frac{\mathrm{cap}\left(B(x_0,\frac{1}{4^{n+1}}r_0)\backslash U,X\backslash B(x_0,\frac{1}{4^{n}}r_0)\right)}{\mathrm{cap}\left(B(x_0,\frac{1}{4^{n+1}}r_0),X\backslash B(x_0,\frac{1}{4^{n}}r_0)\right)}\right)^{1/(p-1)}=+\infty.
\end{equation}
Then for any $h\in \mathcal{F}(U)$ which is bounded and harmonic in $B(x_0,2r_0)\cap U$, we have
$$\lim_{U\ni y\to x_0}h(y)=0.$$
In particular, if $X\backslash U$ has a corkscrew at $x_0$, that is, there exist $c, r_0>0$ such that for any $r\in(0,r_0)$, $B(x_0,r)\backslash U$ contains a ball with radius $cr$, then (\ref{eq_Wiener}) holds.
\end{proposition}

\begin{proof}
Let $u\in \mathcal{F}(B(x_0,r_0))$ be the capacitory potential of $(\overline{B(x_0,\frac{1}{4}r_0)}\backslash U,X\backslash B(x_0,r_0))$ as given in Lemma \ref{lem_cap_dual} (1). For notational convenience, we write $1-u$ to denote $\psi-u$, where $\psi\in \mathcal{F}\cap C_c(X)$ satisfies $\psi=1$ on $\overline{B(x_0,2r_0)}$. Since $h$ is bounded in $B(x_0,2r_0)\cap U$, let $M\in(0,+\infty)$ satisfy $\lvert h\rvert\le M$ in $B(x_0,2r_0)\cap U$, then $\widetilde{h}\le M=M(1-\widetilde{u})$ q.e. on $(B(x_0,2r_0)\cap U)\backslash B(x_0,r_0)$. Since $h\in \mathcal{F}(U)$, we have $\widetilde{h}=0\le M(1-\widetilde{u})$ q.e. on $B(x_0,2r_0)\backslash U$. Since $h, 1-u$ are harmonic in $B(x_0,r_0)\cap U$, by Proposition \ref{prop_comparison}, we have $h\le M(1-u)$ in $B(x_0,r_0)\cap U$. Similarly, we have $-h\le M(1-u)$ in $B(x_0,r_0)\cap U$. Hence by Lemma \ref{lem_Wiener2}, for any $n\ge0$, we have
\begin{align*}
&\esup_{B(x_0,\frac{2}{4^{n+1}}r_0)\cap U}\lvert h\rvert\le M\esup_{B(x_0,\frac{2}{4^{n+1}}r_0)}(1-u)\\
&\le M\exp \left(-C\sum_{k=0}^n\left(\frac{\mathrm{cap}\left(\overline{B(x_0,\frac{1}{4^{k+1}}r_0)}\backslash U,X\backslash B(x_0,\frac{1}{4^{k}}r_0)\right)}{\mathrm{cap}\left(B(x_0,\frac{1}{4^{k+1}}r_0),X\backslash B(x_0,\frac{1}{4^{k}}r_0)\right)}\right)^{1/(p-1)}\right)\\
&\le M\exp \left(-C\sum_{k=0}^n\left(\frac{\mathrm{cap}\left(B(x_0,\frac{1}{4^{k+1}}r_0)\backslash U,X\backslash B(x_0,\frac{1}{4^{k}}r_0)\right)}{\mathrm{cap}\left(B(x_0,\frac{1}{4^{k+1}}r_0),X\backslash B(x_0,\frac{1}{4^{k}}r_0)\right)}\right)^{1/(p-1)}\right)\to0
\end{align*}
as $n\to+\infty$, which gives $\lim_{U\ni y\to x_0}h(y)=0$.

If $X\backslash U$ has a corkscrew at $x_0$, then for any $n\ge0$, there exists $B(y_n,\frac{c}{4^{n+1}}r_0)\subseteq B(x_0,\frac{1}{4^{n+1}}r_0)\backslash U$, by \ref{eq_VD}, \hyperlink{eq_cap}{$\text{cap}(\Psi)$}, we have
$$\frac{\mathrm{cap}\left(B(x_0,\frac{1}{4^{n+1}}r_0)\backslash U,X\backslash B(x_0,\frac{1}{4^{n}}r_0)\right)}{\mathrm{cap}\left(B(x_0,\frac{1}{4^{n+1}}r_0),X\backslash B(x_0,\frac{1}{4^{n}}r_0)\right)}\ge\frac{\mathrm{cap}\left(B(y_n,\frac{c}{4^{n+1}}r_0),X\backslash B(x_0,\frac{1}{4^{n}}r_0)\right)}{\mathrm{cap}\left(B(x_0,\frac{1}{4^{n+1}}r_0),X\backslash B(x_0,\frac{1}{4^{n}}r_0)\right)}\asymp1,$$
hence (\ref{eq_Wiener}) holds.
\end{proof}

We give the proof of Proposition \ref{prop_RBA} as follows.

\begin{proof}[Proof of Proposition \ref{prop_RBA}]
Let $B=B(z_0,R_0)$ be a ball and $\delta\in(0,1)$. By Proposition \ref{prop_Wiener}, we only need to show that there exists a bounded open set $U$ with $\delta B\subseteq U\subseteq B$ such that for any $x_0\in \partial U$, $X\backslash U$ has a corkscrew at $x_0$. Then for any $u\in \mathcal{F}(U)$ which is bounded and harmonic in $U\backslash K$ for some compact set $K\subseteq U$, there exists $r_0\in(0,+\infty)$ such that $u$ is bounded and harmonic in $B(x_0,2r_0)\cap U$ and (\ref{eq_Wiener}) holds, hence $\lim_{U\ni y\to x_0}u(y)=0$. By \ref{eq_EHI}, we have $u\in C(U\backslash K)$, hence $u\in C(\overline{U}\backslash K)$, which gives $\lim_{U\ni y\to x}u(y)=0$ uniformly for all $x\in \partial U$. Therefore, $U$ is regular, \hyperlink{eq_RBA}{\text{RBA}} holds.

By \ref{eq_CC}, there exist a geodesic metric $\rho$ on $X$ and $L\ge1$ such that
$$d(x,y)\le \rho(x,y)\le Ld(x,y)\text{ for any }x,y\in X,$$
see the proof of \cite[PROPOSITION A.1 (a)$\Rightarrow$(b)]{KM20}, then the remainder of the proof follows the same argument as in \cite[Lemma 14.2]{BB11}, with the inner metric $d_{\mathrm{in}}$ therein replaced by $\rho$, and is therefore omitted.
\end{proof}

\bibliographystyle{plain}

\begin{thebibliography}{10}

\bibitem{ABKY11}
Daniel Aalto, Lauri Berkovits, Outi~Elina Kansanen, and Hong Yue.
\newblock John-{N}irenberg lemmas for a doubling measure.
\newblock {\em Studia Math.}, 204(1):21--37, 2011.

\bibitem{AB25}
Patricia Alonso-Ruiz and Fabrice Baudoin.
\newblock Korevaar-{S}choen {$p$}-energies and their {$\Gamma $}-limits on
  {C}heeger spaces.
\newblock {\em Nonlinear Anal.}, 256:Paper No. 113779, 22, 2025.

\bibitem{AB15}
Sebastian Andres and Martin~T. Barlow.
\newblock Energy inequalities for cutoff functions and some applications.
\newblock {\em J. Reine Angew. Math.}, 699:183--215, 2015.

\bibitem{AES25a}
Riku {Anttila}, Sylvester {Eriksson-Bique}, and Ryosuke {Shimizu}.
\newblock {Construction of self-similar energy forms and singularity of Sobolev
  spaces on Laakso-type fractal spaces}.
\newblock {\em arXiv e-prints}, page arXiv:2503.13258, March 2025.

\bibitem{BCLS95}
Dominique Bakry, Thierry Coulhon, Michel Ledoux, and Laurent Saloff-Coste.
\newblock Sobolev inequalities in disguise.
\newblock {\em Indiana Univ. Math. J.}, 44(4):1033--1074, 1995.

\bibitem{Bar13}
Martin~T. Barlow.
\newblock Analysis on the {S}ierpinski carpet.
\newblock In {\em Analysis and geometry of metric measure spaces}, volume~56 of
  {\em CRM Proc. Lecture Notes}, pages 27--53. Amer. Math. Soc., Providence,
  RI, 2013.

\bibitem{BB04}
Martin~T. Barlow and Richard~F. Bass.
\newblock Stability of parabolic {H}arnack inequalities.
\newblock {\em Trans. Amer. Math. Soc.}, 356(4):1501--1533, 2004.

\bibitem{BBK06}
Martin~T. Barlow, Richard~F. Bass, and Takashi Kumagai.
\newblock Stability of parabolic {H}arnack inequalities on metric measure
  spaces.
\newblock {\em J. Math. Soc. Japan}, 58(2):485--519, 2006.

\bibitem{BC23}
Fabrice Baudoin and Li~Chen.
\newblock Sobolev spaces and {P}oincar\'e{} inequalities on the {V}icsek
  fractal.
\newblock {\em Ann. Fenn. Math.}, 48(1):3--26, 2023.

\bibitem{BM95AMPA}
Marco Biroli and Umberto Mosco.
\newblock A {S}aint-{V}enant type principle for {D}irichlet forms on
  discontinuous media.
\newblock {\em Ann. Mat. Pura Appl. (4)}, 169:125--181, 1995.

\bibitem{BM95}
Marco Biroli and Umberto Mosco.
\newblock Sobolev inequalities on homogeneous spaces.
\newblock {\em Potential Anal.}, 4(4):311--324, 1995.
\newblock Potential theory and degenerate partial differential operators
  (Parma).

\bibitem{BB11}
Anders Bj\"orn and Jana Bj\"orn.
\newblock {\em Nonlinear potential theory on metric spaces}, volume~17 of {\em
  EMS Tracts in Mathematics}.
\newblock European Mathematical Society (EMS), Z\"urich, 2011.

\bibitem{Bla01}
S\'ebastien Blach\`ere.
\newblock Harmonic functions on annuli of graphs.
\newblock {\em Ann. Math. Blaise Pascal}, 8(2):47--59, 2001.

\bibitem{BBI01}
Dmitri Burago, Yuri Burago, and Sergei Ivanov.
\newblock {\em A course in metric geometry}, volume~33 of {\em Graduate Studies
  in Mathematics}.
\newblock American Mathematical Society, Providence, RI, 2001.

\bibitem{CGQ22}
Shiping Cao, Qingsong Gu, and Hua Qiu.
\newblock {$p$}-energies on p.c.f. self-similar sets.
\newblock {\em Adv. Math.}, 405:Paper No. 108517, 58, 2022.

\bibitem{CGYZ26}
Aobo Chen, Jin Gao, Zhenyu Yu, and Junda Zhang.
\newblock Besov--{L}ipschitz norm and {$p$}-energy measure on scale-irregular
  {V}icsek sets.
\newblock {\em J. Fractal Geom.}, 13(1-2):37--85, 2026.

\bibitem{CKW21}
Zhen-Qing Chen, Takashi Kumagai, and Jian Wang.
\newblock Stability of heat kernel estimates for symmetric non-local
  {D}irichlet forms.
\newblock {\em Mem. Amer. Math. Soc.}, 271(1330):v+89, 2021.

\bibitem{CG98}
Thierry Coulhon and Alexander Grigor'yan.
\newblock Random walks on graphs with regular volume growth.
\newblock {\em Geom. Funct. Anal.}, 8(4):656--701, 1998.

\bibitem{CHS01}
Thierry Coulhon, Ilkka Holopainen, and Laurent Saloff-Coste.
\newblock Harnack inequality and hyperbolicity for subelliptic
  {$p$}-{L}aplacians with applications to {P}icard type theorems.
\newblock {\em Geom. Funct. Anal.}, 11(6):1139--1191, 2001.

\bibitem{DRY23}
Baptiste Devyver, Emmanuel Russ, and Meng Yang.
\newblock {Gradient Estimate for the Heat Kernel on Some Fractal-Like Cable
  Systems and Quasi-Riesz Transforms}.
\newblock {\em International Mathematics Research Notices},
  2023(18):15537--15583, September 2023.

\bibitem{Eri26}
Sylvester {Eriksson-Bique}.
\newblock {On the Resistance Conjecture}.
\newblock {\em arXiv e-prints}, page arXiv:2602.05477v2, February 2026.

\bibitem{FOT11}
Masatoshi Fukushima, Yoichi Oshima, and Masayoshi Takeda.
\newblock {\em Dirichlet forms and symmetric {M}arkov processes}, volume~19 of
  {\em De Gruyter Studies in Mathematics}.
\newblock Walter de Gruyter \& Co., Berlin, extended edition, 2011.

\bibitem{Gri92}
Alexander Grigor'yan.
\newblock The heat equation on noncompact {R}iemannian manifolds.
\newblock {\em Mat. Sb.}, 182(1):55--87, 1991.

\bibitem{Gri22}
Alexander Grigor'yan.
\newblock Analysis on fractal spaces and heat kernels.
\newblock In {\em Dirichlet forms and related topics}, volume 394 of {\em
  Springer Proc. Math. Stat.}, pages 143--159. Springer, Singapore, [2022]
  \copyright 2022.

\bibitem{GHH18}
Alexander Grigor'yan, Eryan Hu, and Jiaxin Hu.
\newblock Two-sided estimates of heat kernels of jump type {D}irichlet forms.
\newblock {\em Adv. Math.}, 330:433--515, 2018.

\bibitem{GH14}
Alexander Grigor'yan and Jiaxin Hu.
\newblock Heat kernels and {G}reen functions on metric measure spaces.
\newblock {\em Canad. J. Math.}, 66(3):641--699, 2014.

\bibitem{GHL14}
Alexander Grigor'yan, Jiaxin Hu, and Ka-Sing Lau.
\newblock Heat kernels on metric measure spaces.
\newblock In {\em Geometry and analysis of fractals}, volume~88 of {\em
  Springer Proc. Math. Stat.}, pages 147--207. Springer, Heidelberg, 2014.

\bibitem{GHL15}
Alexander Grigor'yan, Jiaxin Hu, and Ka-Sing Lau.
\newblock Generalized capacity, {H}arnack inequality and heat kernels of
  {D}irichlet forms on metric measure spaces.
\newblock {\em J. Math. Soc. Japan}, 67(4):1485--1549, 2015.

\bibitem{GS05}
Alexander Grigor'yan and Laurent Saloff-Coste.
\newblock Stability results for {H}arnack inequalities.
\newblock {\em Ann. Inst. Fourier (Grenoble)}, 55(3):825--890, 2005.

\bibitem{GT12}
Alexander Grigor'yan and Andras Telcs.
\newblock Two-sided estimates of heat kernels on metric measure spaces.
\newblock {\em Ann. Probab.}, 40(3):1212--1284, 2012.

\bibitem{HK95}
Piotr Haj{\l}asz and Pekka Koskela.
\newblock Sobolev meets {P}oincar\'e.
\newblock {\em C. R. Acad. Sci. Paris S\'er. I Math.}, 320(10):1211--1215,
  1995.

\bibitem{Hei01}
Juha Heinonen.
\newblock {\em Lectures on analysis on metric spaces}.
\newblock Universitext. Springer-Verlag, New York, 2001.

\bibitem{HKM06}
Juha Heinonen, Tero Kilpel\"ainen, and Olli Martio.
\newblock {\em Nonlinear potential theory of degenerate elliptic equations}.
\newblock Dover Publications, Inc., Mineola, NY, 2006.
\newblock Unabridged republication of the 1993 original.

\bibitem{HK98}
Juha Heinonen and Pekka Koskela.
\newblock Quasiconformal maps in metric spaces with controlled geometry.
\newblock {\em Acta Math.}, 181(1):1--61, 1998.

\bibitem{HKST15}
Juha Heinonen, Pekka Koskela, Nageswari Shanmugalingam, and Jeremy~T. Tyson.
\newblock {\em Sobolev spaces on metric measure spaces}, volume~27 of {\em New
  Mathematical Monographs}.
\newblock Cambridge University Press, Cambridge, 2015.
\newblock An approach based on upper gradients.

\bibitem{HPS04}
Paul~Edward Herman, Roberto Peirone, and Robert~S. Strichartz.
\newblock {$p$}-energy and {$p$}-harmonic functions on {S}ierpinski gasket type
  fractals.
\newblock {\em Potential Anal.}, 20(2):125--148, 2004.

\bibitem{Hol03}
Ilkka Holopainen.
\newblock Quasiregular mappings and the {$p$}-{L}aplace operator.
\newblock In {\em Heat kernels and analysis on manifolds, graphs, and metric
  spaces ({P}aris, 2002)}, volume 338 of {\em Contemp. Math.}, pages 219--239.
  Amer. Math. Soc., Providence, RI, 2003.

\bibitem{HY23}
Jiaxin Hu and Zhenyu Yu.
\newblock The weak elliptic {H}arnack inequality revisited.
\newblock {\em Asian J. Math.}, 27(5):771--828, 2023.

\bibitem{KM20}
Naotaka Kajino and Mathav Murugan.
\newblock On singularity of energy measures for symmetric diffusions with full
  off-diagonal heat kernel estimates.
\newblock {\em Ann. Probab.}, 48(6):2920--2951, 2020.

\bibitem{KS24a}
Naotaka {Kajino} and Ryosuke {Shimizu}.
\newblock {Contraction properties and differentiability of $p$-energy forms
  with applications to nonlinear potential theory on self-similar sets}.
\newblock {\em arXiv e-prints}, page arXiv:2404.13668v3, April 2024.

\bibitem{Kig23}
Jun Kigami.
\newblock {\em Conductive homogeneity of compact metric spaces and construction
  of {$p$}-energy}, volume~5 of {\em Memoirs of the European Mathematical
  Society}.
\newblock European Mathematical Society (EMS), Berlin, [2023] \copyright 2023.

\bibitem{KM92Pisa}
Tero Kilpel\"ainen and Jan Mal\'y.
\newblock Degenerate elliptic equations with measure data and nonlinear
  potentials.
\newblock {\em Ann. Scuola Norm. Sup. Pisa Cl. Sci. (4)}, 19(4):591--613, 1992.

\bibitem{KM94Acta}
Tero Kilpel\"ainen and Jan Mal\'y.
\newblock The {W}iener test and potential estimates for quasilinear elliptic
  equations.
\newblock {\em Acta Math.}, 172(1):137--161, 1994.

\bibitem{KM03}
Juha Kinnunen and Olli Martio.
\newblock Sobolev space properties of superharmonic functions on metric spaces.
\newblock {\em Results Math.}, 44(1-2):114--129, 2003.

\bibitem{Kor07}
Riikka Korte.
\newblock Geometric implications of the {P}oincar\'e{} inequality.
\newblock {\em Results Math.}, 50(1-2):93--107, 2007.

\bibitem{Lax02}
Peter~D. Lax.
\newblock {\em Functional analysis}.
\newblock Pure and Applied Mathematics (New York). Wiley-Interscience [John
  Wiley \& Sons], New York, 2002.

\bibitem{MS95}
Patrick Maheux and Laurent Saloff-Coste.
\newblock Analyse sur les boules d'un op\'erateur sous-elliptique.
\newblock {\em Math. Ann.}, 303(4):713--740, 1995.

\bibitem{Mur19}
Mathav Murugan.
\newblock Quasisymmetric uniformization and heat kernel estimates.
\newblock {\em Trans. Amer. Math. Soc.}, 372(6):4177--4209, 2019.

\bibitem{Mur23a}
Mathav Murugan.
\newblock A note on heat kernel estimates, resistance bounds and {P}oincar\'e{}
  inequality.
\newblock {\em Asian J. Math.}, 27(6):853--866, 2023.

\bibitem{MS25}
Mathav Murugan and Ryosuke Shimizu.
\newblock First-order {S}obolev spaces, self-similar energies and energy
  measures on the {S}ierpi\'nski carpet.
\newblock {\em Comm. Pure Appl. Math.}, 78(9):1523--1608, 2025.

\bibitem{Sal92}
Laurent Saloff-Coste.
\newblock A note on {P}oincar\'{e}, {S}obolev, and {H}arnack inequalities.
\newblock {\em Internat. Math. Res. Notices}, 1992(2):27--38, 1992.

\bibitem{Sal95}
Laurent Saloff-Coste.
\newblock Parabolic {H}arnack inequality for divergence-form second-order
  differential operators.
\newblock {\em Potential Anal.}, 4(4):429--467, 1995.
\newblock Potential theory and degenerate partial differential operators
  (Parma).

\bibitem{Sal02}
Laurent Saloff-Coste.
\newblock {\em Aspects of {S}obolev-type inequalities}, volume 289 of {\em
  London Mathematical Society Lecture Note Series}.
\newblock Cambridge University Press, Cambridge, 2002.

\bibitem{Sas26}
K{\^o}hei {Sasaya}.
\newblock {Construction of $p$-energy measures associated with strongly local
  $p$-energy forms}.
\newblock {\em J. Funct. Anal.}, 291(1):Paper No. 111489, 2026.

\bibitem{Shi24}
Ryosuke Shimizu.
\newblock Construction of {$p$}-energy and associated energy measures on
  {S}ierpi\'{n}ski carpets.
\newblock {\em Trans. Amer. Math. Soc.}, 377(2):951--1032, 2024.

\bibitem{TW02AJM}
Neil~S. Trudinger and Xu-Jia Wang.
\newblock On the weak continuity of elliptic operators and applications to
  potential theory.
\newblock {\em Amer. J. Math.}, 124(2):369--410, 2002.

\bibitem{Yan25c}
Meng {Yang}.
\newblock {On singularity of $p$-energy measures on metric measure spaces}.
\newblock {\em arXiv e-prints}, page arXiv:2505.12468v2, May 2025.

\bibitem{Yan25e}
Meng {Yang}.
\newblock {On the dichotomy of $p$-walk dimensions on metric measure spaces}.
\newblock {\em arXiv e-prints}, page arXiv:2509.08641v2, September 2025.

\bibitem{Yan25a}
Meng {Yang}.
\newblock {$p$-Poincar{\'e} inequalities and cutoff Sobolev inequalities on
  metric measure spaces}.
\newblock {\em arXiv e-prints}, page arXiv:2504.09503v4, April 2025.

\bibitem{Yos95}
K{\=o}saku Yosida.
\newblock {\em Functional analysis}.
\newblock Classics in Mathematics. Springer-Verlag, Berlin, 1995.
\newblock Reprint of the sixth (1980) edition.

\end{thebibliography}

\end{document}